\documentclass[11pt]{article}
\usepackage{amsmath,amsthm,amssymb}
\usepackage{amsfonts}
\usepackage{graphics}
\newtheorem{intro}{Theorem}[section]

\newtheorem{cintro}[intro]{Corollary}

\newtheorem{theorem}{Theorem}[subsection]
\newtheorem{lemma}[theorem]{Lemma}
\newtheorem{cor}[theorem]{Corollary}
\newtheorem{proposition}[theorem]{Proposition}
\theoremstyle{definition}
\newtheorem{definition}[theorem]{Definition}
\newtheorem{definitions}[theorem]{Definitions}
\newtheorem{remarks}[theorem]{Remarks}
\newtheorem{remark}[theorem]{Remark}
\newtheorem{ques}[intro]{Question}
\def\pri{{\mathcal{P}}{\mathbb{Z}}}%% primitive
\def\slnr{SL(n,{\mathbb{R}})}   %% SL(n,R)
\def\slsr{SL(s,{\mathbb{R}})}   %% SL(s,R)
\def\slnpr{SL(n+1,{\mathbb{R}})}   %% SL(n+1,R)
\def\slsz{SL(s,{\mathbb{Z}})}   %% SL(s,Z)
\def\slnpz{SL(n+1,{\mathbb{Z}})}   %% SL(n+1,Z)
\def\unu{{\mathbf{1}}} %% fct caracteristica
\def\calp{\mathcal{P}}   %% P caligrafic
\def\calr{\mathcal{R}}   %% R caligrafic
\def\calv{\mathcal{V}}   %% V caligrafic
\def\ff{\mathcal{F}}   %% F caligrafic
\def\dd{\mathcal{D}}   %% D caligrafic
\def\oo{\mathcal{O}}   %% O caligrafic
\def\pp{\mathcal{P}}   %% P caligrafic
\def\mm{\mathcal{M}}   %% M caligrafic
\def\hh{\mathcal{H}}   %% H caligrafic
\def\nn{\mathcal{N}}   %% N caligrafic
\def\cc{\mathcal{C}}   %% C caligrafic
\def\bc{\mathbf{C}}   %% C ingrosat
\def\cl{\mathcal{L}}   %% L caligrafic
\def\tl{\widetilde{\mathcal{L}}}   %% L caligrafic cu tilde
\def\qq{\mathfrak{Q}} %% Q gotic
\def\ss{\widetilde{\mathcal{S}}} %% S caligrafic si tilde
\def\ls{\overline{\mathcal{S}}} %% S caligrafic si linie
\def\ux{\mathbf{u}_{\bar{x}}} %% unipotent indice x
\def\1n{\frac{1}{n}} %% fractia 1/n

\newcommand {\iv}{^{-1}}
\newcommand {\dist}{\mathrm{dist}}%% distance
\newcommand {\pr}{\mathrm{proj}}%%whole projection
\newcommand {\slo}{\mathrm{sl}}%%slope projection
\newcommand {\N}{\mathbb{N}} %% positive integers
\newcommand {\Z}{\mathbb{Z}}            %% integers
\newcommand {\R}{\mathbb{R}} %% reals
\newcommand {\Q}{\mathbb{Q}} %% rationals
\newcommand {\C}{\mathbb{C}} %% complex
\newcommand {\K}{\mathbb{K}} %% generic
\newcommand {\proj}{\mathbb{P}} %% projective
\newcommand {\hip}{\mathbb{H}} %% hyperbolic plane
\newcommand {\sph}{\mathbb{S}} %% sphere
\newcommand {\q}{\frak q} %% sphere
\newcommand {\me}{\medskip}
\newcommand {\sm}{\smallskip}
\newcommand {\bi}{\bigskip}

\newcommand {\Notat}{\noindent {\it{Notation}}:\thickspace } %% notatie
\begin{document}
\makeatletter
\title{ Diophantine approximation on rational quadrics}
\author{
Cornelia DRU\c{T}U\\ \\
UFR de Math\'ematiques et UMR CNRS no. 8524,\\
Universit\'e de Lille I,
F--59655 Villeneuve d'Ascq, France \\
 fax no. (33)3.20.43.43.02\\
Cornelia.Drutu@math.univ-lille1.fr}
\date{ }
\maketitle
\footnotetext[1]{
 {\it $2000$ Mathematics Subject Classification}: 11J83, 22E40, 53C35.}
\makeatother

\begin{abstract}
\noindent We compute the Hausdorff dimension of sets of very well
approximable vectors on rational quadrics. We use ubiquitous
systems and the geometry of locally symmetric spaces. As a
byproduct we obtain the Hausdorff dimension of the set of rays
with a fixed maximal singular direction, which move away into one
end of a locally symmetric space at linear depth, infinitely many
times.
\end{abstract}

\tableofcontents

\section{Introduction}

The main result in the present paper is the computation of the
Hausdorff dimensions of the sets of very well approximable vectors
on a rational quadric $\qq$. The method is to consider the
rational non-degenerate quadratic form $\q :\R^n \to \R $ such
that the quadric $\qq$ is defined by $\q=1$ and the quadratic form
$L_\q : \R^{n+1} \to \R \, ,\, L_\q(x_1,\dots ,
x_{n+1})=x_{n+1}^2-\q (x_1,\dots , x_{n})$. The connected
component of the identity $SO_{I} (L_\q)$ of the stabilizer
$SO(L_\q)$ of the form $L_\q$ is a semisimple group (simple if
$n\neq 3$). The integer points of this group compose a lattice.
One can consider the symmetric space associated to $SO_{I} (L_\q)$
and its quotient by the lattice, which is a locally symmetric
space. The set of very well approximable vectors on $\qq$ can be
defined in terms of the geometry of the locally symmetric space,
and its Hausdorff dimension can be estimated using an ubiquitous
system which appears in this context and the general properties of
ubiquitous systems.

\subsection{Hausdorff dimension of sets of very well
approximable vectors in $\R^n$}

We denote by $\|\cdot \|_e$ the Euclidean norm and by $\|\cdot \|$
the max-norm in $\R^n$,
$$\|x \| = \max \{|x_1|, |x_2|, \ldots , |x_n| \}\, .$$

Throughout $\psi :\R_+ \to \R_+$ denotes a decreasing function
satisfying $\lim_{x\to \infty }\psi(x)=0$, also called \textit{an
approximating function}. Rational vectors are always written in the form $\frac{1}{q}\bar{p}$, where $q\in \N,\,
\bar{p}=(p_1,\dots ,p_n)\in \Z^n$ and $\gcd (q,p_1,\dots ,p_n)=1$.

Let $M$ be a submanifold of $\R^n$. The set of {\it simultaneously
$\psi$-approximable vectors in} $M$ is defined by
$$ \mathcal{S}_\psi (M)=\left\{\bar{x}\in M \: ;\:
\|q\bar{x}-\bar{p}\|\leq \psi (q) \mbox{ for infinitely many }q\in
\N \, ,\, \bar{p}\in \Z^n\right\}\, .
$$

In the particular case when $\psi (x)=\frac{1}{x^{\alpha }}$ with
$\alpha
> \frac{1}{n}$, the set is also denoted by $\mathcal{S}_\alpha (M)$ and it is called the set
 of {\it simultaneously $\alpha$-very well
approximable vectors in} $M$. A subset of it is the set of {\it
simultaneously exactly-$\alpha$-very well approximable vectors in}
$M$,
$$
\mathcal{ES}_\alpha(M)=\{\bar{x}\in M \: ;\: \bar{x}\in
\mathcal{S}_\alpha(M) \mbox{ and }\bar{x}\not\in \mathcal{S}_\beta
(M)\, ,\, \forall \beta >\alpha \}\, .
$$

Likewise is defined the set of {\it linearly $\psi $-approximable
vectors in $M$}, by $$ \mathcal{L}_\psi(M)=\left\{\bar{x}\in M \:
;\: |\bar{q}\cdot \bar{x}-p|\leq \psi(\| \bar{q }\|) \mbox{ for
infinitely many }\bar{q} \in \Z^n \, ,\, p\in \N \right\}\, , $$
where $\bar{q}\cdot \bar{x}=\sum_{i=1}^n q_ix_i$. In particular
when $\psi (x) =x^{-\beta }$, with $\beta > n$, the previous set
is denoted by $\mathcal{L}_\beta(M)$ and it is called the set
 of {\it linearly $\beta$-very well
approximable vectors in} $M$.

Khintchine's transference principle \cite[$\S 1.3.1$]{BD} implies
that
$$
\bigcup_{\alpha
> 1/n } \mathcal{S}_\alpha(\R^n )=\bigcup_{\beta > n}
\mathcal{L}_\beta(\R^n)\: .
$$

The vectors in this set are called \textit{very well approximable
vectors}. A consequence of the Khintchine-Groshev Theorem
\cite[$\S 1.3.4$]{BD} is that the set of very well approximable
vectors in $\R^n$ is of Lebesgue measure zero. Thus, in order to
study the sets of type $\mathcal{S}_\psi$ and $\mathcal{L}_\psi $
when $\psi $ decreases sufficiently quickly at infinity so that they are of
Lebesgue measure zero, a more appropriate
tool is the Hausdorff dimension and the Hausdorff measure. In the
sequel we denote by $\dim_H A$ the Hausdorff dimension of a subset
$A$ in a metric space. We denote by $\hh^s$ the Hausdorff measure
corresponding to the parameter $s$ (see Section \ref{ghm} for
definitions). It has been proved in \cite{Ja} that given $s\in
[0,n)$,
\begin{equation}\label{jar}
    \hh^s (\mathcal{S}_\psi(\R^n ))=\left\{%
\begin{array}{ll}
    0\, , & \mbox{ if } \sum_{k=1}^\infty k^{n-s}\psi (k)^s < \infty \, , \\
    \infty \, , & \hbox{ if } \sum_{k=1}^\infty k^{n-s}\psi (k)^s = \infty \, . \\
\end{array}%
\right.
\end{equation}

In particular, for any $\alpha >\frac{1}{n}$

\begin{equation}\label{jard}
    d=\dim_H \mathcal{S}_\alpha(\R^n )=\frac{n+1}{\alpha +1}\; \;
\mbox{ and }\; \; \hh^d (\mathcal{S}_\alpha(\R^n ))=\infty  \, .
\end{equation}

This implies that both relations also hold for
$\mathcal{ES}_\alpha(\R^n )$ instead of $\mathcal{S}_\alpha(\R^n
)$.

In \cite{BoD} it was shown that
\begin{equation}\label{linr}
   \dim_H \mathcal{L}_\beta (\R^n )=n-1+\frac{n+1}{\beta +1},\; \;
\forall \beta >n.
\end{equation}

Moreover the following holds \cite{DV}. Let $s\in (n-1, n)$ and
let $\psi $ be an approximating function. Then
\begin{equation}\label{linrh}
    \hh^s (\mathcal{L}_\psi(\R^n ))=\left\{%
\begin{array}{ll}
    0\, , & \mbox{ if } \sum_{k=1}^\infty k^{2n-1-s}\psi (k)^{s-(n-1)} < \infty \, ,\\
    \infty \, , & \hbox{ if } \sum_{k=1}^\infty k^{2n-1-s}\psi (k)^{s-(n-1)} = \infty \, . \\
\end{array}%
\right.
\end{equation}

\subsection{Known results on very well approximable vectors on manifolds}

 The general question to ask is under what
 conditions the vectors in a submanifold $M$ of $\R^n$ behave similarly to
  the vectors in $\R^n$, with respect to Diophantine approximation. If $M$ is a rational affine subspace of dimension $k<n$ in $\R^n$ then $M=\mathcal{S}_{\frac{1}{k}}(M)$. Therefore, rational affine subspaces
must be avoided.

M.M. Dodson, B.P. Rynne and J.A.G. Vickers have shown in
\cite{DRV3} and in \cite{DRV4} that under some non-zero curvature
condition, the set of very well approximable vectors in $M$ is of
measure $0$. D. Kleinbock and G.A. Margulis have shown in
\cite{KM2} the same result in a submanifold $M$ of $\R^n $
non-degenerate almost everywhere (they have actually shown that a
larger set, the set of very well multiplicatively approximable
vectors in $M$, has measure $0$ in this case). A point $\bar{x}\in
M$ is \textit{non-degenerate} if in a neighborhood of $\bar{x}$,
$M$ is not near to any affine subspace. More precisely, in a
neighborhood of $\bar{x}$ the submanifold $M$ is parameterized by
a function $f$ which is $l$ times continuously differentiable and
such that its partial derivatives in $\bar{x}$ up to order $l$
span $\R^n$. A submanifold $M$ \textit{non-degenerate almost
everywhere} is a submanifold in which almost every point is
non-degenerate.

A Khintchine-Groshev type theory was equally developed in the
setting of manifolds. Concerning the linear approximation (that
is, the Groshev type theory) it has been shown that any
submanifold non-degenerate almost everywhere is of Groshev type
(see \cite{Ber}, \cite{BKM}, \cite{BBKM}, \cite{BDV2}). For the
known results in simultaneous approximation, that is for a
Khintchine type theory on manifolds, we refer to \cite{BD},
\cite{BDV1}, \cite{BDV2}, \cite{DRV3}, \cite{DRV4} and references
therein.

Consider a submanifold $M$ with the set of very well approximable
vectors of measure zero. Such a submanifold is also called
\textit{extremal}. One can ask what is the Hausdorff dimension of
each set $\mathcal{S}_\alpha(M)$ with $\alpha > \frac{1}{n}$ and
$\mathcal{L}_\beta (M)$ with $\beta
>n$. More is known about the sets $\mathcal{L}_\beta(M)$. R.C. Baker
\cite{Ba} proved that if $M$ is a planar curve of class $C^3$
whose curvature is zero at most in a set of points of Hausdorff
dimension zero, then
\begin{equation}\label{bak}
\dim_H \mathcal{L}_\beta(M)=\frac{3}{1+\beta } \; ,\forall \beta
\geq 2 \, .
\end{equation}

M.M. Dodson, B.P. Rynne and J.A.G. Vickers \cite{DRV1} later
proved that if $M$ is a $C^3$-submanifold of dimension $m\geq 2$
in $\R^n$ such that at least two principal curvatures are not zero
except on a set of Hausdorff dimension at most $m-1$, then
\begin{equation}\label{drv}
    \dim_H \mathcal{L}_\beta(M)=m-1+\frac{n+1}{1+\beta }\; ,\forall
\beta \geq n \, .
\end{equation}

H. Dickinson and M.M. Dodson have shown in \cite{DD2} that if $M$
is extremal then
$$
\dim_H \mathcal{L}_\beta(M)\geq m-1+\frac{n+1}{1+\beta }\;
,\forall \beta \geq n .
$$

Finally, in \cite{BDV2}, V. Beresnevich, D. Dickinson and S.
Velani have shown that, given $M$ an $m$-dimensional submanifold
in $\R^n$, with $n\geq 2$, $M$ non-degenerate almost everywhere,
the following holds. Consider $s\in (m-1,m)$. If
$$
\sum_{k=1}^\infty \psi(k)^{s-(m-1)}k^{n+m-1-s} = \infty \mbox{
then }\hh^s (\mathcal{L}_\psi (M))=\infty \, .
$$

In particular for $\psi (x)=x^{-\beta }$ this implies the result
of H. Dickinson and M.M. Dodson, under the given hypotheses for
$M$, and moreover it shows that for $d=m-1+\frac{n+1}{1+\beta }$,
the Hausdorff measure $\hh^d \left( \mathcal{L}_\beta (M)\right)$
is $\infty$.

These results and Khintchine's transference principle can be used
to obtain upper and lower bounds for the Hausdorff dimensions of
the sets $\mathcal{S}_\alpha(M)$. As far as the exact Hausdorff
dimension for sets $\mathcal{S}_\alpha $ goes, the known results
are the following. In \cite{BDV2} it is shown that, given $\psi $
an approximating function such that $\lim_{x\to \infty }x \psi
(x)=0 $ and $s\in (0,1)$, the following holds.
\begin{displaymath}
\hh^s (\mathcal{S}_\psi(\sph^1 ))=\left\{%
\begin{array}{ll}
    0\, , & \mbox{ if } \sum_{k=1}^\infty \left( \frac{\psi (k)}{k} \right)^s < \infty \, , \\
    \infty \, , & \hbox{ if } \sum_{k=1}^\infty \left( \frac{\psi (k)}{k} \right)^s = \infty \, . \\
\end{array}%
\right.
\end{displaymath}

In particular this implies that
\begin{equation}\label{s1}
\dim_H \mathcal{S}_\alpha(\sph^{1})=\frac{1}{1+\alpha}\; \mbox{
and }\;   \hh^{1/(1+\alpha )}(\mathcal{S}_\alpha(\sph^1 ))=\infty
\, ,\, \forall \alpha > 1 .
\end{equation}

The first equality in (\ref{s1}) had already been proved in
\cite{DD}.

In \cite[Theorem 4.8]{BD} it is proved that if $k\in \N \, ,\,
k\geq 3$, and
$$
\mathcal{C}_k=\{ (x,y)\in \R^2 \; ;\;  x^k+y^k=1 \}\, ,
$$ then $\mathcal{S}_\alpha(\mathcal{C}_k)$ contains at most four
points for $\alpha >k-1$, hence $\dim_H
\mathcal{S}_\alpha(\mathcal{C}_k) =0$ for $\alpha >k-1$.

The examples of $\sph^1$ and $\mathcal{C}_k\, ,\, k\geq 3$,
already emphasize that, unlike in the case of linear
approximation, a condition of non-zero curvature is not enough to
deduce the Hausdorff dimensions of the sets $\mathcal{S}_\alpha$.
B.P. Rynne \cite{Ry2} showed moreover that for every
$C^k$-submanifold $M$ of $\R^n$ of dimension $m$ there exist
$C^k$-submanifolds $M_z$ and $M_p$ arbitrarily $C^k$-close to $M$
(in a suitable sense) such that for $\alpha$ sufficiently large
$\mathcal{S}_\alpha (M_z)=\emptyset $ and $\dim_H
\mathcal{S}_\alpha (M_p)
>\frac{m+1}{k(\alpha +1)}$. It follows that conditions taking into
account only the structure of differential submanifold and
depending continuously on this structure cannot suffice to obtain
information about $\mathcal{S}_\alpha (M)$, at least not for large
values of $\alpha$. The following result from \cite{BDV1} on the
other hand seems to indicate that for values of $\alpha $ near to
$\frac{1}{n}$, where $n$ is the dimension of the ambient space
$\R^n$, there should exists however a formula holding for any
non-degenerate submanifold of $\R^n$. More precisely, in
\cite{BDV1} it is shown the following. Let $f\in C^3 ([a,b])$,
$a<b$, let $\mathcal{C}_f=\{ (t,f(t))\; ;\;  t\in [a,b]\}$, let
$s\in (1/2 ,1)$ and let $\psi$ be an approximating function.
\begin{itemize}
    \item If $\sum_{k=1}^\infty k^{1-s} \psi (k)^{s+1}=\infty$
    then $\hh^s(\mathcal{S}_\psi (\cc_f))=\infty$;
    \item Let $\lambda_\psi =\liminf_{x\to \infty }\frac{-\ln \psi (x)}{\ln
    x}\, $. If the Hausdorff dimension of the set $\{ t\in [a,b]\; ;\;
    f''(t)=0\}$ is at most $\frac{2-\lambda_\psi }{1+\lambda_\psi
    }$ then $d=\dim_H (\mathcal{S}_\psi (\mathcal{C}_f))=\frac{2-\lambda_\psi }{1+\lambda_\psi
    }$. Assume moreover that $\lambda_\psi \in (1/2 ,1)$. Then
    $\limsup_{x\to \infty }x^{2-d}\psi (x)^{d+1}>0$ implies that $\hh^d(\mathcal{S}_\psi
    (\cc_f))=\infty$.
\end{itemize}

In the particular case when $\psi(x)=x^{-\alpha }$ with $\alpha
\in (1/2 ,1)$ this gives the following.
\begin{itemize}
    \item  $\dim_H
\mathcal{S}_\alpha (\mathcal{C}_f)\geq d=\frac{2-\alpha }{1+\alpha
    }$ and $\hh^d(\mathcal{S}_\alpha
    (\cc_f))=\infty$;
    \item  If moreover the Hausdorff dimension of the set $\{
t\in [a,b]\; ;\;
    f''(t)=0\}$ is at most $\frac{2-\alpha }{1+\alpha
    }$ then $\dim_H \mathcal{S}_\alpha (\mathcal{C}_f)=d$.
\end{itemize}

In the particular case of a rational quadric ${\mathfrak{Q}}$ in $\R^2$ one obtains $\dim_H
\mathcal{S}_\alpha ({\mathfrak{Q}})=\frac{2-\alpha }{1+\alpha
    }$ for $\alpha
\in [1/2 ,1)$. Note that for ${\mathfrak{Q}}=\sph^1$ this differs from the formula for $\alpha
>1$ given in (\ref{s1}). Thus in this case, unlike
in the cases treated in (\ref{jard}), (\ref{linr}), (\ref{bak})
and (\ref{drv}), the Hausdorff dimension of the sets of very well
approximable vectors is not a rational function in $\alpha$ but a
\textit{piecewise} rational function in $\alpha$, with different
expressions for $\alpha \in [1/2 ,1)$ and for $\alpha
>1$.

In \cite{DL} the Hausdorff dimension of $\mathcal{S}_\alpha (M)$
has been computed for large values of $\alpha $ and for $M$ a
manifold parameterized by polynomials with integer coefficients.

\subsection{Very well approximable vectors on rational
quadrics}\label{main}

The purpose of the present paper is to compute the Hausdorff
dimension of the sets $\mathcal{S}_\alpha({\mathfrak{Q}}_\q )$ for
$\alpha >1$, where ${\mathfrak{Q}}_\q $ is a quadric defined by
the equation $\q =1$, for a given non-degenerate rational
quadratic form $\q :\R^n \to \R $. Obviously $\frak q$ cannot be
negative definite. The main result of the paper, formulated not in
the most general form, is the following.

\newpage

\begin{intro}\label{T2} Let $\psi $ be an approximating function such that $\lim_{x\to
\infty } x \psi (x)=0$.
\begin{itemize}
  \item[\textbf{(1)}] If ${\mathfrak{Q}}_\q \cap \Q^n = \emptyset$ then
  $\mathcal{S}_\psi({\mathfrak{Q}}_\q ) = \emptyset $.

 \item[\textbf{(2)}] If ${\mathfrak{Q}}_\q \cap \Q^n \neq
 \emptyset$ then
    $$
\dim_H \mathcal{S}_\psi (\qq_\q ) = \sigma (n-1) \, ,
    $$ where $\sigma = \limsup_{x\to \infty } \frac{\ln x}{\ln x-\ln \psi
    \left( x\right)}\, $.

    Moreover, if $\limsup_{x\to \infty } x^{1-\sigma }\psi
    \left( x\right)^\sigma >0 $ then $\hh^{\sigma (n-1)} \left(\mathcal{S}_\psi (\qq_\q )\right)=\infty
    $.

     In particular the set $\mathcal{S}_\alpha (\qq_\q )$ has Hausdorff
dimension $d=\frac{n-1}{1+\alpha}$ for any $\alpha >1$ and $\hh^d
(\mathcal{S}_\alpha (\qq_\q ))=\infty$. Both statements also hold
for the set  $\mathcal{E}\mathcal{S}_\alpha (\qq_\q )$.
\end{itemize}
\end{intro}

\sm

 According to
\cite[Chapter 1, $\S 7$]{BSh} a rational non-degenerate quadratic
form in $n\geq 5$ variables takes the zero value on $\Z^n
\setminus \{ (0,\dots ,0 ) \}$ if and only if it is not defined.
This theorem applied to the form $L_\q$ implies that for $n\geq
4$, ${\mathfrak{Q}}_\q \cap \Q^n \neq \emptyset$ for any rational
quadratic form $\q \, $. For $n=2,3$ see \cite[Chapter 1, $\S
7$]{BSh}.

\me

\noindent \textit{Outline of proof of Theorem \ref{T2}.}

Statement (1) is a straightforward consequence of Lemma
\ref{aprox}. Therefore we may assume that ${\mathfrak{Q}}_\q \cap
\Q^n \neq \emptyset$. The symmetric space corresponding to the
semisimple group $SO_I (L_\q)$ is the space $\mathcal{P}_{n+1}
(L_\q )$ of minimal positive definite quadratic forms $Q$ such
that $|L_\q (\bar{x})|\leq Q(\bar{x})\, ,\, \forall \bar{x}$ (see
\cite{Bo} or Section \ref{parab}). The boundary at infinity of it,
$\partial_\infty \mathcal{P}_{n+1} (L_\q )$, is a spherical
building which can be canonically identified with the spherical
building of flags of $\R^{n+1}$ composed of subspaces totally
isotropic with respect to $L_\q$ (\cite[$\S 15, \S 16$]{Mo},
\cite[$\S 4.G$]{Wi}). In particular $\partial_\infty
\mathcal{P}_{n+1} (L_\q )$ contains a maximal singular stratum
corresponding to the 1-dimensional subspaces totally isotropic
with respect to $L_\q$. We call it \textit{the stratum} $\wp$ and
the points composing it \textit{points of type }$\wp $.
Correspondingly we say that a geodesic ray in $\mathcal{P}_{n+1}
(L_\q )$ is \textit{of type }$\wp $ if its point at infinity is.

\me

\noindent \textit{Convention}: Throughout, a semisimple group acts
by isometries \textit{on the right} on the symmetric space
associated to it and on its boundary at infinity.

\me

The quadric ${\mathfrak{Q}}_\q $ can be identified with an open
Zariski dense subset of the stratum $\wp$ in $\partial_\infty
\mathcal{P}_{n+1} (L_\q )$. On the other hand, for any geodesic
ray $\varrho$ in $\mathcal{P}_{n+1} (L_\q )$ of type $\wp$, the
opposite unipotent $U_+(\varrho )$ of $\varrho $ (see Section
\ref{bdry0} for a definition) can be identified with an open Zariski
dense subset of the stratum $\wp$ via the bijection $\mathbf{u}
\mapsto \varrho(\infty )\mathbf{u}$. With a countable covering
argument we can replace in our study $\mathcal{S}_\psi (\qq_\q)$
by $\mathcal{S}_\psi (\Omega )$, where $\Omega $ is a relatively
compact open subset whose closure is contained in the image of
$U_+(\varrho)$ for some $\varrho$. The set $\Omega $ can be
identified with a relatively compact open subset of $U_+(\varrho )$.

Let $\Gamma = SO_I (L_\q )\cap \slnpz$. The locally symmetric
space $\mathcal{V} = \mathcal{P}_{n+1} (L_\q )/ \Gamma $ has ends
if and only if ${\mathfrak{Q}}_\q \cap \Q^n \neq \emptyset$.
Moreover there exist finitely many geodesic rays $\bar{r}_i,\,
i\in \{ 1,2,\dots , k \}$, in $\calv$ such that their lifts $r_i$
in $\mathcal{P}_{n+1} (L_\q )$ are of type $\wp $ and moreover the
following holds. Let $r_i\Gamma $ be the $\Gamma$-orbit of $r_i$,
let $Csp_i$ be the corresponding set of points at infinity
$r_i(\infty )\Gamma $ and let $Csp$ be the set of points at
infinity $\bigcup_{i=1}^k Csp_i$. Then $Csp$ intersected with
${\mathfrak{Q}}_\q$ seen as a subset of $\partial_\infty
\mathcal{P}_{n+1} (L_\q )$ is ${\mathfrak{Q}}_\q \cap \Q^n$. In
particular $\Omega \cap \Q^n = \Omega \cap Csp \, $, and it can be
seen as a subset of $U_+(\varrho )$. For each $w\in Csp\cap \Omega
$ we use $\mathbf{u}_w$ to denote its corresponding unipotent
element in $U_+ (\varrho )$.

Note that to every point $w=r_i(\infty )\gamma$ in $Csp$ it is
naturally associated a horoball $Hb_w=Hb(r_i \gamma)$ having it as
a basepoint (see Section \ref{not} for the definition of a
horoball). Let $\varrho^{op}$ be the geodesic ray opposite to
$\varrho \, $. To every element $w\in Csp\cap \Omega $ one can
associate a weight $d_w\in \R_+$ which is the distance from the
horoball $Hb(\varrho^{op})$ to the horoball $Hb_w$. In $\Omega $
seen as a subset of $U_+(\varrho )$ one can then consider the set
$\ss_\Psi^0 (\Omega )$ of elements $\mathbf{u}$ such that
\begin{equation}\label{inegspsi}
    \dist (\mathbf{u}\, ,\, \mathbf{u}_w )\leq \Psi (d_w)\, ,\mbox{
for infinitely many }w\in Csp\, ,
\end{equation} where $\dist $ is a left invariant metric on $U_+(\varrho )$
 and $\Psi $ is an approximating function.

It turns out that, due to Lemma \ref{aprox}, the sets
$\mathcal{S}_\psi (\Omega )$ and $\ss_\Psi^0 (\Omega )$ are
closely related, for an appropriate choice of the function $\Psi$.
This relation is established using some explicit formulas obtained
in Sections \ref{parab} and \ref{gho}. See the double inclusion
(\ref{dincls}) and the whole discussion in Section \ref{trei} for
details.

It suffices to study the set $\ss_\Psi^0 (\Omega )$ from the point
of view of the Hausdorff dimension. Moreover, it is not difficult
to see that one can restrict the study to a subset $\ss_\Psi^i
(\Omega )$ defined by replacing in (\ref{inegspsi}) the set $Csp$
by the subset $Csp_i$.

In the particular cases when $\q$ is positive definite or of
signature $(1,n-1)$, $\mathcal{P}_{n+1}(L_\q)$ is isometric to the
hyperbolic space $\hip^n$, and all the results in this paper
follow from the results in \cite[$\S 8.3$]{BDV2}, generalizing
previous results from \cite{HV}. We give an argument for the
remaining cases. This argument actually works for the two previous
cases too, with some slight modification. The inequality $ \dim_H
\ss_\Psi^i (\Omega ) \leq  \sigma (n-1) $ is not difficult to
obtain. The main ingredient in its proof is the counting result
Corollary \ref{lSu}, which gives an estimate of the number of
balls $B(\mathbf{u}_w \, ,\, \Psi (d_w ) )$ in $U_+(\varrho)$ of a
given size. This counting result itself follows from the
equidistribution result Proposition \ref{equiduc}.

For the converse inequality we use ubiquitous systems. We deduce
from the equidistribution result Proposition \ref{equidunip} and
the counting result Corollary \ref{lSu} that the set of points
$\Re = \{ \mathbf{u}_w\; ;\; w\in Csp_i\cap \Omega \}$ together
with the weight function $\varpi : Csp_i\cap \Omega \to \R_+$,
$\varpi (w)=d_w$, compose a local ubiquitous system with respect
to an appropriate ubiquitous function and an appropriate
increasing sequence of positive numbers, in the terminology of
Section \ref{us}. We then use the properties of ubiquitous systems
as developed in \cite{BDV2} to deduce the lower bound of the
Hausdorff dimension of $\ss_\Psi^i (\Omega )$, as well as the
other results.

Some comments are necessary concerning the counting result
Corollary \ref{lSu}. This statement corresponds in our case to the
result in \cite[$\S 6$, Proposition 4]{Su}, given for the rank one
case. A generalization of Sullivan's result in the setting of
geometrically finite Kleinian groups has been given in \cite{HV}.
 A consequence of Corollary \ref{lSu} is the following statement.

\begin{cintro}[equidistribution of rational vectors on rational quadrics]\label{counting3}

Suppose that $\qq_\q \cap \Q \neq \emptyset$. Let $\Omega$ be a
relatively compact open subset of ${\mathfrak{Q}}_{\q}$ such that
its closure does not intersect $T_{\bar{x}_0}{\mathfrak{Q}}_{\q}$
for some $\bar{x}_0\in {\mathfrak{Q}}_{\q}$. Let $a>1$. For every
open subset $\oo$ of $\Omega$ we denote by $N(k\, ;\, \oo )$  the
cardinal of the set of rational vectors
$$\left\{ \frac{1}{q}\bar{p} \in \Q^n \cap \oo \; ;\; | q |
\in [a^k,a^{k+1})\right\}\, .$$

For any $a\geq a_0(\q\, ,\, \Omega)$ we have that
$$
\mathbf{K}_1\, a^{(k+1)(n-1)}\nu (\oo )\; \leq\;  N(k\, ;\, \oo
)\; \leq \; \mathbf{K}_2\, a^{(k+1)(n-1)}\nu (\oo )\, ,\; \mbox{
for every }\; k\geq k_0(\oo ,\Omega )\, ,
$$ where $\nu $ is the canonical
measure on $\qq_\q$ and $\mathbf{K}_i=\mathbf{K}_i(\q \, ,\,
\Omega )$.
\end{cintro}

It is worth mentioning that our methods cannot be
used to obtain either Khintchine type results or results on badly
approximable vectors in ${\mathfrak{Q}}_\q$ or any other type of
results concerning vectors approximable nearly as well as the
generic vectors in $\R^n$. This can be seen for instance by
applying the logarithm law (\cite{Su}, \cite{KM}) in our setting.
It implies that for every $\varepsilon
>0$, for almost every $\bar{x}\in {\mathfrak{Q}}_\q $, we have
$$
\left\|\bar{x} - \frac{1}{q}\bar{p} \right\| \geq
\frac{c_1(\bar{x})}{q(\ln q )^{\frac{1}{n-1}+\varepsilon }} >
\frac{1}{q^{1+\frac{1}{n}}}, \, \forall\: \frac{1}{q}\bar{p}\in
{\mathfrak{Q}}_\q \, .
$$

Consequently, for almost all $\bar{x}\in {\mathfrak{Q}}_\q $ the
rational approximants are outside ${\mathfrak{Q}}_\q $. It seems
that in order to study badly approximable and Khintchine type
approximable vectors in ${\mathfrak{Q}}_\q $, the study of the
intrinsic geometry of $\mathcal{V}$ is not sufficient, and one has
to consider also the ``ambient'' geometry of
$\mathcal{T}_{n+1}=\mathcal{P}_{n+1} / SL(n+1,{\mathbb{Z}})$,
where $\mathcal{P}_{n+1}$ is the symmetric space of positive
definite quadratic forms on $\R^{n+1}$ of determinant 1 in the
canonical basis. The locally symmetric space $\mathcal{T}_{n+1}$
contains an embedding of $\mathcal{V}$ \cite[$\S 5$]{Bo}.

\subsection{Rays moving away in the cusp}

We consider again the set $\ss_\Psi^i (\Omega )$ defined in
Section \ref{main}. Without loss of generality we may assume that
$\Omega = U_+ (\varrho )$ and we may denote the corresponding set
simply by $\ss_\Psi^i$. This set can be related to a set of
geodesic rays moving away in the cusp for infinitely many times
$t$ at depth at least $t-\phi (t)$, where the depth is measured
with respect to the ray $\bar{r}_i$ and $\phi :[a, +\infty ) \to
[b,+\infty )$ is a function depending on the function $\Psi$. The
results on the Hausdorff dimension of the set $\ss_\Psi^i $ can be
thereby translated in terms of this set of rays. To simplify the
exposition we present here a particular case of the results that
can be obtained with such an argument, the general statements can
be found in Section \ref{sectrays1}.

Let $\beta \in (0,1)$ and define
\begin{displaymath}
\mathcal{R}_{\beta } = \left\{ \mathbf{u}\in U_+(\varrho ) \; ;\;
f_{\bar{r}_i}\left(\pr \left(\varrho (t)\mathbf{u}\right)\right)
\leq -\beta t \mbox{ infinitely many times as }t\to \infty
\right\}\, .
\end{displaymath} Above we have resumed the notation in Section
\ref{main}, $\pr $ denotes the projection of $\calp_{n+1} (L_\q )$
onto $\calv$ and $f_{\bar{r}_i}$ denotes the Busemann function of
the ray $\bar{r}_i$ in $\calv$, as defined in Section \ref{not}.

Note that $\mathcal{R}_{\beta }$ can also be seen as a set of
geodesic rays, by identifying each $\mathbf{u}$ to the ray
$\varrho \, \mathbf{u}$. The condition defining $\mathcal{R}_\beta
$ means that for infinitely many times $t$ the projection onto
$\calv$ of the geodesic ray $\varrho\, \mathbf{u}$ goes into the
cusp at depth at least $\beta t$, the depth into the cusp being
measured with respect to the ray $\bar{r}_i$. We also consider a
subset of $\mathcal{R}_{\beta }$, representing the rays which in
some sense do not go deeper than $\beta t$ in the cusp as $t\to
\infty$:
$$
\mathcal{E}\mathcal{R}_{\beta } =\calr_\beta \setminus
\bigcup_{\beta' > \beta
    }\calr_{\beta'} =\left\{ \mathbf{u}\in
\mathcal{R}_\beta \; ;\; \limsup_{t\to +\infty }
\frac{-f_{\bar{r}_i}\left(\pr
\left(\varrho(t)\mathbf{u}\right)\right)}{t} =\beta \right\}\, .
$$

\begin{intro} [Corollary \ref{corrbeta}]\label{TR1}
For any $\beta \in (0,1)$,
$$
\dim_H \mathcal{R}_{\beta }=\dim_H \mathcal{E}\mathcal{R}_{\beta
}=(1-\beta )\dim U_+ (\varrho )=d\; \mbox{ and }\; \hh^d \left(
\mathcal{R}_{\beta } \right)=\hh^d \left(
\mathcal{E}\mathcal{R}_{\beta } \right)=\infty\, .$$
\end{intro}

For a discussion of the cases $\beta =0,1$ see Section
\ref{sectrays1}.

A natural question to ask is whether other results on the
Hausdorff dimension and measure of sets of very well approximable vectors have
an interpretation in terms of rays moving away in the cusp of some
locally symmetric space. We establish such an interpretation for the
formulas (\ref{jar}) and (\ref{linrh}). Most likely this can be
done in other cases as well. For the two formulas that we discuss
the appropriate symmetric space is $\calp_{n+1}$, and the
appropriate locally symmetric space is
$\mathcal{T}_{n+1}=\calp_{n+1}/\slnpz$. Let $\pr$ be the
projection of $\calp_{n+1}$ onto $\mathcal{T}_{n+1}\, $. Let $r_1$
and $r_n$ be the geodesic rays in $\calp_{n+1}$ defined as in
(\ref{rays}). The ray $r_i,\, i=1,n,$ projects onto a geodesic ray
$\bar{r}_i$ in $\mathcal{T}_{n+1}$. We define the set
\begin{displaymath}
\mathcal{R}_{\beta }^i = \left\{ \mathbf{u}\in U_+(r_i) \; ;\;
f_{\bar{r}_i}\left(\pr \left(r_i(t)\mathbf{u}\right)\right)\leq
-\beta t \mbox{ infinitely many times as }t\to \infty \right\}\, ,
\end{displaymath} where $i=1,n$, and $\beta \in (0,1)$. We also
consider the subset
$$
    \mathcal{E}\calr_\beta^i = \calr_\beta^i \setminus \bigcup_{\beta' > \beta
    }\calr_{\beta'}^i =\left\{ \mathbf{u}\in
\calr_\beta^i \; ;\; \limsup_{t\to \infty }
\frac{-f_{\bar{r}_i}\left(\pr
\left(r_i(t)\mathbf{u}\right)\right)}{t}=\beta \right\}\, .
    $$

Formula (\ref{jar}) implies the following.
\begin{intro}[Corollary \ref{corr11}]\label{TR2}
For any $\beta \in (0,1)$,
$$
\dim_H \calr_\beta^i=\dim_H \mathcal{E}\calr_\beta^i=(1-\beta
    )\dim U_+(r_i)=d\; \mbox{ and }\; \hh^d (\calr_\beta^i )=\hh^d (\mathcal{E}\calr_\beta^i )= \infty\,
    ,\, i=1,n \, .$$
\end{intro}

Formula (\ref{linrh}) also can be expressed in terms of sets of
rays moving away in the cusp, but the situation changes slightly.
In this case the ray in the cusp with respect to which the depth
is measured and the rays whose behavior is studied are not in the
same orbit of $\slnpr$, or in the terminology of Section
\ref{bdry}, they do not have the same slope. This explains why in
this case the parameter $\beta $ does not get near to $1$, but is
bounded by a smaller constant depending on the two slopes. More
precisely, we define for every $\beta \in \left( 0,\frac{1}{n}
\right)$
\begin{displaymath}
\mathcal{R}_{\beta }^{1n} = \left\{ \mathbf{u}\in U_+(r_1) \; ;\;
f_{\bar{r}_n}\left(\pr \left(r_1(t)\mathbf{u}\right)\right)\leq -
\beta t \mbox{ infinitely many times as }t\to \infty \right\}\, .
\end{displaymath}
Let $\mathcal{E}\calr_\beta^{1n} = \calr_\beta^{1n} \setminus
\bigcup_{\beta' > \beta
    }\calr_{\beta'}^{1n}\, $. The sets $\mathcal{R}_{\beta }^{n1}$ and
    $\mathcal{E}\mathcal{R}_{\beta }^{n1}$ can be defined
similarly by intertwining $1$ and $n$.

\begin{intro}[Corollary \ref{corr1n}]\label{TR3}
For any $\beta \in \left( 0,\frac{1}{n} \right)$,
$$
\dim_H \calr_\beta^{ij}=\dim_H \mathcal{E}\calr_\beta^{ij}=
(1-\beta
    )\dim U_+(r_i)=d\; \mbox{ and }\; \hh^d (\calr_\beta^{ij} )=\hh^d (\mathcal{E}\calr_\beta^{ij} )=
    \infty\, ,\, \{ i,j \} = \{ 1,n \}\, .
$$
\end{intro}

For the cases $\beta = 0\, ,\, \frac{1}{n}\, $, see Section
\ref{cds}.

\subsection{Open question}

Theorems \ref{TR1}, \ref{TR2} and \ref{TR3} suggest that there
might be a general formula for the Hausdorff dimension of the set
of rays moving away into the cusp at linear depth. This justifies
the following question.

Let $X$ be a symmetric space of non-compact type without Euclidean
factors, let $G$ be the connected semisimple group of isometries
of $X$, let $\Gamma$ be a non-uniform irreducible lattice of
isometries of $X$, let $\calv =X/\Gamma $ and let $\pr$ be the
projection of $X$ onto $\calv$. Consider $\varrho $ a geodesic ray
in $X$, $\bar{r}$ a geodesic ray in $\calv$ and $r$ a lift of
$\bar{r}$ in $X$. The ray $r$ is contained in some Weyl chamber
of vertex $r(0)$. In this same Weyl chamber there exists a unique
ray $\varrho_1$ of vertex $r(0)$ and contained in the orbit
$\varrho G$. The Busemann function $f_r$ restricted to $\varrho_1$
has the form $-\beta_0 t$ for some $\beta_0\geq 0$. This implies
that, as soon as $\beta_0>0$, $\pr (\varrho_1)$ moves away in the
cusp of $\calv$ and the depth at which it moves away at time $t$
measured with respect to the ray $\bar{r}$ is $\beta_0 t$. Note
that among all the geodesic rays in $\varrho G$ with origin on the
horosphere $H(\bar{r})$, the ray $\varrho_1$ has the maximal depth
at moment $t$ with respect to $\bar{r}$.

\begin{ques}\label{oq}
For every $\beta \in (0,\beta_0)$, consider the set
$$
\calr_\beta =\left\{ \mathbf{u}\in U_+(\varrho)\; ;\;
-f_{\bar{r}}\left(\varrho (t) \mathbf{u} \right)\geq \beta
t\mbox{ infinitely many times as }t\to \infty \right\}\, .
$$

Is it true that $d=\dim_H \calr_\beta = (1-\beta )\dim
U_+(\varrho)$ and that $\hh^d \left( \calr_\beta \right)=\infty$ ?
\end{ques}

\subsection{Organization of the paper}

Section \ref{Pr} contains preliminaries on horoballs, symmetric
spaces and semisimple groups. The equidistribution results
Proposition \ref{equidunip} and Proposition \ref{equiduc} in
Section \ref{eq} play an important part in our arguments. In
particular the latter implies the counting results Proposition
\ref{cpct} and Corollary \ref{lSu}.

In Section \ref{symquadr} are given the formulas for the Busemann
functions in the ambient symmetric space ${\mathcal{P}}_{n+1}$ as
well as in the symmetric space associated to the quadric,
${\mathcal{P}}_{n+1}( L_\q )\, $. In Sections \ref{gho} and
\ref{trc} we study the geometry of horoballs of
$\mathcal{P}_{n+1}(L_{\frak q})$. The obtained results together
with the counting result Corollary \ref{lSu} yield the
equidistribution of rational vectors on rational quadrics as
formulated in Corollary \ref{counting3}, and also a more general
result, Proposition \ref{counting1}.

Section \ref{elip} contains the proof of Theorem \ref{T2}. The
notion of ubiquitous system is recalled in Section \ref{us}. In
Section \ref{trei} we show the relation between the set
$\mathcal{S}_\psi (\qq_\q)$ and the set of unipotents
$\ss_\Psi^0$. We end our argument in Section \ref{nonhip} by
exhibiting a local ubiquitous system and applying results from
\cite{BDV2}.

In Section \ref{sectrays} we prove results on the Hausdorff
dimension and measure of sets of locally geodesic rays moving away
in the cusp of a locally symmetric space. In Section
\ref{sectrays1} we study the case of rays of type $\wp $ in the
locally symmetric space $\mathcal{P}_{n+1}(L_{\frak q})/\Gamma$,
where $\Gamma$ is an arbitrary lattice in $SO_I(L_\q )$. In the
other two sections we deduce from (\ref{jar}) and (\ref{linrh})
respectively results about rays in the locally symmetric space
$\mathcal{T}_{n+1}$.

\me

\noindent {\bf{Acknowledgements}} : I am grateful to the referee
for pointing out to me the references \cite{BDV1} and \cite{BDV2},
which allowed me to improve the results that I had in the first
draft of this paper. I wish to thank Livio Flaminio
 for having explained to me the equidistribution result Proposition
  \ref{equidunip} in the case of $SL(2,\R )$ as well as for many useful
  conversations. I also thank Maurice Dodson for inspiring conversations
   and for providing useful references.

\me

\section{Preliminaries on (locally) symmetric
spaces}\label{Pr}

The reader acquainted with semisimple groups and symmetric spaces
may skip Sections \ref{not} to \ref{lssp} and refer to them only
when needed.

\subsection{Notation and conventions}\label{not}

We denote by $\pri^n $ {\it{the set of primitive integer vectors
in }}$\R^n$,
$$\{ (p_1,p_2,\dots ,p_n)\in
\Z^n \setminus \{ (0,\dots 0) \} \; ;\;  {\rm{gcd}}(p_1,p_2,\dots
,p_n)=1 \}\, ,$$ and we denote by $\pri_+^n $ the subset
$$\{ (p_1,p_2,\dots ,p_n)\in
\pri^n  \; ;\;  p_n>0 \mbox{ or } p_i >0\, ,\, p_{i+1}=\cdots
=p_n=0 \} \; .$$

In a metric space $(X,\dist)$, for any subset $A$ of $X$, we
denote by ${\mathcal{N}}_{a}(A)$ the set $$\{ x\in X \: ;\: \dist
(x,A)< a \}\, .$$ When $A=\{ x_0 \}$ then ${\mathcal{N}}_{a}(A)$
becomes an open ball and we use the notation $B(x_0, a)$.

\me

We denote by diag$(a_1,a_2,\dots ,a_n)$ the diagonal matrix having
entries $a_1,a_2,\dots ,a_n$ on the diagonal. In the particular
case when $a_1=\cdots =a_k=1$ and $a_{k+1}= \cdots =a_{k+\ell
}=-1\, ,\, k+\ell =n$, we denote by $I_{k,\ell }$ the diagonal
matrix. We denote by $Id_n$ the identity matrix.

Throughout by \textit{line} we mean a 1-dimensional linear
(sub)space.

Let $A$ be a subset of $\R^n$. We denote by $\R A$ the union of
all the lines intersecting $A$. We denote by $\proj A$ the image
of $\R A$ in $\proj^{n-1} \R$. If $B\subset \proj^{n-1} \R $ we
denote by $\R B$ the subset in $\R^n$ which is union of all lines
contained in $B$.

We denote by $\langle v_1,\dots ,v_k \rangle$ the linear subspace
generated by the vectors $v_1,\dots ,v_k$.

\me

Given two functions $f$ and $g$ with values in $\R$, we write
$f\ll g$ if $f(x)\leq C\cdot g(x)$, for every $x$, where $C >0$ is
a universal constant. We write $f\asymp g$ if both $f\ll g$ and
$f\gg g$ hold. We write $f \sim g$ if $\frac{f(x)}{g(x)}\to 1$
when $x\to \infty $. We denote by $\| f\|_\infty $ the supremum
norm of the function $f$.

\me

If $G$ is a group, we denote by $Z(G)$ its center $\{ z\in G \;
;\; zg=gz\, ,\, \forall g\in G\}$. If $H$ is a subgroup of $G$ we
denote by $C_G (H)$ the center of $H$ in $G$, that is the group
$\{ z\in G \; ;\;  zh=hz\, ,\, \forall h\in H\}$.

If $G$ is a topological group, we denote by $G_e$ its connected
component containing the identity.

Let $G$ be a Lie group. {\it{A lattice}} in $G$ is a discrete
subgroup $\Gamma $ of $G$ such that $G/ \Gamma $ has a finite
$G$-invariant measure induced by the Haar measure on $G$. If $G/
\Gamma $ is compact, the lattice is called {\it{uniform}},
otherwise it is called {\it{non-uniform}}.

If a group $G$ acts on a space $X$, for every point $x\in X$ we
denote by $G_x$ the stabilizer of $x$ in $G$, that is the subgroup
$\{g\in G \; ;\;  gx=x\}$.

\me

Let $X$ be a complete Riemannian manifold of non-positive
curvature. Two geodesic rays in $X$ are called {\it{asymptotic}}
if they are at finite Hausdorff distance one from the other. This
defines an equivalence relation $\sim $ on the set $\mathcal{R}$
of all geodesic rays in $X$. The boundary at infinity of $X$ is
the quotient $\mathcal{R}/\sim $. It is usually denoted by
$\partial_\infty X$. Given $\xi\in \partial_\infty X$ and a
geodesic ray in the equivalence class $\xi$, one writes $r(\infty
)=\xi$.

%The set of geodesic rays $\mathcal{R}$ can be endowed with a
%natural topology induced by the modified Hausdorff pseudo-distance
%(or the Gromov-Hausdorff pseudo-distance) as it is defined in
%\cite{BH}.\fn{de verificat}. This induces a topology on
%$\partial_\infty X$, sometimes called \textit{the cone topology}.

Let $r$ be a geodesic ray in $X$. {\it The Busemann function
associated to $r$} is the function
$$
f_r:X\to \R \, ,\; f_r(x)=\lim_{t\to \infty}[\dist (x,r(t))-t]\; .
$$

Since the function $t\to \dist (x,r(t))-t$ is non-increasing and
bounded, the limit exists. The level hypersurfaces $H_a(r)=\lbrace
x\in X \; ;\;  f_r(x)= a \rbrace$ are called {\it horospheres},
the sublevel sets $Hb_a(r)=\lbrace x\in X \; ;\;  f_r(x)\leq a
\rbrace$ are called {\it{closed horoballs}} and their interiors,
$Hbo_a(r)$, are called {\it open horoballs}. For $a=0$ we use the
notation $H(r)$ for the horosphere, and $Hb(r),\; Hbo(r)$ for the
closed and open horoball, respectively.

Suppose moreover that $X$ is simply connected.

Given an arbitrary point $x\in X$ and an arbitrary point at
infinity $\xi \in \partial_\infty X$, there exists a unique
geodesic ray $r$ with $r(0)=x$ and $r(\infty )=\xi$.
%On the other
%hand every oriented geodesic segment extends to a geodesic ray by
%the Hopf-Rinow Theorem \cite{GLP}. Hence one can define a
%bijection $\phi $ between the unitary tangent space in $x$,
%$S_xX$, and the boundary at infinity, by associating to each
%tangent vector $v\in S_xX$ the unique geodesic ray $r_v$ with
%$r_v(0)=x$ and tangent to $v$, and then putting $\phi (v) =r_v
%(\infty )$. The topology on $\partial_\infty X$ which makes $\phi $
% a homeomorphism is called \textit{the cone topology}.

The Busemann functions of two asymptotic rays in $X$ differ by a
constant \cite{BH}. Therefore we shall sometimes call them
{\it{Busemann functions of basepoint }}$\xi$, where $\xi $ is the
common point at infinity of the two rays. The families of
horoballs and horospheres are the same for the two rays. We shall
say that they are horoballs and horospheres {\it of basepoint
}$\xi$.

Two points $\xi$ and $\zeta $ in $\partial_{\infty }X$ are said to
be \textit{opposite} if there exists a complete geodesic $\frak G$
such that the point at infinity of ${\frak G}|_{[0,+\infty )}$ is
$\xi$ and the point at infinity of ${\frak G}|_{(-\infty ,0]}$ is
$\zeta$.

\begin{definition}
The {\textit{oriented distance}} odist$(Hb(r_1)\, ,\, Hb(r_2))$
{\textit{between two horoballs}} $Hb(r_1)$ and $Hb(r_2)$ of
opposite basepoints is $ \inf_{x\in Hb(r_2)} f_{r_1}(x) $.
\end{definition}

\subsection{Semisimple groups and symmetric spaces}\label{ssg}

Henceforth by $X$ we denote a symmetric space of non-compact type
without Euclidean factors, and by $G$ the connected component of
the identity in its group of isometries. Then $G$ is a semisimple
Lie group. We identify the symmetric space $X$ with $K\backslash
G$, where $K$ is a maximal compact subgroup of $G$. Hence we
consider the action of $G$ on $X$ by isometries \textit{to the
right}, and correspondingly we consider the action of $G$ on
itself by isometries to the right (with respect to the proper
metric), and the action by conjugation also to the right, i.e. $a
: G\to Aut\, (G)\, ,\, a(g_0)(g)= g_0^{-1}gg_0$. For the theory of
symmetric spaces and associated semisimple groups we refer to
\cite{He}.

We recall that every connected semisimple real Lie group is
isomorphic to the identity component of the real Lie group of real
points of a semisimple algebraic group. Therefore, one can always
talk about polynomial, rational and bi-rational maps on $G$.
Moreover $G$ has a faithful embedding $f: G \to \slnr $ such that
$f(G)^T=f(G)$ and $f(K)=f(G) \cap O(n,\R )$. Details can be found
for instance in \cite{OV}, \cite{Mo} or in \cite{Ra}.

\me

%%%%%%%%%%%%%%%%%add

\Notat We denote $\dist (e,g)$ by $|g|$, where $\dist $ is the
right invariant metric on $G$.

\me

An element $g_0$ in $\slnr$ is \textit{hyperbolic} if there exists
$g\in GL(n,\R)$ such that $g g_0 g\iv$ is diagonalizable with all
the eigenvalues real positive.

All the Lie groups considered in the sequel are real Lie groups,
unless otherwise stated.

Consider a field $\K \subset \R$. We say that a Lie group $G$ is
\textit{defined over $\K$} if $G$ has finitely many connected
components and if its connected component of the identity
coincides with the connected component of the identity of a real
algebraic group defined over $\K$ \cite[Definition 6.2]{Wi}.

A \textit{torus} is a closed connected Lie subgroup of $\slnr$
which is diagonalizable over $\C$, i.e. such that there exists
$g\in GL(n,\C)$ with the property that $g\, T\, g\iv $ is
diagonal. A torus is called $\K$-\textit{split} if it is defined
over $\K$ and diagonalizable over $\K$, that is if there exists
$g\in GL(n,\K)$ with the property that $g\, T\, g\iv $ is
diagonal.

A torus (and more generally a reductive group) is called
$\K$-\textit{anisotropic} if it is defined over $\K$ and if it
does not contain any non-trivial $\K$-split torus. Note that a
$\Q$-anisotropic torus $T$ has the property that the set of its
integer points $T_\Z$ is a lattice in it \cite{Bo}.

\me

\noindent \textit{Conventions}: Henceforth by \textit{torus} we
mean a non-trivial $\R$-split torus. The only exception is when we
talk about $\K$-anisotropic torus, in which case the word keeps
its general meaning. By wall/Weyl chamber we mean a
\textit{closed} wall/Weyl chamber.
 By its \textit{relative interior} we mean the open wall/Weyl chamber.

\me

 We call {\it{singular torus in}} $G$ a torus $A_0$ which,
 in every maximal torus $A$ containing it,
  can be written as $\bigcap_{\lambda \in \Lambda}\ker \lambda$,
   where $\Lambda $ is a non-empty set of roots on $A$.
    Any such torus is a union of walls of Weyl chambers.

Let $\triangleleft A_0$ be a wall or a Weyl chamber in the torus
$A_0$, and let $\triangleleft A_0^{op}$ be the opposite wall. We
consider the parabolic group corresponding to $\triangleleft
  A_0$,
$$
P(\triangleleft A_0 ) = \{ g\in G \: ; \: \sup_{n\in \N }
|\mathbf{a}^ng\mathbf{a}^{-n}| < + \infty \, ,\, \forall
\mathbf{a}\in \triangleleft A_0 \}\, ,
$$ and the unipotent group corresponding to $\triangleleft A_0$,
$$
U(\triangleleft A_0 ) = \{ g\in G \: ; \: \lim_{n\to \infty }
\mathbf{a}^ng\mathbf{a}^{-n} =e \, ,\, \forall \mathbf{a}\mbox{ in
the relative interior of } \triangleleft A_0 \}\, .
$$ We denote $U(\triangleleft A_0^{op} )$ by $U_+(\triangleleft A_0)$.

We have that $P(\triangleleft A_0 )=C_G(A_0)U(\triangleleft A_0
)=U(\triangleleft A_0 )C_G(A_0)$, $U(\triangleleft A_0 )$ is the
unipotent radical of $P(\triangleleft A_0 )$, and $P(\triangleleft
A_0 )$ is the normalizer of $U(\triangleleft A_0 )$ in $G$. The
center decomposes as $C_G(A_0)= M A_0=A_0M$, where $Z(M)$ is
compact and $M/Z(M)$ is semisimple. It follows that
$$
P(\triangleleft A_0 )=M A_0U(\triangleleft A_0 )=U(\triangleleft
A_0 )A_0M\, ,
$$ which is
called {\it{the Langlands decomposition of }}$P(\triangleleft A_0
)$.

\begin{remark}
The action of $M$ on $U(\triangleleft A_0 )$ by conjugation
preserves the Haar measure.
\end{remark}

\proof Any semisimple connected Lie group coincides with its
commutator subgroup (see for instance \cite[$\S 1.4.1$ and $\S 4.1.3$]{OV}), hence any linear representation of a semisimple group
    preserves the volume. Consequently $Ad(M)$ restricted to
the Lie algebra $\frak u$ of $U$ preserves the volume, which
yields the conclusion.\endproof

The Lie algebras ${\frak u}$ and ${\frak u}_{+}$ of
$U(\triangleleft A_0 )$ and $U_{+}(\triangleleft A_0)$ decompose
into proper subspaces for $Ad(A_0)$, ${\frak u}=\bigoplus_{\lambda
(\triangleleft A_0)>1}{\frak u}_{\lambda }$ and ${\frak
u}_{+}=\bigoplus_{\lambda (\triangleleft A_0)<1}{\frak u}_{\lambda
}$, respectively. Here $\lambda (\triangleleft A_0)>1$ signifies
that $\lambda >1$ when restricted to the relative interior of
$\triangleleft A_0$.

%Consequently, in the particular case when $A_0$ and $\triangleleft A_0$ are one-dimensional, for each element $\mathbf{a}$ in
%$\triangleleft A_0 \setminus \{ e\}$, the restriction of the
%conjugation map $g\mapsto \mathbf{a} g \mathbf{a}\iv$ to
%$U(\triangleleft A_0 )$ is a contracting homothety, and the
%restriction to $U_+(\triangleleft A_0 )$ is a dilating homothety,
%of factors inverse to each other.

The sets $P(\triangleleft A_0 )U_{+}(\triangleleft A_0 )$ and $U_{+}(\triangleleft A_0 )P(\triangleleft A_0 ) $ are open and Zariski dense in
$G$. Therefore they both give coordinate systems in $G$ near $e$.

Suppose that the group $G$ is defined over $\Q$, that $A_0$ is a
$\Q$-split torus and that $\triangleleft A_0$ is a $\Q$-wall or a
$\Q$-Weyl chamber in it. Then $C_G(A_0)$ and $U(\triangleleft A_0
)$ are also defined over $\Q$. Moreover $C_G(A_0) = M' A_0 = A_0
M'$, where $M'$ is defined over $\Q $, $Z(M')_e$ is a
$\Q$-anisotropic torus and $M'/Z(M')_e$ is semisimple. Recall that
in this case $\Gamma = G_\Z$ is a lattice in $G$, that
$U(\triangleleft A_0 )\cap \Gamma$ is a uniform lattice and that
$M'\cap \Gamma $ is a lattice in $M'$.
%Moreover the embedding
%\begin{equation}\label{emb}
%M'/M'\cap \Gamma \rightarrow G/\Gamma
%\end{equation}
%is proper.

For details on the previous results we refer to \cite{Bo},
\cite{Ra} and \cite{Wi}.

We recall that a \textit{flat} in $X$ is a totally geodesically
embedded copy of an Euclidean space in $X$, and that a
\textit{maximal flat} is a flat which is maximal with respect to
the inclusion. Every maximal flat $F$ is the orbit of a maximal
torus $A$. Given a point $x\in F$, a \textit{Weyl chamber} or
\textit{a wall} with vertex $x$ is a set of type $x \triangleleft
A_0$, where $\triangleleft A_0$ is a Weyl chamber or respectively
a wall in the torus $A$. A \textit{singular flat through} $x$ is
an orbit $x A_0$, where $A_0$ is a singular torus in $A$. In the
particular case when $G$ is defined over $\Q$, $A$, $A_0$ are
$\Q$-split, $\triangleleft A_0$ is a $\Q$-Weyl chamber or wall,
the corresponding maximal/singular flat, Weyl chamber or wall are
called $\Q$-maximal/singular flat, $\Q$-Weyl chamber and
$\Q$-wall, respectively.

The group $G$ acts transitively on the collection of maximal
flats, as well as on the collection of Weyl chambers in $X$. This
is equivalent to saying that it acts transitively by conjugation
on the collection of maximal tori and on the collection of Weyl
chambers in $G$. The stabilizer in $G$ of a Weyl chamber $W_0$ in
$X$ is a compact subgroup $K_0$. Therefore $K_0\backslash G$ can
be identified with the fiber bundle of the Weyl chambers in $X$.

\subsection{Geodesic rays, Busemann functions}\label{bdry0}

 Let ${\mathcal{A}}=(\mathbf{a}_t)$ be a one-parameter
subgroup of $G$ composed of hyperbolic elements and let
${\mathcal{A}}^{+}$ be the positive sub-semigroup
$(\mathbf{a}_t)_{t\geq 0}$. Let $r$ be a geodesic ray in $X$ such
that $r(t)=r(0)\mathbf{a}_t$ for every $t\geq 0$. We consider
$A_0$ either the minimal singular torus containing ${\mathcal{A}}$
or, if no such torus exists, the unique maximal
 torus containing ${\mathcal{A}}$. We have the equality $C_G({\mathcal{A}})=C_G(A_0)$.
  If $A_0$ has dimension
one we call the one-parameter group ${\mathcal{A}}$, the semigroup
${\mathcal{A}}^{+}$ and the geodesic ray $r$ \emph{maximal
singular}.

 Let $\triangleleft A_0$ be the wall/Weyl chamber
  containing ${\mathcal{A}}^{+}\setminus \{ e \}$ in its relative
  interior. We denote $P(\triangleleft A_0 )$, $U(\triangleleft A_0 )$
  and $U_+(\triangleleft A_0 )$ also by $P(r)$, $U(r)$ and $U_+(r)$,
  respectively, and we call them the {\it{parabolic}}, the {\it{unipotent}} and the
\textit{opposite (expanding) unipotent group} of the ray $r$.
%$$P(\triangleleft A_0 ) = \{ g\in G \: ; \: \sup_{t\in [0, +\infty )} |\mathbf{a}_tg\mathbf{a}_{-t}| < + \infty \}\: ,
% $$
% $$
% U(\triangleleft A_0 ) = \{ g\in G \: ; \: \lim_{t\to +\infty } \mathbf{a}_tg\mathbf{a}_{-t} =e \}
% $$ and
% $$ U_+(\triangleleft A_0 ) = \{ g\in G \: ; \: \lim_{t\to -\infty } \mathbf{a}_tg\mathbf{a}_{-t} =e \}\: .
% $$
%We denote $P(\triangleleft A_0 )$, $U(\triangleleft A_0 )$
  %and $U_+(\triangleleft A_0 )$ also by
  %$P({\mathcal{A}}^{+})$, $U({\mathcal{A}}^{+})$ and
  %$U_+({\mathcal{A}}^{+})$,
  %and we call them {\it{the parabolic, unipotent}} and {\it{opposite (or expanding)
   %unipotent subgroup of}} ${\mathcal{A}}^{+}$, respectively.
%In other words, $r([0, +\infty ))$ is an orbit of the positive sub-semigroup
%${\mathcal{A}}^{+}=(\mathbf{a}_t)_{t\geq 0}$.
The parabolic group $P(r)$ decomposes as $P(r)=\mathcal{A}P^0(r)$,
where $P^0(r)$ is a codimension 1 subgroup acting transitively
with compact stabilizer on every horosphere $H_a(r)$. We call
$P^0(r)$ {\it the horospherical group} of $r$.
%For every Weyl chamber $\triangleleft
%A_0 \supset {\mathcal{A}}^{+}$, $P^0(\triangleleft
%A_0,r)=P(\triangleleft A_0)\cap P^0(r)$ is a solvable group acting
%simply transitively on every horosphere $H_a(r)$. We call
%$P^0(\triangleleft A_0,r)$ a {\it solvable horospherical group of
%$r$}.

The following simple lemma will be useful in the future.

\begin{lemma}\label{Binv}
Let $r$ be a geodesic ray in the symmetric space $X$ and let
$\frak G$ be the unique geodesic containing it, parameterized by
arc length such that $r= {\frak G}|_{[0, +\infty )}$. Let $P^0$ be
the horospherical group of $r$. A function $\Psi :X\to \R$
 which is invariant with respect to $P^0$ and such that
  $\Psi ( {\frak G}(t) ) =-t\, ,\, \forall t\in \R $, coincides with $f_r$.
\end{lemma}

\proof For every $x\in X$ there exists a unique $t\in \R $ and
$p\in P^0$ such that $x={\frak G}(t)p$. We have $\Psi (x)=\Psi
({\frak G}(t)p)=\Psi ({\frak G}(t))=-t= f_r({\frak G}(t))=f_r(x)$.
\endproof

Consider the particular case when $G$ is defined over $\Q$ and when $\mathcal{A}^+ = \triangleleft A_0$ is a $\Q$-wall.
By the discussion in the end of Section \ref{ssg}, the horospherical group $P^0$ equals $M' U(\triangleleft A_0)=U(\triangleleft A_0) M'$,
 where $M'$ and $U(\triangleleft A_0)$ are defined over $\Q$.

\me

\subsection{Boundary at infinity}\label{bdry}

If $W$ is a Weyl chamber or a wall in $X$ then its boundary at
infinity $W(\infty )$ is a spherical simplex in $\partial_\infty
X$, also called \textit{spherical chamber} or respectively
\textit{spherical wall}. These simplices cover $\partial_\infty X$
and determine a structure of spherical building on it
(\cite[Chapters 15,16]{Mo}, \cite[Appendix 5]{BGS}).

Let $W_0$ be an arbitrary Weyl chamber in $X$. The group $G$ acts
on $\partial_\infty X$ on the right with fundamental domain
$W_0(\infty )$. Given a point $\xi$ in the relative interior of a
spherical wall $W(\infty )$, where $W=x \triangleleft A_0 $, the
stabilizer of $\xi$ is the stabilizer of the whole wall $W(\infty
)$, and it is the parabolic group $P(\triangleleft A_0)$. Since
any parabolic group acts transitively on $X$, it follows that for
every point $x\in X$ there exists a wall $W_x$ of vertex $x$ and
such that $W_x(\infty )=W(\infty )$.

Given a fixed a Weyl chamber $W_0$, $\partial_\infty X /G$ can be
identified with $W_0 (\infty )$, and one can define a projection
$\mathrm{sl}:
\partial_\infty X \to W_0 (\infty )$. The image $\mathrm{sl} (\xi )$ of every point
$\xi $ in $\partial_\infty X$ is called \textit{the slope of
}$\xi$. The \textit{slope of a geodesic ray} $r$ is the slope of
$r(\infty )$.

Let $x_0$ be an arbitrary point in $X$ and let $K$ be the maximal
compact subgroup fixing $x_0$. Given a wall $W$ with vertex $x_0$,
its stabilizer $K_W$ in $K$ is contained in the stabilizer $K_F$
of the minimal singular flat containing $W$, and it fixes both $W$
and $F$ pointwise. The group $K$ acts transitively on the set of
Weyl chambers of vertex $x_0$. Hence, given the stabilizer
$K_{W_0}$ of a Weyl chamber $W_0$ of vertex $x_0$, the quotient
$K_{W_0}\backslash K$ can be identified with the set of Weyl
chambers of vertex $x_0$. In particular, by the previous remarks,
$K$ acts transitively on the set of spherical chambers of
$\partial_\infty X$, and every spherical chamber $W_0(\infty )$
can be seen as the quotient $\partial_\infty X /K$.

\me

\subsection{Locally symmetric spaces}\label{lssp}

Let $\Gamma $ be a lattice in $G$. Here we shall be mainly
interested in non-uniform irreducible lattices in semisimple
groups of real rank at least $2$. By Margulis Arithmeticity
Theorem \cite{Wi} such a lattice $\Gamma$ is an arithmetic lattice
of $\Q$-rank $\mathbf{r}\geq 1$. The quotient space $\mathcal{V}=X
/\Gamma $ is a locally symmetric space. It contains finitely many
totally geodesic Euclidean sectors $W_1, \cdots ,W_m$, of
dimension $\mathbf{r}$, eventually glued to each other along
faces, such that $\calv$ is at finite Hausdorff distance of the
union $W_1\cup \cdots \cup W_m$. Every sector $W_1,\dots ,W_m$ is
the projection of a $\Q$-Weyl chamber. The quotient $\mathcal{V}$
can have several topological ends if and only if $\mathbf{r}=1$.
For details see \cite{BoS} and \cite{Le}.

\me

\noindent \textit{Notation}:  We denote by $\pr $ the projection
of $X$ onto $\mathcal{V}$ and by $\pr_G$ the projection of $G$
onto $G/\Gamma$.

\me

Given a geodesic ray $\bar{r}$ entering one of the sectors
$W_i, i\in \{ 1,\dots ,m \}$, the depth into the end containing
$W_i$ can be measured by the Busemann function $f_{\bar{r}}$ of
$\bar{r}$. If $\bar{r}$ is a face of dimension one of $W_i$, $i\in
\{1,\dots ,m \}$, then we call it a \textit{maximal singular cusp
ray.} Let $r$ be a lift of $\bar{r}$ in $X$.

\begin{remarks}\label{horob}
\begin{itemize}
    \item[(1)] For $a<0$ with $|a|$ large enough, the projection $\pr (Hb_a(r))$
is $Hb_a(\bar{r} )$.
    \item[(2)] There exists $\alpha =\alpha (\bar{r}) >0 $ such that
   $$ \lim_{a\to -\infty } \frac{\ln vol Hb_a (\bar{r})}{a}
   =\alpha \, .$$
\end{itemize}
\end{remarks}

\proof  (1) Since the projection $\pr $ is a contraction,
$f_{\bar{r}}(\pr (x))\leq f_{r}(x)\, ,\; \forall x\in X$. This
implies that $\pr (Hb_a(r)) \subset Hb_a(\bar{r} )$.

One can identify $\calv$ with a fundamental domain of $\Gamma $ in
$X$, contained in a Siegel set as in \cite[Theorem 15.5]{Bo}.
Suppose that $r$ is chosen so that under this identification
$\bar{r} $ becomes $r$. Obviously for $a<0$ with $|a|$ large
enough, $Hb_a (\bar{r} )$ coincides with the trace of $Hb_a(r)$ on
the fundamental domain. This implies that $Hb_a(\bar{r} )\subset
\pr (Hb_a(r))$.

(2) follows by looking at the form of the Siegel set as given in
\cite[Theorem 15.5]{Bo}.
\endproof

Suppose that $\Gamma$ is arithmetic. Then without loss of
generality we may suppose that $G$ admits a $\Q$-structure such
that $\bar{r}$ is the projection of a $\Q$-wall $r$. The
horospherical group $P^0(r)$ can be written as $M'U(\triangleleft
A_0)=U(\triangleleft A_0)M'$,
  with both $M'$ and $U(\triangleleft A_0)$ defined over $\Q$.
  In what follows, we denote $P^0(r)$ and $U(\triangleleft A_0)$ simply by $P^0$
   and respectively $U$.
  According to \cite[Corollary 7.13]{Bo}, $P^0\cap \Gamma $
   is commensurable to the semidirect product $(U\cap \Gamma )(M'\cap \Gamma )$. Therefore $P^0/(P^0\cap \Gamma)$
    and $P^0/(U\cap \Gamma )(M'\cap \Gamma )$ have
   a common finite covering. Now given $\mathcal{D}$ a
   fundamental domain of $U$ with respect to $U\cap \Gamma $ and
   $\mathcal{F}$ a fundamental domain of $M'$ with respect to $M'\cap \Gamma
   $, the set $\mathcal{F}\mathcal{D}$ is a fundamental domain of
   $P^0$ with respect to $(U\cap \Gamma )(M'\cap \Gamma )$.
   Indeed:
   \begin{itemize}
    \item $\ff \dd  (U\cap \Gamma )(M'\cap \Gamma )=\ff U (M'\cap \Gamma)=\ff (M'\cap
    \Gamma)U=M'U$;
    \item if $u\in U\cap \Gamma $ and $m\in M'\cap \Gamma$ are
    such that for some $f_i\in \ff$ and $d_i\in \dd$, $i=1,2$, $f_1d_1
    um=f_2d_2$, then $f_1m m\iv (d_1 u)m=f_2d_2$, whence
    $f_1m=f_2$ and $m\iv (d_1 u)m = d_2$. The former equality
    implies that $m=e$, the latter implies that $u=e$.
   \end{itemize}

%There exists $\tau_i <0$ such that
%if $f_i(\pr (x))\leq \tau_i$ then $f_i(\pi (x))= f_{\rho_i}(x)$.
%This is because the function $f_{\rho_i}$ is invariant with
%respect to the stabilizer $\Gamma_i$ of $\rho_i (\infty )$ in
%$\Gamma $ and for large $\tau_i $ the projection of $Hb_{-\tau_i}(\rho_i )$ into the quotient space $X/ \Gamma$ coincides with
%$Hb_{-\tau_i }(\rho_i ) / \Gamma_i$.

%We denote by $\mathcal{V}_i^{\tau_i}$ the projection $Hb_{-\tau_i
%}(\rho_i ) / \Gamma_i$. For $\tau_i$ large enough the sets
%$\mathcal{V}_i^{\tau_i}, i\in \{ 1,2,\dots ,k \}$, are pairwise
%disjoint and we can write $\mathcal{V}= K \cup \bigsqcup_{i=1}^k
%\mathcal{V}_i^{\tau_i}$, where $K$ is a compact subset. It is not
%difficult to deduce from this that there exists a constant
%$C=C(\mathcal{V})$ such that for every $i\neq j$ the inequality $
%f_i+f_j\leq C$ holds on $\mathcal{V}$. The pre-image of the set
%$\bigsqcup_{i=1}^k \mathcal{V}_i^{\tau_i}$ is composed of finitely
%many $\Gamma$-orbits of disjoint horoballs, $\bigsqcup_{i=1}^k
%\bigsqcup_{\gamma \in \Gamma_i \backslash \Gamma} Hb_{-\tau_i}(\rho_i ) \gamma $.

\me

%%%%%%%%%%%%%%%%%%%%%%
%%%%%%%%%%%%%%%%%%%%%%%%%%%%%%%
\subsection{Equidistribution results}\label{eq}

Let $G$ be a connected semisimple Lie group without compact
factors and with trivial center. Let
${\mathcal{A}}=(\mathbf{a}_t)$ be a one-parameter subgroup of $G$
composed of hyperbolic elements, and
 ${\mathcal{A}}^{+}=(\mathbf{a}_t)_{t\geq 0}$. Let $A_0$ be either the minimal singular torus containing ${\mathcal{A}}$
  or the unique maximal torus containing ${\mathcal{A}}$,
   and $\triangleleft A_0$ its unique wall/Weyl chamber containing
    ${\mathcal{A}}^{+}\setminus \{ e \}$ in its relative interior.
 Let $\mathfrak{C}=C_G(A_0)=C_G({\mathcal{A}})$, $P=P(\triangleleft A_0)$,
 $U=U(\triangleleft A_0)$
 and $U_{+}=U_+(\triangleleft A_0)$, endowed with their Haar measures.

\smallskip

\noindent {\it{Notation :}} For every subset $S$ of $G$, we
denote by $S_t$ the subset $a(\mathbf{a}_{-t}) S $. We denote by
$S^{-1}$ the image of $S$ under the inversion in $G$.

\smallskip

For $p\in P$ fixed, we consider the (partially defined) map
$\Psi_{p}$ from $ U_{+}$ to $U_{+}$, defined by

\begin{displaymath}
\Psi_{p} (\mathbf{u}_{+})=\mathbf{u}_{+}' \mbox{ such
that }P\mathbf{u}_{+} =P\mathbf{u}_{+}'p \; .
\end{displaymath}

\begin{figure}[!ht]
\centering
\includegraphics{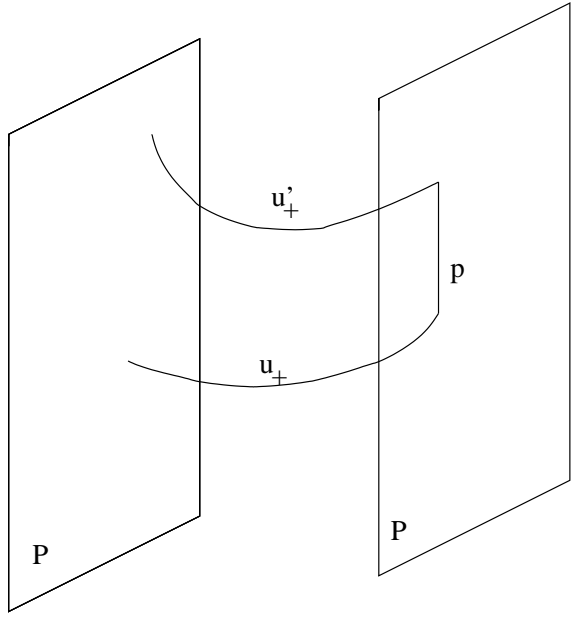}
\caption{The map $\Psi_{p}$}
\label{fig1}
\end{figure}

Let $D_p$ be its maximal domain of definition. Associated to this
map, we have the maps $\pi_p :D_p \to P$, $\upsilon_{p} : D_p\to
U$ and $\kappa_{p} : D_p\to {\mathfrak{C}}$ defined by the
relations
\begin{equation}\label{iotaz}
\Psi_{p} (\mathbf{u}_{+})p = \pi_p (\mathbf{u}_{+}) \mathbf{u}_{+}
\mbox{ and }\pi_p (\mathbf{u}_{+})= \upsilon_{p}(\mathbf{u}_{+})
\kappa_{p}(\mathbf{u}_{+}) \; .
\end{equation}

\Notat For $\alpha >0$ we define ${\mathcal{Q}}_\alpha=\{ p\in P
\: ;\: p=uc\, ,\, u\in U\, ,\, |u|\leq \alpha\, ,\, c\in
{\mathfrak{C}}\, ,\, |c|\leq \alpha \}$.

\begin{lemma}\label{psizu}
\begin{itemize}
  \item[(i)] For any $p\in P$, the domain $D_p$ is
  an open Zariski dense subset of $U_{+}$ and the map $\Psi_{p}$ is bi-rational.
   It satisfies the relation $a(\mathbf{a}_{t}) \circ \Psi_{p} = \Psi_{a(\mathbf{a}_{t})(p)}\circ
  a(\mathbf{a}_{t})$, where $a:G \to Aut(G)$ is the action to the right of G on
  itself by conjugation.
  \item[(ii)] Let $K$ be a compact subset in $P$. The set $\bigcap_{p\in K} D_p$ contains a neighborhood of $e$.
  \item[(iii)] Let $K$ be a compact subset in $P$ and $\Omega$ an open subset in $U_+$.
  The intersection $\bigcap_{p\in K} \Psi_{p} (\Omega \cap D_p)$ is open.

  \item[(iv)] The map $\Lambda_{p}  : D_p\to U_{+}$, $\Lambda_{p}(\mathbf{u}_{+})= \mathbf{u}_{+}^{-1} \Psi_{p}(\mathbf{u}_{+})$
   tends to the constant map equal to the identity element $e$ uniformly on
  compact subsets as $p\in {\mathcal{Q}}_\alpha$ and $\alpha \to 0$.
  \item[(v)] The Jacobian of the map $\Psi_{p}$, which we
  denote by $\left| \frac{d \Psi_{p}}{d\mathbf{u}_{+}} \right|$, tends to the constant
   map equal to $1$ uniformly on
  compact subsets as $p\in {\mathcal{Q}}_\alpha$ and $\alpha \to 0$.
  \item[(vi)] The map $S_\alpha (\mathbf{u}_{+})=\sup_{p\in {\mathcal{Q}}_\alpha }
 \left( |\upsilon_{p}(\mathbf{u}_{+})| + | \kappa_{p}(\mathbf{u}_{+})| \right) $
  tends to zero uniformly on
  compact subsets as $p\in {\mathcal{Q}}_\alpha$ and $\alpha \to 0$.
\end{itemize}
\end{lemma}

\proof (i) Let $D_p = U_+ \cap PU_+ p$, which is an open Zariski
dense subset of $U_+$. For every $\mathbf{u}_{+} \in D_p$,
$\mathbf{u}_{+}p^{-1}=\overline{p}\, \overline{\mathbf{u}}_{+}\in
PU_+$. It follows that
$\Psi_{p}(\mathbf{u}_{+})=\overline{\mathbf{u}}_{+}$. The map from
$PU_{+}$ to $U_{+}$ defined by ${\overline{p}}\,
{\overline{\mathbf{u}}}_{+} \to {\overline{\mathbf{u}}}_{+}$ is a
rational map. Hence the map $\Psi_{p}$ is rational. Moreover,
since the converse map is $\Psi_{p^{-1}}$, the map $\Psi_{p}$ is
bi-rational. The behavior of $\Psi_{p}$ with respect to the action
of the group $(\mathbf{a}_t)$ can be deduced by applying
$a(\mathbf{a}_t)$ in the first equality in (\ref{iotaz}).

\me

(ii) Suppose that $\bigcap_{p\in K } D_p$ does not contain a
neighborhood of $e$. Then there exists a sequence
$\mathbf{u}^{+}_n \to e$ and a sequence $p_n \in K$ such that
$\mathbf{u}^{+}_n p_n^{-1} \not\in PU_+ $ for any $n\in \N $. A
subsequence of $\mathbf{u}^{+}_n p_n^{-1}$ converges to some
$p_0\in P \subset PU_+$. This contradicts the fact that $PU_+$ is
open.

\me

(iii) We prove that $\bigcup_{p\in K }\complement \Psi_{p} (\Omega
\cap D_p)$ is closed. Let $\mathbf{u}^{+}_n$ be a sequence in this
set, converging to $\mathbf{u}_{+}$. For every $n\in \N $ there
exists $p_n\in K$ such that $\mathbf{u}^{+}_n \not\in \Psi_{p_n}
(\Omega \cap D_{p_n})$. Up to taking a subsequence, $p_n$
converges to $p\in K$. Suppose that $\mathbf{u}_{+}\in \Psi_{p}
(\Omega \cap D_p)$, which is equivalent to the fact that
$\mathbf{u}_{+} p\in P\Omega $. Since $P\Omega $ is an open set in
$G$ and $\mathbf{u}^{+}_n p_n \to \mathbf{u}_{+} p$, for some $n$
large enough we have $\mathbf{u}^{+}_n p_n \in P\Omega $, that is $ \mathbf{u}^{+}_n \in \Psi_{p_n} (\Omega \cap D_{p_n})$. This
contradicts the hypothesis.

From its definition it is straightforward that when $p\in
{\mathcal{Q}}_\alpha $ and $\alpha \to 0$ the map $\Psi_{p}\to Id$
uniformly in the $C^1$ topology on compact subsets. This implies
(iv), (v) and (vi).
\endproof

\me

Let $\Omega $ be a relatively compact open neighborhood of $e$ in
$U_{+}$ and let $\alpha $ be a small positive number. We suppose
that $\alpha \leq 1$ and that $\Omega$ is small enough to be
contained in $\bigcap_{p\in {\mathcal{Q}}_1 \cup
{\mathcal{Q}}_1^{-1} } D_p$.

\begin{definitions}
The $(\alpha\, ,\, \Omega )$-\textit{box of basepoint} $g_0$ is
the set
\begin{displaymath}
  Box_{\alpha,\Omega } (g_0) = \bigcup_{p\in {\mathcal{Q}}_\alpha }\; \bigcup_{\mathbf{u}_{+}\in \Psi_{p}(\Omega )} \mathbf{u}_{+}pg_0\;
  .
\end{displaymath}

We call
\begin{displaymath}
\mathbf{a}_tBox_{\alpha,\Omega } (g_0)= \bigcup_{p\in
({\mathcal{Q}}_\alpha)_t}\; \bigcup_{\mathbf{u}_{+}\in
\Psi_{p}(\Omega_t)} \mathbf{u}_{+}pg_t \; ,
\end{displaymath} where $g_t=\mathbf{a}_tg_0$, \textit{the
$t$-pushed } $(\alpha\, ,\, \Omega )$-box.
\end{definitions}

We denote by $\delta (\alpha\, ,\, \Omega )$ and $S(\alpha\, ,\,
\Omega )$ the maximal values, for $p\in {\mathcal{Q}}_\alpha $, of
$\sup_{\mathbf{u}_{+}\in \Omega } \left|\, \left| \frac{d
\Psi_{p}}{d\mathbf{u}_{+}} \right| -1 \right|$ and of
$\sup_{\mathbf{u}_{+}\in \Omega }S_\alpha (\mathbf{u}_{+})$
respectively.

\begin{definitions}
We call $\varepsilon $-\textit{base of }$\Omega $ an open
relatively compact neighborhood $\Omega_{\varepsilon }$ of $e$ in
$U_{+}$ such that $\nu (\Omega \vartriangle \Omega
\mathcal{K})\leq \varepsilon \nu (\Omega )$, for every nonempty
$\mathcal{K} \subset \Omega_{\varepsilon }$, where $\nu $ is the
Haar measure on $U_{+}$ and $A \vartriangle B = (A\setminus B
)\cup (B\setminus A)$.

We call $(\varepsilon ,\alpha) $-\textit{base of }$\Omega $ any
subset of $U_+$ of the form $\Omega_{\varepsilon , \alpha
}=\bigcap_{p\in {\mathcal{Q}}_\alpha}\Psi_{p^{-1}}\left(
\Omega_{\varepsilon } \cap D_{p^{-1}} \right)$, where
$\Omega_{\varepsilon }$ is an $\varepsilon $-base of $\Omega $.
According to Lemma \ref{psizu}, (iii), $\Omega_{\varepsilon ,
\alpha }$ is an open relatively compact neighborhood of $e$.

We call $Box_{\alpha,\Omega_{\varepsilon , \alpha }} (g_0)$ an
$\varepsilon$-\textit{base of the box }$Box_{\alpha,\Omega }
(g_0)$.
\end{definitions}

\begin{figure}[!ht]
\centering
\includegraphics{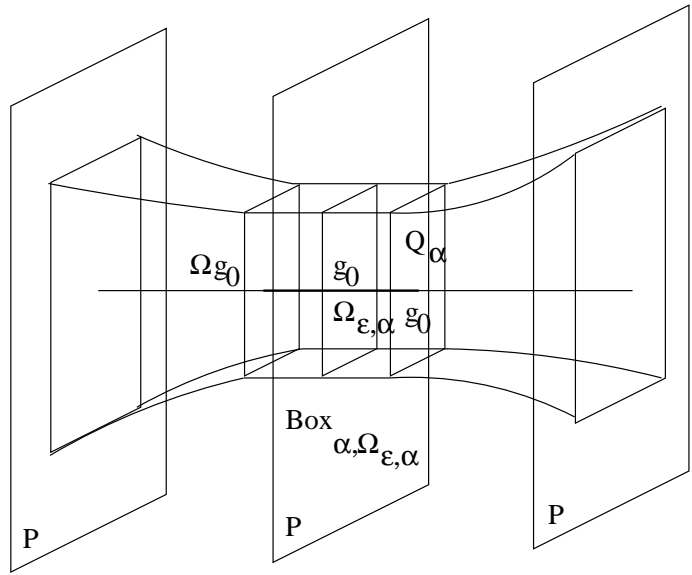}
\caption{A box and an $\varepsilon$-base of it} \label{fig2}
\end{figure}

\Notat We denote $\oint_D f\, d\mu = \frac{1}{\mu (D)}\int_D f\,
d\mu $.

\me

\begin{lemma}\label{ebox}
Let $f$ be a bounded uniformly continuous function on $G$, of
modulus of continuity $\omega $. Let $\Omega $ be a sufficiently
 small relatively compact open neighborhood of $e$ in $U_{+}$,
  let $\varepsilon $ be a small positive number and let $\Omega_\varepsilon$
   be an $\varepsilon$-base of $\Omega$. For $\alpha >0$
sufficiently small we have that
\begin{equation}\label{diff}
  \left| \oint_{\Omega_t} f(\mathbf{u}_{+}y)\, d\nu (\mathbf{u}_{+})- \oint_{\Omega_t} f(\mathbf{u}_{+}y')\, d\nu (\mathbf{u}_{+})
  \right|= O\left( \omega \left( S(\alpha , \Omega )\right) +\| f\|_\infty \left[\varepsilon + \delta \left( \alpha , \Omega \Omega_\varepsilon \right) \right]
  \right)\, ,
\end{equation} for every $g_0\in G$, every $t\in [0,+\infty )$ and every $y,y'$ in the $t$-pushed
$\varepsilon$-base $\mathbf{a}_tBox_{\alpha,\Omega_{\varepsilon
,\alpha } } (g_0)$.
\end{lemma}

\proof We suppose that $\alpha <1$ and that $\Omega \subset
\bigcap_{p\in \mathcal{Q}_1\cup \mathcal{Q}_1^{-1}}D_p$. It
suffices to prove (\ref{diff}) for $y'=g_t$. Since $y\in
\mathbf{a}_t Box_{\alpha,\Omega_{\varepsilon , \alpha } } (g_0)$,
we may write $y= \mathbf{a}_t\widetilde{\mathbf{u}}_{+}pg_0$,
where $p\in {\mathcal{Q}}_\alpha$ and
$\widetilde{\mathbf{u}}_{+}\in \Psi_{p}(\Omega_{\varepsilon ,
\alpha } )\subset \Omega_{\varepsilon }$. Then
$$
\oint_{\Omega_t} f(\mathbf{u}_{+}y)\, d\nu (\mathbf{u}_{+}) =
\oint_{\Omega_t}
f(\mathbf{u}_{+}\mathbf{a}_t\widetilde{\mathbf{u}}_{+}pg_0)\, d\nu
(\mathbf{u}_{+}) =
$$
$$
\oint_{\Omega} f(\mathbf{a}_t \eta_{+}
\widetilde{\mathbf{u}}_{+}pg_0)\, d\nu (\eta_{+} )=\oint_{\Omega
\widetilde{\mathbf{u}}_{+}} f(\mathbf{a}_t \eta_{+} pg_0)\, d\nu
(\eta_{+} )\, .
$$

By the definition of $\Omega_{\varepsilon }$, we have
that $\nu \left( \Omega \vartriangle \Omega
\widetilde{\mathbf{u}}_{+}\right) \leq \varepsilon \nu (\Omega )
$. It follows that
$$
\oint_{\Omega \widetilde{\mathbf{u}}_{+}} f(\mathbf{a}_t \eta_{+}
pg_0)\, d\nu (\eta_{+} ) =
\oint_{\Omega } f(\mathbf{a}_t \eta_{+}
pg_0)\, d\nu (\eta_{+} ) + O\left(
\varepsilon \| f\|_\infty \right)\, .
$$

We want to change the integration domain from $\Omega $ to
$\Psi_{p}(\Omega )$. We consider $\alpha $ small enough so that
$\Lambda_{p}(\Omega )\subset \Omega_\varepsilon$, for any $p\in
\mathcal{Q}_\alpha \cup \mathcal{Q}_\alpha^{-1}$.

We have
$$
\Psi_{p}(\Omega ) = \{\Psi_{p}(\mathbf{u}_{+})\; ;\;
\mathbf{u}_{+}\in \Omega \} =\{
\mathbf{u}_{+}\Lambda_{p}(\mathbf{u}_{+})\; ;\;  \mathbf{u}_{+}\in
\Omega \}\subset \Omega \Omega_\varepsilon \, .
$$

Therefore $\nu (\Psi_{p}(\Omega
)\setminus \Omega )\leq \varepsilon \nu (\Omega )$.

We write $\Omega \setminus \Psi_{p}(\Omega )$
  as the image under $\Psi_{p}$ of $\Psi_{p^{-1}}(\Omega )\setminus
\Omega$. An argument similar to the previous implies that
$\Psi_{p^{-1}}(\Omega )\subset \Omega \Omega_\varepsilon $, hence $\nu (\Psi_{p^{-1}}(\Omega )\setminus \Omega )\leq
\varepsilon \nu (\Omega )$. Since the Jacobian of $\Psi_{p}$
differs from $1$ by $O\left(\delta\left(\alpha\, ,\, \Omega
\Omega_\varepsilon \right)\right)$ on $\Psi_{p^{-1}}(\Omega )$, we
may conclude that $\nu (\Omega \setminus \Psi_{p}(\Omega ))\leq
\varepsilon \nu (\Omega )\left[ 1+ \delta\left(\alpha\, ,\, \Omega
\Omega_\varepsilon \right) \right]$. Consequently
$$
\oint_{\Omega } f(\mathbf{a}_t \eta_{+} pg_0)\, d\nu (\eta_{+} ) =
\frac{1}{\nu (\Omega)}\int_{\Psi_{p}(\Omega )} f(\mathbf{a}_t
\eta_{+} pg_0)\, d\nu (\eta_{+} ) + O\left(\varepsilon \|f
\|_\infty \left[ 1+ \delta\left(\alpha\, ,\, \Omega
\Omega_\varepsilon \right)\right] \right)\, .
$$

With the change $\eta_{+} = \Psi_{p}(\eta_{+}' )$ we may write
$$
\frac{1}{\nu (\Omega)}\int_{ \Psi_{ p }(\Omega )}
f\left(\mathbf{a}_t \eta_{+} pg_0\right)\, d\nu (\eta_{+} ) =
\oint_{\Omega } f\left(\mathbf{a}_t \Psi_{p }(\eta_{+}'
)pg_0\right)\, d\nu (\eta_{+}')+ O\left(\, \delta(\alpha ,\Omega
)\; \|f \|_\infty \, \right)\, .
$$

Using the notation in (\ref{iotaz}) we may write
$$
\oint_{\Omega } f\left(\mathbf{a}_t \Psi_{p }(\eta_{+}' )
pg_0\right)\, d\nu (\eta_{+}') = \oint_{\Omega }
f\left(\mathbf{a}_t \upsilon_{p}(\eta_{+}' ) \kappa_{p}(\eta_{+}'
) \eta_{+}' g_0\right)\, d\nu (\eta_{+}')\, .
$$

By the right invariance of the metric on $G$ we have that
$$
\dist \left( \mathbf{a}_t \upsilon_{p}(\eta_{+}' )
\kappa_{p}(\eta_{+}' ) \eta_{+}' g_0 \, ,\, \mathbf{a}_t\eta_{+}'
g_0 \right) = \dist \left(a(\mathbf{a}_{-t})\left(
\upsilon_{p}(\eta_{+}' ) \right) \kappa_{p}(\eta_{+}' )\, ,\, e
\right)\leq S(\alpha , \Omega )\, .
$$

Therefore the last integral is equal to
$$
\oint_{\Omega } f(\mathbf{a}_t\eta_{+}' g_0 )\, d\nu (\eta_{+}') +
O( \omega (\, S(\alpha , \Omega ))\, ) = \oint_{\Omega_t }
f(\mathbf{u}_{+} g_t )\, d\nu (\mathbf{u}_{+}) + O( \omega (\,
S(\alpha , \Omega ))\, )\; .
$$ \endproof

\begin{proposition}\label{equidunip}
Let $\Gamma $ be an irreducible lattice in $G$ and let $f: G/
\Gamma \to \R$ be a function which is uniformly continuous and
bounded. Let $(\mathbf{a}_t)_{t\in \R }$ be a one-parameter group
composed of hyperbolic elements and let $U_+$ be the expanding
unipotent subgroup corresponding to $(\mathbf{a}_t)_{t\geq 0}$.
Let $\Omega $ be an open relatively compact set in $U_{+}$, with
the property that there exists $\mathbf{u}_0\in \Omega\iv$ such
that for any $t_0\in \R$, the family of sets $a(\mathbf{a}_{-t})
\left( \Omega \mathbf{u}_0\right) \, ,\, t\in [t_0 , +\infty )$,
is a summing net for $U_{+}$, in the sense of \cite[$\S
4.15$]{Pa}. For any $\bar{g}_0 \in G/ \Gamma $,
\begin{displaymath}
  \oint_{\Omega } f(\mathbf{a}_t\mathbf{u}_{+}\bar{g}_0 )\, d\nu (\mathbf{u}_{+})
  \rightarrow \oint_{G/\Gamma } f\, d \mu \mbox{ as } t \to +\infty \; ,
\end{displaymath} where $\mu $ is the measure on $G/\Gamma$ induced
by the Haar measure on $G$.
\end{proposition}

\proof \textit{Step 1}.\quad We suppose that $\Omega$ is a
neighborhood of $e$ contained in $\bigcap_{p\in \mathcal{Q}_1\cup
\mathcal{Q}_1^{-1}}D_p$ and that $\mathbf{u}_0=e$. We denote by
$\omega $ the modulus of continuity of $f$. We fix arbitrary
$\bar{g}_0 \in G/ \Gamma $ and $\varepsilon $ small positive
number. We consider $\Omega_\varepsilon $ an $\varepsilon$-base of
$\Omega$ and $\alpha >0$ sufficiently small so that the conclusion
of Lemma \ref{ebox} holds, and also so that $\omega (S(\alpha ,
\Omega))\leq \varepsilon $ and $\delta \left( \alpha , \Omega
\Omega_\varepsilon \right)\leq \frac{\varepsilon }{\|
f\|_{\infty}}$.

\smallskip

The group $U_{+}$ acts ergodically on $G/\Gamma $ \cite[$\S
2.2$]{Zi}. This and the fact that the family of sets $\Omega_t =
a(\mathbf{a}_{-t})\left( \Omega \right)$ is a summing net implies
that $\oint_{\Omega_t } f(\mathbf{u}_{+}\bar{g})\, d\nu
(\mathbf{u}_{+})$, seen as a function of $\bar{g} \in G/\Gamma $,
converges to $\oint_{G/\Gamma} f\, d\mu $ in $L^2(G/\Gamma)$ as
$t\to \infty$ \cite[$\S 5$]{Pa}. This implies that for the given
$\varepsilon $ and for any small $\lambda >0$ there exists
$T=T(\varepsilon , \lambda , \Omega )$ such that for every $t\geq
T$ the set of points $\bar{g} \in G/\Gamma $ satisfying
\begin{equation}\label{e7}
\left| \oint_{\Omega_t } f(\mathbf{u}_{+}\bar{g})\, d\nu
(\mathbf{u}_{+}) \; -\; \oint_{G/\Gamma} f\, d\mu   \right| \geq
\varepsilon
\end{equation} has measure strictly smaller than $\lambda$. We take $\lambda $ to be the measure
 of the projection in $G/\Gamma$ of the $\varepsilon$-base
 $Box_{\alpha , \Omega_{\varepsilon , \alpha } }(g_0)$. Hence, for every $t\geq T$,
  at least one point $\bar{y}$
 in the projection of the $t$-pushed $\varepsilon$-base satisfies the inequality opposite to (\ref{e7}). This
and Lemma \ref{ebox} imply that for every $t\geq T$,
\begin{displaymath}
\left| \oint_{\Omega_t } f(\mathbf{u}_{+}\bar{g}_t)\, d\nu
(\mathbf{u}_{+}) \; -\; \oint_{G/\Gamma} f\, d\mu   \right| < C
\varepsilon \; ,
\end{displaymath} where $C=C(\|f\|_\infty )$ and $\bar{g}_t = \mathbf{a}_{t}\bar{g}_0$.

\medskip

\textit{Step 2}.\quad We suppose that $\mathbf{u}_0=e$ and that
$\Omega$ is an arbitrary open relatively compact neighborhood of
$e$ in $U_{+}$ satisfying the hypothesis and not necessarily
contained in $\bigcap_{p\in \mathcal{Q}_1\cup
\mathcal{Q}_1^{-1}}D_p$. There exists $\tau \in (0,+\infty )$ so
that $\Omega_{-\tau} = a(\mathbf{a}_{\tau})\left( \Omega \right)
\subset \bigcap_{p\in \mathcal{Q}_1\cup \mathcal{Q}_1^{-1}}D_p$.
We apply the result obtained in Step 1 to $\Omega_{-\tau}$ and
with a change of variables we obtain the same result for $\Omega$.

\medskip

\textit{Step 3}.\quad We consider the general case. By Step 2 we
have the conclusion of the Proposition for $\Omega \mathbf{u}_0$.
This implies the conclusion for $\Omega$.\endproof

\begin{proposition}\label{equiduc}
Let $\Gamma $ be a lattice in $G$ and $f: G/ \Gamma \to \R$ a
bounded uniformly continuous function. Let $(\mathbf{a}_t)_{t\in
\R }$ be a one-parameter group composed of hyperbolic elements and
let $U_+$ be the expanding unipotent subgroup corresponding to
$(\mathbf{a}_t)_{t\geq 0}$ and $\frak C$ be the center of the
group in $G$. Let $\Omega $ be an open relatively compact set in
$U_{+}$, with the same property as in Proposition \ref{equidunip}
and let $\Phi $ be a finite volume submanifold in $\frak C$. We
denote by $\vartheta$ the volume on $\Phi $. For any $\bar{g}_0
\in G/ \Gamma $,
\begin{displaymath}
  \oint_{\Phi \Omega } f(\mathbf{a}_t \mathbf{c} \mathbf{u}_{+} \bar{g}_0 )\, d\nu
  (\mathbf{u}_{+})d\vartheta (\mathbf{c})
  \rightarrow \oint_{G/\Gamma } f\, d \mu \mbox{ as } t \to +\infty \; ,
\end{displaymath} where $\mu $ is the measure on $G/\Gamma$ induced
by the Haar measure on $G$.
\end{proposition}

\proof We fix an arbitrary small positive number $\varepsilon$.
There exists a compact subset $K=K(\varepsilon )$ in $\Phi $ such
that $\vartheta(\Phi  \setminus K)< \frac{\varepsilon }{2 \|
f\|_\infty }\, \vartheta (\Phi)$. It follows that
$$
\left| \oint_{\Phi \Omega  } f(\mathbf{a}_t \mathbf{c}
\mathbf{u}_{+}\bar{g}_0 )\, d\nu
  (\mathbf{u}_{+})d\vartheta (\mathbf{c}) - \oint_{ K \Omega } f(\mathbf{a}_t \mathbf{c} \mathbf{u}_{+} \bar{g}_0 )\, d\nu
  (\mathbf{u}_{+})d\vartheta (\mathbf{c}) \right|< \varepsilon \,
  .
$$

Let $\omega$ be the modulus of continuity of $f$. Let $\delta
>0$ be such that $\omega( \delta ) < \varepsilon $. By compactness
of $K$, there exist $\mathbf{k}_1,\dots ,\mathbf{k}_m$ in $K$ so
that the set of balls $B(\mathbf{k}_i, \delta ),\, i\in \{ 1,\dots
,m \}$, is a cover of $K$. Thus, for every $\mathbf{c}\in K$ there
exists $i\in \{ 1,\dots ,m \}$ such that $\dist ( \mathbf{c},
\mathbf{k}_i)< \delta $, which by the right invariance of the
distance and the fact that the projection is a contraction implies
that $\dist ( \mathbf{c} \mathbf{a}_t\mathbf{u}_+ \bar{g}_0\, ,\,
\mathbf{k}_i\mathbf{a}_t\mathbf{u}_+ \bar{g}_0)< \delta$, for any
$t\in \R , \mathbf{u}_+\in U_+$ and $g_0\in G$. It follows that
$$
\left| f(\mathbf{c} \mathbf{a}_t\mathbf{u}_+ \bar{g}_0 ) -
f(\mathbf{k}_i \mathbf{a}_t\mathbf{u}_+ \bar{g}_0 )\right| <
\varepsilon \, ,
$$ whence
$$
\left| \oint_{\Omega } f(\mathbf{c} \mathbf{a}_t \mathbf{u}_+
\bar{g}_0) d\nu (\mathbf{u}_+) - \oint_{\Omega }
f(\mathbf{k}_i\mathbf{a}_t \mathbf{u}_+ \bar{g}_0) d\nu
(\mathbf{u}_+) \right| <\varepsilon \, ,
$$ for any $t\in \R $ and $g_0\in G$.

Now
$$\oint_{\Omega } f(\mathbf{k}_i\mathbf{a}_t \mathbf{u}_+
\bar{g}_0) d\nu (\mathbf{u}_+)=\oint_{\mathbf{k}_i \Omega
\mathbf{k}_i\iv } f(\mathbf{a}_t \mathbf{u}_+ \mathbf{k}_i
\bar{g}_0) d\nu (\mathbf{u}_+)\, .$$
 The set $\mathbf{k}_i \Omega
\mathbf{k}_i\iv $ also satisfies the hypothesis of Proposition
\ref{equidunip}. It follows that there exists $T>0$ such that for
any $t\geq T$ and any $i\in \{ 1,\dots ,m \}$,
$$
\left|\oint_{\mathbf{k}_i \Omega \mathbf{k}_i\iv  } f(\mathbf{a}_t
\mathbf{u}_+ \mathbf{k}_i\bar{g}_0) d\nu (\mathbf{u}_+) -
\oint_{G/\Gamma } f\, d \mu \right| <\varepsilon \, .
$$

We conclude that for every $t\geq T$
$$
\left|\oint_{K \Omega } f(\mathbf{c} \mathbf{a}_t \mathbf{u}_+
\bar{g}_0) d\nu (\mathbf{u}_+) d\vartheta (\mathbf{c}) -
\oint_{G/\Gamma } f\, d \mu \right| <2\varepsilon
$$ and that
$$
\left|\oint_{\Phi \Omega } f(\mathbf{c} \mathbf{a}_t \mathbf{u}_+
\bar{g}_0) d\nu (\mathbf{u}_+) d\vartheta (\mathbf{c}) -
\oint_{G/\Gamma } f\, d \mu \right| <3\varepsilon \, .
$$

\endproof

\subsection{Counting results}\label{count}

Throughout the whole of this section we work in the following
setting. We consider $\Gamma $ an irreducible lattice in $G$ and
$\mathcal{V}=X/\Gamma $ the corresponding locally symmetric space.
Let $\bar{r}$ be a maximal singular cusp ray in $\mathcal{V}$, let
$r$ be a lift of it and let $\xi= r(\infty)$. We denote by $Hb_t$
the horoball $Hb_t(r)$ and by $H_t$ its boundary horosphere. When
$t=0$ we drop the index. Consider $P=P(r)$, $P^0=P^0(r)$ and
$U=U(r)$ the parabolic, horospherical and respectively the
unipotent group of $r$. Let $A_0=(\mathbf{a}_t)_{t\in \R}$ be the
one-dimensional singular torus such that $r(t)=r(0)\mathbf{a}_t$
for every $t\geq 0$. Assume that $\mathbf{u} \mapsto
a(\mathbf{a}_{t})(\mathbf{u})=\mathbf{a}_{-t} \mathbf{u}
\mathbf{a}_{t}$ is a dilating homothety on $U$ of factor
$e^{\lambda t}$, with $\lambda
>0$. We denote the topological dimension of $U$ by $\Delta $.

\me

\noindent \textit{Notation}: For every $y\in X$ we denote by
$K(y)$ the maximal compact subgroup of $G$ stabilizing $y$.
According to Section \ref{bdry}, we can identify $\xi G$ with $\xi
K(y)$ and with $K(y)_\xi \backslash K(y)$. Every such
identification endows $\xi G$ with a natural measure coming from
the measure on $K(y)_\xi \backslash K(y)$, which we denote by
$\theta_y$.

Every open set $\Omega$ in  $\xi G$ can be identified with an open
set in $K(y)_\xi \backslash K (y)$. We denote by $\Omega_{K(y)}$
its pre-image in $K(y)$, also open. The set $\Omega_{K(y)}$ is the
maximal set in $K(y)$ such that $\Omega =\xi \Omega_{K(y)}$.
%It represents the set of Weyl chambers of vertex $y$ containing one of the points in $\Omega $ in their boundaries at infinity.

\me

%%%%%%%%%%%%%%
%%%%%%%%%%%%%%

\begin{proposition}\label{cpct}
Let ${\mathcal{O}}$ be an open set in $\xi G$, let $x$ be a point
in $ X$ and let $T>0$. For every $k\in \N$, let $N_x\left( k ,\oo
\right)$ be the number of horoballs $Hb \gamma$, $\gamma \in
\Gamma $, with basepoint in ${\mathcal{O}}$ and such that $\dist
(x, Hb \gamma )$ is in $[kT , (k+1)T)$. For any $T\geq T_0(G)$ and
any $x\in X$ we have that
$$
\bc_1 e^{\lambda (k+1)T \Delta}\theta_x (\oo ) \leq  N_x\left( k
,\oo \right) \leq \bc_2 e^{\lambda (k+1)T \Delta}\theta_x (\oo )\;
\mbox{ for every }k\geq k_0(x,\oo , Hb)\, ,
$$ where $\bc_i=\bc_i (G,\Gamma )$ for $i=1,2$.
\end{proposition}

\proof We fix an arbitrary point $x$ in $X$ and an open set $\oo$
in $\xi G$. We put $K$ for $K(x)$. We also fix a Weyl chamber
$W_0$ having $r$ as a face and we denote its stabilizer in $G$ by
$K_0$. Since $P$ acts transitively on $X$ it follows that there
exists $p\in P$ such that $r(0)p=x$. Then $W_0p$ is a Weyl chamber
of vertex $x$, containing $\xi$ in its boundary at infinity, and
$W_0p \oo_{K}$ is the set of Weyl chambers of vertex $x$
containing one of the points in $\oo $ in their boundaries at
infinity.

The set of Weyl chambers with vertices on $H_s$ containing $\xi$
in their boundaries at infinity is $W_0\mathbf{a}_{-s}P^0$.
Consider the horoball $Hb \gamma $. It has its basepoint $\xi
\gamma$ in ${\mathcal{O}} $ and it is at distance at least $kT$
and smaller than $(k+1)T$ from $x$ if and only if for some $s\in
[kT,(k+1)T)$, $W_0 \mathbf{a}_{-s}P^0 \gamma \cap W_0p \oo_{K}
\neq \emptyset $. This is equivalent to $K_0\mathbf{a}_{-s}P^0
\gamma \cap p \oo_{K} \neq \emptyset $. Since $K_0$ commutes with
$A_0$ and it is contained in $P^0$, we have that
$K_0\mathbf{a}_{-s}P^0=\mathbf{a}_{-s}P^0$. Thus, $N_x\left( k
,\oo \right)$ is the cardinal of the set
$$
\left\{ \gamma \in (\Gamma \cap P^0)\backslash \Gamma \; ;\;
\exists s\in [kT,(k+1)T) \mbox{ such that }p {\mathcal{O}}_K \cap
\mathbf{a}_{-s} P^0 \gamma \neq \emptyset \right\}=
$$
$$
\left\{ \gamma \in \Gamma/(\Gamma \cap P^0) \; ;\;  \exists s\in
[kT,(k+1)T) \mbox{ such that }p{\mathcal{O}}_K \gamma \cap
\mathbf{a}_{-s} P^0 \neq \emptyset \right\}=
$$
$$
\left\{ \gamma \in \Gamma \; ;\;  \exists s\in [kT,(k+1)T) \mbox{
such that }p{\mathcal{O}}_K \gamma \cap \mathbf{a}_{-s}
P^0/(\Gamma \cap P^0) \neq \emptyset \right\}\, ,
$$ where in the last set either we may suppose that we are in $G/(\Gamma \cap
P^0)$, or we may suppose that we are in $G$, in which case
$P^0/(\Gamma \cap P^0)$ signifies a fundamental domain of $P^0$
with respect to $\Gamma \cap P^0$.

\textbf{Case $1$.} Suppose that $G$ has real rank at least $2$.
Then $\Gamma $ is an arithmetic lattice and $r$ is a maximal
singular $\Q$-wall. According to Section \ref{bdry0}, the
horospherical group $P^0$ is equal to $UM'=M'U$, where $U/U\cap
\Gamma $ is compact and $M'/M'\cap \Gamma$ has finite volume. By
Section \ref{lssp}, $P^0/(\Gamma \cap P^0 )$ and $P^0/(U\cap
\Gamma )(M'\cap \Gamma )$ have a common finite covering.
Consequently, we can replace in the counting problem above the
former by the latter. Also according to Section \ref{lssp}, if
$\mathcal{F}$ is a fundamental domain of $M'$ with respect to
$M'\cap \Gamma$, and $\mathcal{D}$ is a fundamental domain of $U$
with respect to $U\cap \Gamma$ then $\ff \dd $ is a fundamental
domain of $P^0$ with respect to $(U\cap \Gamma)(M'\cap \Gamma )$.
The counting problem above becomes the counting of the number of
times when $\pr_G \left( \mathbf{a}_{-s}\ff\dd \right) $
intersects $\pr_G \left(p \oo_K \right)$ in $G/\Gamma $, for $s\in
[kT,(k+1)T)$.

\me

\noindent \textit{Notation}: For a small positive number
$\varepsilon$ we denote by $B_\varepsilon $ the open ball $ B(x,
\varepsilon )$ in $X$.

We also consider $\oo^{\varepsilon -}\subset \oo \subset
\oo^{\varepsilon +}$ with $\oo^{\varepsilon -}$ and
$\oo^{\varepsilon +}$ two open subsets of $\xi G$ very near $\oo$.

\me

\begin{figure}[!ht]
\centering
\includegraphics{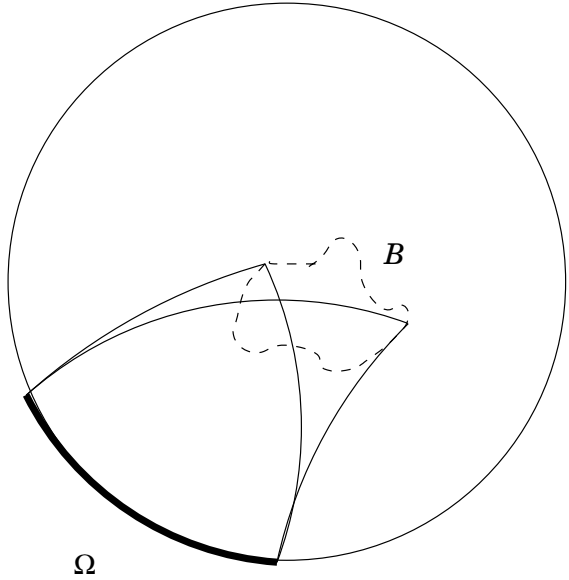}
\caption{Example of a set $W_0\Omega_{\mathcal{B}}$} \label{fig3}
\end{figure}

The map $$\Pi :K_0\backslash G \to K(r(0))\backslash G \times
P\backslash G\simeq X\times \xi G$$ defined by the two projections
is $C^\infty$. Therefore for two open sets $\mathcal{B}$ in $X$
and $\Omega $ in $\xi G$, $\Pi\iv (\mathcal{B}\times \Omega)$ is
an open set in $K_0\backslash G$. We denote its pre-image in $G$
by $\Omega_{\mathcal{B}}$. The set $\Omega_{\mathcal{B}}$ is the
maximal set in $G$ with the property that
$W_0\Omega_{\mathcal{B}}$ is the set of Weyl chambers with
vertices in $\mathcal{B}$ and containing a point from $\Omega$ in
their boundary at infinity. A picture of a set
$W_0\Omega_{\mathcal{B}}$ can be seen in Figure \ref{fig3}, in the
particular case when $X$ is the hyperbolic disk $\mathbb{D}^2$.

\me

\noindent \textbf{Upper estimate.} We want to define a continuous
function $f_1$ on $G$ taking values in $[0,1]$ and such that $f_1$
is $1$ on $\oo_{B_\varepsilon}$ and $0$ outside $\oo^{\varepsilon
+}_{B_{2\varepsilon }}$. There exists a continuous function
$f_1^1:X\to [0,1]$, $f_1^1=1$ on $B_\varepsilon$ and $f_1^1=0$
outside $B_{2\varepsilon }$. If $\oo^{\varepsilon+}$ is well
chosen then there exists a continuous function $f_1^2:\xi G \to
[0,1]$, $f_1^2=1$ on $\oo$ and $f_1^2=0$ outside
$\oo^{\varepsilon+}$. The function $f:X\times \xi G \to [0,1]$
defined by $f(x,\xi g)=f_1^1(x)f_1^2(\xi g)$ is 1 on
$B_\varepsilon \times \oo$ and 0 outside $B_{2\varepsilon }\times
\oo^{\varepsilon+} $, and by means of it and of the projection
$\Pi$ one can obtain the function $f_1$.

If $\varepsilon $ is small enough, $f_1$ can also be seen as a
function in $G/\Gamma$. Note that $f_1$ is a bounded function with
compact support, hence uniformly continuous. Proposition
\ref{equiduc} applied to $f_1$, the semigroup
$(\mathbf{a}_{-s})_{s\geq 0}$, $\dd \subset U$ and $\ff \subset
M'$ (or their relative interiors), and $g_0=e$ gives that
$$
\oint_{\ff \dd } f_1 (\mathbf{a}_{-s} \mathbf{c} \mathbf{u}
\bar{e}) d\nu (\mathbf{u}) d\vartheta (\mathbf{c }) \rightarrow
\oint_{G/\Gamma } f_1 d\mu\; \mbox{ when }\; s\to \infty \, ,
$$ where $\nu $ and $\vartheta $ are the
measures induced from the Haar measures on $U$ and $M'$,
respectively. In particular, for $s\geq s(f_1)$, we have

$$
\oint_{\ff \dd } f_1 (\mathbf{a}_{-s} \mathbf{c} \mathbf{u}
\bar{e}) d\nu (\mathbf{u}) d\vartheta (\mathbf{c }) \leq 2
\oint_{G/\Gamma } f_1 d\mu \leq 2 vol(B_{2\varepsilon })[\nu_K
(\oo_K )+ \chi ] \, ,
$$ where $\nu_K$ is the Haar measure on $K$ and $\chi \to 0$ when $\varepsilon \to 0 $ and $\oo^{\varepsilon +}$ converges to
$\oo$.

This is equivalent to the fact that for $s\geq s(f_1)$,
$$
\oint_{\ff \dd_s} f_1 (\mathbf{c} \mathbf{u}
\bar{\mathbf{a}}_{-s}) d\nu (\mathbf{u}) d\vartheta (\mathbf{c })
\leq 2 vol(B_{2\varepsilon })[\nu_K (\oo_K )+ \chi] \,  ,
$$ where $\dd_s = \mathbf{a}_{-s} \dd \mathbf{a}_{s}$. We can
rewrite the previous inequality as
$$
\int_{\ff \dd_s} f_1 (\mathbf{c} \mathbf{u} \bar{\mathbf{a}}_{-s})
d\nu (\mathbf{u}) d\vartheta (\mathbf{c }) \leq 2\vartheta (\ff )
e^{\lambda s \Delta } \nu (\dd ) vol(B_{2\varepsilon })\left[\nu_K
(\oo_K )+ \chi \right]\, .
$$

By integration we obtain
\begin{equation}\label{integr}
\int_t^{t+\tau }\int_{\ff \dd_s} f_1 (\mathbf{c} \mathbf{u}
\bar{\mathbf{a}}_{-s}) d\nu (\mathbf{u}) d\vartheta (\mathbf{c
})\; ds\, \leq C_1 e^{(t+\tau ) \lambda \Delta } \; ,
\end{equation} where $C_1= \frac{2}{\lambda \Delta}\,\vartheta (\ff ) \nu (\dd
) vol(B_{2\varepsilon })\left[\nu_K (\oo_K )+ \chi \right] $.

We have that
\begin{equation}\label{inequ}
\int_t^{t+\tau }\int_{\ff \dd_s} f_1 (\mathbf{c} \mathbf{u}
\bar{\mathbf{a}}_{-s}) d\nu (\mathbf{u}) d\vartheta (\mathbf{c
})\; ds \geq \int_t^{t+\tau }\int_{\ff \dd_s}
\unu_{\oo_{B_\varepsilon }} (\mathbf{c} \mathbf{u}
\bar{\mathbf{a}}_{-s}) d\nu (\mathbf{u}) d\vartheta (\mathbf{c
})\; ds \; .
\end{equation}

The second term is the same as
 $$
\int_t^{t+\tau }\int_{\ff \dd_s} \unu_{\bigsqcup_{\gamma \in
\Gamma }  \oo_{B_\varepsilon }\gamma }(\mathbf{c} \mathbf{u}
\mathbf{a}_{-s}) d\nu (\mathbf{u}) d\vartheta (\mathbf{c })\; ds =
\int_t^{t+\tau }\int_{\ff \dd} \unu_{\bigsqcup_{\gamma \in \Gamma
} \oo_{B_\varepsilon }\gamma }(\mathbf{a}_{-s} \mathbf{c}
\mathbf{u}) d\nu (\mathbf{u}) d\vartheta (\mathbf{c })
e^{\lambda\Delta s}\; ds =
 $$
$$
\sum_{\gamma \in \Gamma } \int_t^{t+\tau }\int_{\ff \dd}
\unu_{\oo_{B_\varepsilon }\gamma }(\mathbf{a}_{-s} \mathbf{c}
\mathbf{u}) d\nu (\mathbf{u}) d\vartheta (\mathbf{c })
e^{\lambda\Delta s}\; ds \, .
 $$
Since $\ff \dd$ and a fundamental domain $\pp$ of $P^0$ with
respect to $P^0 \cap \Gamma$ have a common finite covering, the
integral above is equivalent (in the sense of the relation
$\asymp$ defined in Section \ref{not}) to the integral
$$
\sum_{\gamma \in \Gamma } \int_t^{t+\tau }\int_{\pp }
\unu_{\oo_{B_\varepsilon }\gamma }(\mathbf{a}_{-s} \mathbf{c}
\mathbf{u}) d\nu (\mathbf{u}) d\vartheta (\mathbf{c })
e^{\lambda\Delta s}\; ds  = $$
$$
\sum_{\gamma \in \Gamma /P^0 \cap \Gamma } \sum_{\gamma_1 \in P^0
\cap \Gamma } \int_t^{t+\tau }\int_{\pp } \unu_{\oo_{B_\varepsilon
}\gamma \gamma_1 }(\mathbf{a}_{-s} \mathbf{c} \mathbf{u}) d\nu
(\mathbf{u}) d\vartheta (\mathbf{c }) e^{\lambda\Delta s}\; ds =
$$
$$
\sum_{\gamma \in \Gamma /P^0 \cap \Gamma } \sum_{\gamma_1 \in P^0
\cap \Gamma } \int_t^{t+\tau }\int_{\pp \gamma_1 }
\unu_{\oo_{B_\varepsilon }\gamma}(\mathbf{a}_{-s} \mathbf{c}
\mathbf{u}) d\nu (\mathbf{u}) d\vartheta (\mathbf{c })
e^{\lambda\Delta s}\; ds =
$$
$$
\sum_{\gamma \in \Gamma /P^0 \cap \Gamma } \int_t^{t+\tau
}\sum_{\gamma_1 \in P^0 \cap \Gamma } \int_{\pp \gamma_1 }
\unu_{\oo_{B_\varepsilon }\gamma}(\mathbf{a}_{-s} \mathbf{c}
\mathbf{u}) d\nu (\mathbf{u}) d\vartheta (\mathbf{c })
e^{\lambda\Delta s}\; ds = $$
\begin{equation}\label{sum}
\sum_{\gamma \in \Gamma /P^0 \cap \Gamma } \int_t^{t+\tau }
\int_{P^0 } \unu_{\oo_{B_\varepsilon }\gamma}(\mathbf{a}_{-s}
\mathbf{c} \mathbf{u}) d\nu (\mathbf{u}) d\vartheta (\mathbf{c })
e^{\lambda\Delta s}\; ds\, .
\end{equation}

The projection $G\to K_0 \backslash G$ sends $K_0\mathbf{a}_{-s}
P^0=\mathbf{a}_{-s} P^0$ onto the set of Weyl chambers containing
the point $\xi $ in their boundaries at infinity and with vertices
on the horosphere $H_s$. The image of the set $K_0
\oo_{B_\varepsilon }\gamma = \oo_{B_\varepsilon }\gamma$ is the
set of Weyl chambers containing each one point from $\oo \gamma$
in their boundaries at infinity, and with vertices in
$B_\varepsilon \gamma$.

This shows that $\unu_{\oo_{B_\varepsilon }\gamma}$ takes the
value $1$ on $\mathbf{a}_{-s} P^0$ for some $s\in [t,t+\tau ] $ if
and only if $\xi \in \oo \gamma$ and $B_\varepsilon \gamma$
intersects the closed strip $Hb_{t+\tau }\setminus Hbo_t$.

Let $\gamma (P^0\cap \Gamma )$ be such that $B_\varepsilon \gamma$
is entirely contained in the closed strip $Hb_{t+\tau }\setminus
Hbo_t$ and such that $\xi \in \oo \gamma$. The first condition is
equivalent to $d=\dist (x\gamma , H)\in [t+\varepsilon ,t+\tau
-\varepsilon ]$. The second condition is equivalent to the fact
that $p\oo_K \gamma \cap \mathbf{a}_{-d} P^0 \neq \emptyset $. We
can write
$$
\int_t^{t+\tau } \int_{P^0 } \unu_{\oo_{B_\varepsilon
}\gamma}(\mathbf{a}_{-s} \mathbf{c} \mathbf{u}) d\nu (\mathbf{u})
d\vartheta (\mathbf{c }) e^{\lambda\Delta s}\; ds \geq
\int_{d-\varepsilon }^{d+\varepsilon } \int_{P^0 }
\unu_{\oo_{B_\varepsilon }\gamma}(\mathbf{a}_{-s} \mathbf{c}
\mathbf{u}) d\nu (\mathbf{u}) d\vartheta (\mathbf{c })
e^{\lambda\Delta s}\; ds
$$

The latter term can be rewritten as
$$
\int_{-\varepsilon }^{+\varepsilon } \int_{P^0 }
\unu_{\oo_{B_\varepsilon }\gamma}( \mathbf{c}
\mathbf{u}\mathbf{a}_{-(d+s)}) d\nu (\mathbf{u}) d\vartheta
(\mathbf{c })\; ds\, = vol (B_\varepsilon ) \nu_{K_\xi} (K_\xi )\,
.
$$

Here $\nu_{K_\xi}$ is the measure induced on $K_\xi $ by the Haar
measure on $P^0$.

It follows that the sum in (\ref{sum}) is larger than $N\, vol
(B_\varepsilon ) \nu_{K_\xi} (K_\xi )$, where $N$ is the number of
$\gamma \in \Gamma /(\Gamma \cap P^0 )$ such that $d=\dist (x
\gamma , H)\in [t+\varepsilon,t+\tau -\varepsilon]$ and $p\oo_K
\gamma \cap \mathbf{a}_{-d} P^0 \neq \emptyset $. If we put
$t+\varepsilon = kT$ and $t+\tau -\varepsilon =(k+1)T$ then $N$ is
nothing else than $N_x\left( k ,\oo \right)$. Inequalities
(\ref{integr}) and (\ref{inequ}) imply that
\begin{displaymath}
N_x\left( k ,\oo \right)\, vol(B_\varepsilon)\leq \mathbf{c}_2
vol(B_{2\varepsilon})\left[ \nu_K (\oo_K ) +\chi \right]
e^{\lambda \Delta \varepsilon } e^{\lambda \Delta (k+1)T }\, ,
\end{displaymath} where $\mathbf{c}_2=\mathbf{c}_2(G,\Gamma )$.

For small $\varepsilon $ the ratio $\frac{vol(B_{2\varepsilon})
}{vol(B_{\varepsilon})}$ is bounded from above. Also, if
$\varepsilon$ is small enough and $\oo^{\varepsilon+}$ is close
enough to $\oo$, we have that $\chi < \nu_K (\oo_K ) $ and that
$e^{\lambda \Delta \varepsilon }<2$. Thus we obtain
$$
N_x\left( k ,\oo \right) \leq \bc_2\, \nu_K (\oo_K ) e^{\lambda
\Delta (k+1)T }\, , \mbox{ for } k\geq k_0(\oo , x, Hb)\, .
$$

\noindent \textbf{Lower estimate.} Consider $f_2:G\to [0,1]$ a
continuous function defined to be $1$ on $\oo^{\varepsilon
-}_{B_{\varepsilon }}$ and $0$ outside $\oo_{B_{2\varepsilon }}$.
Such a function can be constructed for $\oo^{\varepsilon -}$ well
chosen in a way similar to $f_1$. For small $\varepsilon $, $f_2$
can also be seen as a function in $G/\Gamma$. We apply Proposition
\ref{equiduc} as for the upper estimate, but with the function
$f_2$. We get
$$
\oint_{\ff \dd } f_2 (\mathbf{a}_{-s} \mathbf{c} \mathbf{u}
\bar{e}) d\nu (\mathbf{u}) d\vartheta (\mathbf{c }) \rightarrow
\oint_{G/\Gamma } f_2 d\mu\, ,
$$ when $s\to \infty $. As above, we obtain that for $s\geq s'(f_2)$, we have

$$
\oint_{\ff \dd_s} f_2 (\mathbf{c} \mathbf{u}
\bar{\mathbf{a}}_{-s}) d\nu (\mathbf{u}) d\vartheta (\mathbf{c })
\geq \frac{1}{2} \oint_{G/\Gamma } f_2 d\mu \geq \frac{1}{2}
vol(B_{\varepsilon })[\nu_K (\oo_K )- \chi' ] \, ,
$$ where $\chi' \to 0$ when $\varepsilon \to 0 $ and $\oo^{\varepsilon -}$ converges to
$\oo$.

Computations similar to the previous yield
\begin{displaymath}
\int_t^{t+\tau }\int_{\ff \dd_s} f_2 (\mathbf{c} \mathbf{u}
\bar{\mathbf{a}}_{-s}) d\nu (\mathbf{u}) d\vartheta (\mathbf{c
})\; ds\, \geq C_2' e^{t\lambda \Delta } \; ,
\end{displaymath} where $C_2'=
\frac{1}{2}\frac{e^{\tau \lambda \Delta}-1}{\lambda
\Delta}\,\vartheta (\ff ) \nu (\dd )vol(B_{\varepsilon
})\left[\nu_K (\oo_K )- \chi' \right] $.

For $\tau$ large enough $e^{\tau \lambda \Delta}-1 \geq
\frac{1}{2}e^{\tau \lambda \Delta} $ and we have
\begin{equation}\label{integr2}
\int_t^{t+\tau }\int_{\ff \dd_s} f_2 (\mathbf{c} \mathbf{u}
\bar{\mathbf{a}}_{-s}) d\nu (\mathbf{u}) d\vartheta (\mathbf{c
})\; ds\, \geq C_2 e^{(t+\tau )\lambda \Delta } \; ,
\end{equation} where $C_2=
\frac{1}{4\lambda \Delta}\,\vartheta (\ff ) \nu (\dd )
vol(B_{\varepsilon })\left[\nu_K (\oo_K )- \chi'  \right] $.

Now we write
\begin{equation}\label{inequ2}
\int_t^{t+\tau }\int_{\ff \dd_s} f_2 (\mathbf{c} \mathbf{u}
\bar{\mathbf{a}}_{-s}) d\nu (\mathbf{u}) d\vartheta (\mathbf{c
})\; ds \leq \int_t^{t+\tau }\int_{\ff \dd_s}
\unu_{\oo_{B_{2\varepsilon }}} (\mathbf{c} \mathbf{u}
\bar{\mathbf{a}}_{-s}) d\nu (\mathbf{u}) d\vartheta (\mathbf{c
})\; ds \; .
\end{equation}

The previous argument implies that the second term in (\ref{inequ2}) is equivalent to the
 sum

 \begin{equation}\label{sum2} \sum_{\gamma \in \Gamma /P^0 \cap
\Gamma } \int_t^{t+\tau } \int_{P^0 } \unu_{\oo_{B_{2\varepsilon }
}\gamma}(\mathbf{a}_{-s} \mathbf{c} \mathbf{u}) d\nu (\mathbf{u})
d\vartheta (\mathbf{c }) e^{\lambda\Delta s}\; ds \, .
\end{equation}

The considerations above imply that the sum in (\ref{sum2}) is
smaller than $N'\, vol (B_{2\varepsilon} ) \nu_{K_\xi} (K_\xi )$,
where $N'$ is the number of $\gamma \in \Gamma /\Gamma \cap P^0 $
such that $B_{2\varepsilon} \gamma$ intersects $Hb_{t+\tau }
\setminus Hbo_t$ and such that $\xi \in \oo \gamma $. These
conditions are equivalent to the fact that the distance $d=\dist
(x\gamma , H)$ is in $[t-2\varepsilon , t+\tau +2\varepsilon ]$
and to $p\oo_K \gamma \cap \mathbf{a}_{-d} P^0 \neq \emptyset $.
Consequently, if we choose $t-2\varepsilon = kT$ and $t+\tau +
2\varepsilon =(k+1)T$, we obtain $N'=N_x\left( k ,\oo \right)$.
Inequalities (\ref{integr2}) and (\ref{inequ2}) imply that
\begin{displaymath}
N_x\left( k ,\oo \right) \, vol(B_{2\varepsilon})\geq \mathbf{c}_1
vol(B_{\varepsilon})\left[ \nu_K (\oo_K ) -\chi'
\right]e^{-2\varepsilon \lambda \Delta } e^{\lambda \Delta (k+1)T
}\, ,
\end{displaymath} whence we obtain that for $\varepsilon$ small
enough
$$
N_x\left( k ,\oo \right) \geq \bc_1\, \nu_K (\oo_K ) e^{\lambda
\Delta (k+1)T }\, \mbox{ for every }k\geq k_1(\oo ,x,Hb).
$$

We now note that $\nu_K (\oo_K ) \asymp \theta_x(\oo )$.

\noindent \textbf{Case 2.} Suppose that $G$ has real rank $1$.
Then the horospherical group $P^0$ is equal to the unipotent group
$U$, and we have that $U/U\cap \Gamma $ is compact. One can repeat
the same argument as previously, simplified by the fact that there
is no central factor $M'$, hence $\ff$ no longer appears, use
Proposition \ref{equidunip} instead of \ref{equiduc}, and get the
same estimate.\endproof

\begin{cor}\label{lSu}
Let $\rho$ be a geodesic ray such that $\rho (\infty ) $ is
opposite to $\xi$ and let $U=U(\rho )$. The subset $\xi U$ is open
and dense in $\xi G$, it can be identified with $U$ and thus
equipped with a measure induced from the measure on $U$, which we
denote by $\nu_u$. Let $\Omega $ be a relatively compact open
subset of $\xi U$ and let $T>0$. For every open subset $\oo$ of
$\Omega$, we denote by $N_r\left( k ,\oo \right)$ the number of
horoballs $Hb \gamma \, ,\gamma \in \Gamma$, with basepoints $\xi
\gamma \in \oo$ and such that the oriented distance $\mathrm{odist}\, \left (Hb(\rho ),
Hb \gamma \right)$ is in $ [kT, (k+1)T)$. For any $T\geq T_0(G,\rho ,
\Omega )$ we have that
$$
\mathbf{K}_1 e^{\lambda (k+1)T\Delta }\nu_u(\oo ) \leq N_r\left( k
,\oo \right)\leq \mathbf{K}_2 e^{\lambda (k+1)T\Delta }\nu_u(\oo
)\, ,\; \mbox{ for every }\; k\geq k_0(\oo , \rho , Hb)\, ,
$$ where $\mathbf{K}_i=\mathbf{K}_i(G, \Gamma , \rho , \Omega )$
for $i=1,2$.
\end{cor}

\proof Let $\frak G$ be a geodesic line such that ${\frak
G}(+\infty )=\xi $, ${\frak G}(-\infty )=\rho (\infty ) $ and
${\frak G}(0)\in H(\rho )$. Let $\Omega_U$ be the relatively
compact open subset of $U$ such that $\xi \Omega_U =\Omega$. The
set ${\frak G}(0)\Omega_U$ is a relatively compact subset of
$H(\rho )$ of diameter $D$. We choose a point $x$ in it. Let
$\widetilde{Hb}$ be a horoball with basepoint $\widetilde{\xi}$ in
$\Omega$ and such that $\mathrm{odist} \left( Hb(\rho),
\widetilde{Hb}\right)$ is positive and large enough. We have that
$\dist \left(x, \widetilde{Hb}\right)\geq \mathrm{odist}
\left(Hb(\rho ), \widetilde{Hb}\right)$.

Let $\widetilde{u}\in \Omega_U$ be such that $\xi
\widetilde{u}=\widetilde{\xi}$ and let $\widetilde{{\frak
G}}={\frak G}\widetilde{u}$. Let $t>0$ be such that
$\widetilde{{\frak G}}(t)$ is the entrance point of
$\widetilde{{\frak G}}$ into $\widetilde{Hb}$. Then
$\mathrm{odist} \left(Hb(\rho ), \widetilde{Hb}\right) =t$.

We have that $\dist \left(x, \widetilde{Hb}\right) \leq \dist
\left(x, \widetilde{{\frak G}}(t) \right)\leq t+\dist \left(x,
\widetilde{{\frak G}}(0) \right)=t+\dist \left(x, {\frak
G}(0)\widetilde{u} \right)\leq t+D$. Overall we obtain
$$
\mathrm{odist} \left(Hb(\rho ), \widetilde{Hb}\right)\leq \dist
\left(x, \widetilde{Hb}\right)\leq \mathrm{odist} \left(Hb(\rho ),
\widetilde{Hb}\right)+D\, .
$$

This inequality, together with Proposition \ref{cpct} and the fact that on $\Omega$ the two measures $\nu_u$ and
$\theta_x$ are equivalent implies the desired conclusion.\endproof

\section{Symmetric spaces of positive definite quadratic
forms}\label{symquadr}

\subsection{The ambient space}\label{str}

 Throughout the paper we shall
identify a quadratic form $Q$ on $\R^s $ with its matrix $M_Q$ in
the canonical basis of $\R^s$. The matrix of $Q$ in some other
basis $\mathcal{B}$ of $\R^s$ shall be denoted by
$M_Q^{\mathcal{B}}$. We shall denote by $b_Q$ the bilinear form
associated to $Q$.

Let ${\mathcal{P}}_s = SO(s) \backslash SL(s, \R )$. This space
can be identified with the space of positive definite quadratic
forms of determinant one on $\R^s$ by associating to each right
coset $SO(s)\, Y$ the quadratic form $Q_Y$ whose matrix in the
canonical basis is $M_Y=Y^T\cdot Y$.

We recall that ${\mathcal{P}}_s$ is equipped with a canonical
metric defined as follows. Given $Q_1\, ,\, Q_2 \in
{\mathcal{P}}_s$, there exists an orthonormal basis with respect
to $Q_1$ in which $Q_2$ becomes diagonal with coefficients
$\lambda_1,\dots ,\lambda_s\in \R_+^*$. We define
\begin{equation}\label{dist}
d(Q_1,Q_2)=\sqrt{\sum_{i=1}^s (\ln \lambda_i)^2}\; .
\end{equation}

Let $\mathcal{Q}_s$ be the space of quadratic forms on $\R^s$ and
$\mathcal{P}\mathcal{Q}_s$ the space of positive definite
quadratic forms. The group $GL(s,\R)$ acts on the right on
$\mathcal{Q}_s$ by
\begin{displaymath}
\Phi : GL(s,\R)\times \mathcal{Q}_s\to \mathcal{Q}_s\, ,\, \Phi
(B,M)=B^TMB\; .
\end{displaymath} This action can be written in terms of
quadratic forms as $\Phi (B,Q)=Q[B]=Q\circ B$.

The space $\mathcal{P}\mathcal{Q}_s$ is a cone over
$\mathcal{P}_s$. It is composed of strata of the form
\begin{displaymath}
\mathcal{P}_s(\delta )=\{ Q:\R^s \to \R \; ;\;  Q \mbox{ positive
definite quadratic form, }\det M_Q=\delta \},
\end{displaymath} where $\delta \in \R_+^*$. We endow each of these
strata with a metric defined as in (\ref{dist}). For each $\delta
\in \R_+^*$, any $B\in GL(s,\R)$ with $\det B=b$ induces an
isometry from $\mathcal{P}_s(\delta )$ to $\mathcal{P}_s(b^2\delta
)$. In particular, each $\mathcal{P}_s(\delta )$ is an orbit of
$\slsr$.

%%%%%%%
%%%%%%%%%%%%%
%%%%%%%%%%%%%%

The subgroup $A=\{ \mbox{diag}(e^{t_1},e^{t_2},\dots , e^{t_s})\;
;\;  t_1+t_2+\cdots + t_s=0 \}$ is a maximal $\Q $-split torus as
well as a maximal $\R $-split torus. A $\Q $-Weyl chamber (as well
as an $\R $-Weyl chamber) is $\triangleleft A= \{
\mbox{diag}(e^{t_1},e^{t_2},\dots , e^{t_s})\; ;\;  t_1+t_2+\cdots
+t_s=0,\; t_1\geq t_2\geq\cdots \geq t_s \}$.

Let $Q_0$ be the quadratic form of matrix $Id_s$. The maximal flat
$F_0=Q_0[A]$ is the set of positive definite quadratic forms $\{
\mbox{diag }(e^{t_1},e^{t_2},\dots , e^{t_s})\; ;\; t_1+t_2+\cdots
+ t_s=0 \}$. The Weyl chamber $W_0=Q_0[\triangleleft A]$ is the
subset of quadratic forms whose matrices moreover satisfy $t_1\geq
t_2\geq\cdots \geq t_s$.

\me

The dimension 1 walls (singular rays) of $W_0$, parameterized with
respect to the arc length, are the sets of quadratic forms
\begin{equation}\label{rays}
r_i= \{ \mbox{diag }(\underbrace{e^{\lambda_it},\dots
,e^{\lambda_it}}_{s-i\mbox{ times}}, \underbrace{e^{-\mu_i
t},\dots  e^{-\mu_i t}}_{i\mbox{ times}})\; ;\;  t\in \R_+ \}\; ,
\end{equation} where $\lambda_i=\sqrt{\frac{i}{s(s-i)}}$ and $\mu_i=\sqrt{\frac{s-i}{si}}\, $,
$i\in \{ 1,2,\dots s-1 \}$.

\me

The parabolic group of $r_i$ is the group
\begin{displaymath}
P(r_i)=\left\{ \left(
\begin{array}{cc}
M_1 & 0 \\
N & M_2
\end{array}
\right)\in \slsr \; ;\;  M_1\in GL(s-i,\R ),\; M_2\in GL(i,\R ),\;
N\in M_{i\times (s-i)}(\R ) \right\}\; .
\end{displaymath}

The horospherical subgroup is

\begin{displaymath}
P^0(r_i)=\left\{ \left(
\begin{array}{cc}
\epsilon M_1 & 0 \\
N & \epsilon M_2
\end{array}
\right) \; ;\;  M_1\in SL(s-i,\R ),\; M_2\in SL(i,\R ),\; \epsilon
\in \{ \pm 1\},\; N\in M_{i\times (s-i)}(\R ) \right\}\; .
\end{displaymath}

The opposite unipotent group is
\begin{equation}\label{uri}
U_+(r_i)=\left\{ \left(
\begin{array}{cc}
Id_{s-i} & N \\
0 & Id_i
\end{array}
\right) \; ;\; N\in M_{(s-i)\times i}(\R ) \right\}\; .
\end{equation}

The boundary at infinity $\partial_\infty {\mathcal{P}}_s$ can be
identified with the spherical building of flags in $\R^s$. Via this
identification, $r_1(\infty )=\langle e_s \rangle $ and more
generally $r_i(\infty )$ is the subspace $\langle e_{s-i+1},\dots
,e_s\rangle$, for $i\in \{ 1,2,\dots s-1 \}$. The spherical
chamber $W_0 (\infty )$ is identified with the flag $\langle e_s
\rangle \subset \cdots \subset \langle e_{s-i+1},\dots ,e_s\rangle
\subset \cdots \subset \langle e_{2},\dots ,e_s\rangle$.

According to Section \ref{bdry}, we can define a projection  $\slo
:
\partial_\infty {\mathcal{P}}_s \to W_0(\infty )$ and thus define
the slope of a point in $\partial_\infty {\mathcal{P}}_s$ and of a
ray in $\calp_s$. In particular a maximal singular ray $r$ has
slope $r_i(\infty )$ if and only if $r(\infty )$ is a linear
subspace of dimension $i$.

Given a flag $\mathcal{F}: V_1\subset \dots \subset V_k$ in $\R^s$
and a matrix $M\in GL(s,\R )$ we denote by $M \mathcal{F}$ the
flag $M (V_1)\subset \dots \subset M (V_k)$.

\begin{remark}\label{action}
The isometric action to the right $\Phi$ of $\slsr$ on $\calp_s$
induces the action to the right $\Phi $ on $\partial_\infty
{\mathcal{P}}_s$ identified with the spherical building of flags in
$\R^s$, defined by $\Phi (B , \mathcal{F} ) = B\iv \mathcal{F}$,
where $\mathcal{F}$ is an arbitrary flag.
\end{remark}

\subsection{The Busemann functions in the ambient
space}\label{busam}

By means of Lemma \ref{Binv} we can deduce the Busemann function
$f_{r_i}$.

\begin{lemma}\label{Bus}
Let $Q$ be a positive definite quadratic form of determinant 1 on
$\R^s$, let $Q_i$ be its restriction to $\langle e_{s-i+1},\dots
,e_s\rangle $ and let $\det Q_i$ be the determinant of $Q_i$ in
the basis $\{ e_{s-i+1},\dots ,e_s\} $. Then
$$
f_{r_i}(Q)=\sqrt{\frac{s}{(s-i)i}} \ln \det Q_i \; .
$$
\end{lemma}

%The quantity $\det Q_i$ can also be seen as follows. For a
%quadratic form $Q$ let $E_Q$ be the ellipsoid in $\R^s$ defined by
%$Q(X)\leq 1$. The determinant of the matrix of $Q$ with respect to
%a certain basis is equal to $\left( \frac{v_s}{vol\;
%E_Q}\right)^2$, where $v_s$ is the volume of the unitary ball in
%the Euclidean s-space and the volume of $E_Q$ is measured with
%respect to the given basis. It follows that $\det Q_i=\left(
%\frac{v_i}{vol\; (E_Q \cap \langle e_{s-i+1},\dots ,e_s\rangle )}
%\right)^2$.

\me

\proof According to Lemma \ref{Binv}, it is enough to prove that
the function $\Phi : {\mathcal{P}}_s\to \R,\; \Phi
(Q)=\sqrt{\frac{s}{(s-i)i}} \ln \det Q_i$, is invariant with
respect to $P^0(r_i)$ and coincides with the Busemann function on
the geodesic line ${\frak G}_i$ containing $r_i$. The second
property is obvious.

\sm

It suffices to show that the function $\Psi (Q)= \det Q_i$ is
$P^0(r_i)$-invariant. The symmetric matrix $M_Q$ of $Q$ can be written as
\begin{displaymath}
M_Q=\left(
\begin{array}{cc}
E & F \\
F^T & H
\end{array}
\right) \; ,\;  E\in M_{s-i}(\R ),\; E=E^T,\;H\in M_{i}(\R ),\;
H=H^T,\; F\in M_{(s-i)\times i}(\R ) \; .
\end{displaymath}

We have $\Psi (Q)= \det H$. Let $B\in P^0(r_i)$,

\begin{displaymath}
B= \left(
\begin{array}{cc}
\epsilon M_1 & 0 \\
N & \epsilon M_2
\end{array}
\right) \; ,\;  M_1\in SL(s-i,\R ),\; M_2\in SL(i,\R ),\;
\epsilon\in \{ \pm 1\},\; N\in M_{i\times (s-i)}(\R )\; .
\end{displaymath}

The quadratic form $Q \circ B$ restricted to $\langle
e_{s-i+1},\dots ,e_s\rangle $ has the matrix $M_2^T H M_2$ in the
basis $\left\{ e_{s-i+1},\dots ,e_s\right\}$. It follows that
$\Psi (Q \circ B ) =\det M_2^T H M_2=\det H= \Psi (Q)$.\endproof

\bi

In particular we have
$$
f_{r_1}(Q)=\sqrt{\frac{s}{s-1}}\, \ln Q(e_s) \; \mbox{ and }\;
f_{r_{s-1}}(Q)=\sqrt{\frac{s}{s-1}}\, \ln Q^*(e_1)\; ,
$$
where $Q^*$ is the ``dual quadratic form'', that is the quadratic
form whose matrix in the canonical basis is $M_Q^*$, if $M_Q$ is
the matrix of $Q$.

\begin{lemma}\label{bvector}
Let $d$ be a line in $\R^s$ and let $v$ be a non-zero vector on
$d$.
\begin{itemize}
\item[(i)] The function $f_v :
\mathcal{P}_s \to \R $,
$$
f_v(Q)=\sqrt{\frac{s}{s-1}}\, \ln Q(v)\, ,
$$ is a Busemann function of basepoint $d$.
\item[(ii)] Every Busemann function of basepoint $d$ is of the form $f_w$, where $w\in d$, $w\neq 0$.
%\item[(iii)] Let $v$ be a vector and $r$ a geodesic ray such that $f_v=f_r$.
%Then $f_{\lambda v}=f_r+ 2\sqrt{\frac{s}{s-1}}\ln \lambda$.
\end{itemize}
\end{lemma}

\proof  (i) We can write $v=Be_s$ for some $B\in \slsr$. Then
$f_v(Q)=\sqrt{\frac{s}{s-1}} \ln Q(Be_s)=f_{r_1}(\phi
(B)(Q))=f_{\phi (B)\iv r_1}(Q)$. According to Remark \ref{action},
$\phi (B\iv )r_1(\infty ) =B r_1(\infty )=B \langle e_s \rangle =
\langle v \rangle $.

(ii) Let $g$ be a Busemann function of basepoint $d$. Then $g-f_v$
is a constant function $c$. This implies that $g=f_v+c=f_w$, where
$w=e^{\frac{c}{2}\sqrt{\frac{s-1}{s}}}v$.\endproof

%(iii) is a consequence of (i).

A similar argument gives the following.

\begin{lemma}\label{bhiperpl}
Let $\hh$ be a linear hyperplane in $\R^s$ and let $v$ be a
non-zero vector orthogonal to~it.
\begin{itemize}
\item[(i)] The function $f_v^* :
\mathcal{P}_s \to \R $,
$$
f_v^*(Q)=\sqrt{\frac{s}{s-1}}\, \ln Q^*(v)\, ,
$$ is a Busemann function of basepoint $\hh$.
\item[(ii)] Every Busemann function of basepoint $\hh$ is of the form $f_w^*$, where $w\neq 0$ is orthogonal to $\hh$.
%\item[(iii)] Let $v$ be a vector and $r$ a geodesic ray such that $f_v^*=f_r$.
%Then $f_{\lambda v}^*=f_r+ 2\sqrt{\frac{s}{s-1}}\ln \lambda$.
\end{itemize}
\end{lemma}

\me

We have that $f_{r_1}=f_{e_s}$ and $f_{r_{s-1}}=f^*_{e_1}$.

\me

\noindent \textit{Notation}: Given a non-zero vector $v\in \R^n$
we denote by $Hb_v^a$ and $ H_v^a$ the horoball and horosphere
defined respectively by $f_v\leq a$ and $f_v=a$. We denote by
$Hb_{v^*}^a$ and $ H_{v^*}^a$ the horoball and horosphere defined
respectively by $f_v^*\leq a$ and $f_v^*=a$.

\me

\subsection{Totally geodesic symmetric subspaces of ${\mathcal{P}}_{s}$}\label{parab}

For details on the discussion contained in this paragraph, see
\cite[$\S I.5$]{Bo}. Let $L:\R^s \to \R $ be a non-degenerate
quadratic form of signature $(a,b)\, ,\, a+b=s$. Following
\cite[Chapter I, $\S 5$]{Bo}, we denote by $\mathcal{P}_s (L)$ the
set of positive definite quadratic forms $Q$ such that
$|L(\bar{x})|\leq Q(\bar{x})\, ,\, \forall \bar{x}\in \R^s$, and
such that $Q$ is minimal in the ordered set of positive definite
quadratic forms, verifying the previous inequality.

\begin{proposition}[\cite{Bo}, Chapter I, Proposition 5.2]\label{psl}
The following are equivalent:
\begin{itemize}
  \item[(i)] $Q\in \mathcal{P}_s (L)$ ;
  \item[(ii)] There exists a basis $\mathcal{B}$ of $\R^s$ with respect to
  which $M_Q^{\mathcal{B}}=Id_s$ and $M_L^{\mathcal{B}}=I_{a,b}$ (with the notation defined in Section \ref{not}).
\end{itemize}
\end{proposition}

\begin{cor}
If $\det M_L= \delta $ then $\mathcal{P}_s (L) \subset
\mathcal{P}_s(|\delta |)$.
\end{cor}

\proof Let $Q\in \mathcal{P}_s (L)$. By Proposition \ref{psl},
(ii), there exists $P\in GL(s, \R )$ such that $M_L= P^T I_{a,b}P$
and $M_Q=P^TP$.\endproof

\me

We consider $\mathcal{P}_s (L)$ with the metric induced from
$\mathcal{P}_s(|\delta |)$.

\begin{remarks}\label{d}
\begin{itemize}
  \item[(1)] For every $B\in GL(s,\R )$ we have
  $\mathcal{P}_s(L[B])=\mathcal{P}_s(L)[B]$.
  \item[(2)] If $L_1$ and $L_2$ are two non-degenerate quadratic
  forms of the same signature then $\mathcal{P}_s(L_1)$ and
  $\mathcal{P}_s(L_2)$ are isometric.
\end{itemize}
\end{remarks}

\proof Statement (1) is obvious. Statement (2) is a consequence of
(1) and of the discussion in Section \ref{str}.\endproof

\begin{remark}\label{totg}
If $L$ is a non-degenerate quadratic form of determinant $\delta$ then
 $\mathcal{P}_s(L)$ is a totally geodesic subspace of  $\mathcal{P}_s( |\delta |)$.
\end{remark}

\proof  By the previous remarks it suffices to prove the statement for
$\mathcal{P}_s(I_{a,b})\subset \mathcal{P}_s$. The geodesic
symmetry of $\mathcal{P}_s$ with respect to $Id_s$ is $M_Q\to
M_Q^{-1}$. The fact that $Id_s\in \mathcal{P}_s(I_{a,b})$ and that
$\mathcal{P}_s(I_{a,b})$ is invariant with respect to the previous
geodesic symmetry, together with the homogeneity of
$\mathcal{P}_s(I_{a,b})$, imply that it is a totally geodesic
subspace of $\mathcal{P}_s$.\endproof

\me

\noindent{\textit{Notation}}:  For every quadratic form $L:\R^s
\to \R$ we denote by $Con_L$ the set of vectors $\bar{x}$ in
$\R^s$ satisfying the equation $L(\bar{x})=0$.

\me

\begin{proposition}[\cite{Mo}, $\S 15, \S 16$, \cite{Wi}, $\S
4.G$]\label{flagsL}
The boundary at infinity, $\partial_{\infty }
\mathcal{P}_s(L)$, of $\mathcal{P}_s(L)$ can be identified with the
spherical building of flags of $\R^s$ composed of subspaces
totally isotropic with respect to $L$. In particular any line in
$Con_L$ is a maximal singular point in $\partial_{\infty }
\mathcal{P}_s(L)$.
\end{proposition}

\begin{remark}\label{actionL}
The action to the right of $SO_I(L)$ on $\partial_{\infty }
\mathcal{P}_s(L)$ seen as a set of flags, which corresponds to the
action to the right on $\mathcal{P}_s(L)$ as defined in this
paper, is the same as the one given in Remark \ref{action}.
\end{remark}

We study the geometry of $\mathcal{P}_s(L)$ in more detail. By
Remarks \ref{d}, it suffices to study $\mathcal{P}_s(L_0)$, where
$$
L_0= 2x_1x_s + 2x_2 x_{s-1}+\cdots + 2 x_\ell x_{s-\ell +1}+
\epsilon (x_{\ell +1}^2+\cdots + x_{s-\ell}^2)\; , \; \ell =\min
(a,b)\; , \; \epsilon \in \{\pm 1 \}\; .
$$

Let $SO_{I}(L_0)$ be the connected component of the identity of
the stabilizer of $L_0$. A maximal torus in it is \cite[$\S
11.16$]{Bo}
$$
T=\{ \mbox{diag} (e^{-t_1},\dots ,e^{-t_\ell},\underbrace{1,\dots
,1}_{s-2\ell},e^{t_\ell},\dots ,e^{t_1}) \; ;\;  (t_1,\dots
,t_\ell ) \in \R^\ell \}\; ,
$$ and a Weyl chamber is
$$
W=\{ \mbox{diag} (e^{-t_1},\dots ,e^{-t_\ell},\underbrace{1,\dots
,1}_{s-2\ell},e^{t_\ell}\dots e^{t_1}) \; ;\;  t_1 \geq t_2 \geq
\dots \geq t_\ell \}\; .
$$

Consider the one-parameter group
$\mathcal{A}=({\mathbf{a}}_t)_{t\in \R}\, $, with
\begin{equation}\label{at}
{\mathbf{a}}_t = \mbox{diag}\left(e^{-t/2\sqrt{2}}, 1, \dots , 1,
e^{t/2\sqrt{2}}\right)\, ,
\end{equation} and its sub-semigroup
${\mathcal{A}}^+=({\mathbf{a}}_t)_{t\geq 0}$.
%The factor $\frac{1}{2\sqrt{2}}$ is introduced so that the orbit of $Id_s$
%parameterized by $t\in \R $ is a unitary geodesic.
Let $r$ be the geodesic ray defined by $r(t) = {\mathbf{a}}_t^T
{\mathbf{a}}_t,\, \forall t\geq 0$. The parabolic group $P=P(r)$
writes as
\begin{displaymath}
P=\left\{ \left(
\begin{array}{ccc}
a & -a \left( M_{L_0'}X^{-1}\bar{b} \right)^T & -\frac{a}{2} L_0'(\bar{b}) \\
0 & X & \bar{b} \\
0 & 0 & a\iv\\
\end{array}
\right)\in SO_{I}(L_0) \; ;\; a\in \R^*\; ,\; X\in SO(L_0' )\; ,\;
\bar{b}\in \R^{s-2} \right\}\; ,
\end{displaymath}
 where $L_0' :\R^{s-2} \to \R \; ,\; L_0' (x_2 ,\dots , x_{s-1})= 2x_2x_{s-1} +\cdots +2 x_\ell x_{s-\ell +1} + \epsilon (x_{\ell+1}^2 + \cdots + x_{s-\ell
 }^2)$. We note that $P=\{ g\in SO_{I}(L_0) \; ;\;  g\iv (\R e_1)=\R e_1
 \}$. For this reason we also denote it by~$P_{e_1}$.

 The horospherical subgroup of $r$ is
\begin{displaymath}
P_{e_1}^0=\left\{ \left(
\begin{array}{ccc}
\epsilon & -\epsilon \left( M_{L_0'}X^{-1}\bar{b} \right)^T & -\frac{\epsilon}{2} L_0'(\bar{b}) \\
0 & X & \bar{b} \\
0 & 0 & \epsilon\\
\end{array}
\right)\in SO_{I}(L_0) \; ;\; \epsilon \in \{ \pm 1\} \; ,\; X\in
SO(L_0' )\; ,\; \bar{b}\in \R^{s-2} \right\}\; ,
\end{displaymath} and the unipotent subgroup of $r$ is
\begin{displaymath}
U_{e_1}=\left\{ \left(
\begin{array}{ccc}
1 & -\left( M_{L_0'}\bar{b} \right)^T & -\frac{1}{2} L_0'(\bar{b}) \\
0 & Id_{s-2} & \bar{b} \\
0 & 0 & 1\\
\end{array}
\right) \; ;\; \bar{b}\in \R^{s-2}
\right\}\; .
\end{displaymath}

We call a geodesic ray $\rho $ in the orbit $rSO_I(L_0)$
{\it{maximal singular ray of type }}$\wp$. The parabolic group
corresponding to it, $P(\rho )$, can be written as $\{ g\in
SO_{I}(L_0)\; ;\;  g\iv (\R v) = \R v\}$, where $v\in
Con_{L_0}\setminus \{ 0 \}$. It follows that, with the
identification of Proposition \ref{flagsL}, $\rho (\infty ) = \R
v$, that is a point in $\proj Con_{L_0}$. Whence the notation
$\wp$, coming from ``point''. We extend in the natural way the
notion of maximal singular ray of type $\wp$ to the general case
of a non-degenerate quadratic form $L$.

\medskip

\Notat Let $d$ be a line in $Con_L$. We denote by $P(d)$ the
parabolic group corresponding to $d$ seen as a point in $\proj
Con_{L}$. We denote by $U(d)$ the unipotent radical of $P(d)$.

\begin{lemma}\label{geo}
Let $d_1$ and $d_2$ be two lines in $Con_L$. If $b_L(d_1,d_2)\neq
0$ then $d_1$ and $d_2$, seen as maximal singular points in
$\partial_{\infty }\mathcal{P}_s(L)$, are opposite.
\end{lemma}

\proof We show that if $b_L(d_1,d_2)\neq 0$ then there exists a
maximal singular geodesic $\frak G$ such that ${\frak{G}}(-\infty
) =d_1$ and ${\frak{G}}(+\infty ) =d_2$. Let $v\in d_1$ and $w\in
d_2$ be two vectors such that $b_L(v,w)=1$. We consider $V= \ker
b_L(v, \cdot ) \cap \ker b_L(w, \cdot )$ of dimension $s-2$. By
the general theory of non-degenerate quadratic forms (see
\cite{Be}) we may choose a basis $\{ w_1,\dots w_{s-2} \}$ of $V$
such that in the coordinates with respect to the basis
$\mathcal{B}=\{ v,w_1,\dots w_{s-2},w \}$ the form $L$ writes as
\begin{equation}\label{l0}
L= 2x_1x_s + 2x_2 x_{s-1}+\cdots + 2 x_\ell x_{s-\ell +1}+
\epsilon (x_{\ell +1}^2+\cdots + x_{s-\ell}^2)\; , \; \ell =\min
(a,b)\; , \; \epsilon \in \{\pm 1 \}\; .
\end{equation}

The geodesic ${\frak G}(t)=Q_t$ with
$M^{\mathcal{B}}_{Q_t}=\mbox{diag }(e^{\frac{t}{\sqrt{2}}},1,\dots
,1,e^{-\frac{t}{\sqrt{2}}})$ satisfies ${\frak G}(-\infty)=d_1$
and ${\frak G}(+\infty)=d_2$.\endproof

\begin{remark}\label{ud}
Let $d_1$ and $d_2$ be two lines in $Con_L$ such that
$b_L(d_1,d_2)\neq 0$ and let $\mathcal{H}_1$ be the hyperplane
defined by $b_L(d_1, \cdot )=0$.
%Lemma \ref{geo} implies that $d_1$ and $d_2$ are opposite. Consequently $d_2 U(d_1)$ is the set
%of points in $\partial_{\infty }\mathcal{P}_s(L)$ opposite to $d_1$.
The map $$\begin{array}{ccc}
  U(d_1) & \to & \proj \left( Con_{L} \setminus \mathcal{H}_1 \right) \\
  \mathbf{u} & \mapsto & \mathbf{u}\iv d_2 \\
\end{array}$$
 is a bijection.
\end{remark}

\proof  By Remarks \ref{d} and the argument in the proof of Lemma
\ref{geo}, we may suppose that $L=L_0$, $d_1=\R e_1$, $d_2=\R e_s$
and $\mathcal{H}_1=\ker e_s^*$. It follows that $U(d_1)=U_{e_1}$.
Consider a unipotent element in $U_{e_1}$,
\begin{displaymath}
\mathbf{u}(\bar{b})= \left(
\begin{array}{ccc}
1 & -\left( M_{L_0'}\bar{b} \right)^T & -\frac{1}{2} L_0'(\bar{b}) \\
0 & Id_{s-2} & \bar{b} \\
0 & 0 & 1\\
\end{array}
\right).
\end{displaymath}

To it corresponds the line $d(\bar{b})$ in $Con_L\setminus
\mathcal{H}_1$ containing the vector $\left( -\frac{1}{2}
L_0'(\bar{b}),-\bar{b},1 \right)$ and this establishes a bijection
between $U_{e_1}$ and the set of lines in $Con_L\setminus
\mathcal{H}_1$.\endproof

\me

\noindent \textit{Notation}: For every line $d\in \proj \left(
Con_{L} \setminus \mathcal{H}_1 \right)$, we denote by
$\mathbf{u}_d$ the unipotent in $U(d_1)$ corresponding to it by
the previous bijection.

\subsection{Horoballs in $\mathcal{P}_s(L)$ and counting result}\label{gho}

\begin{lemma}\label{bup}
Let $d$ be a line in $Con_L$ and let $v$ be a non-zero vector on
$d$.
\begin{itemize}
\item[(i)] The function $f_v :
\mathcal{P}_s(L) \to \R $,
$$
f_v(Q)=\sqrt{2} \ln Q(v)\, ,
$$ is a Busemann function of basepoint $d$.
\item[(ii)] Every Busemann function of basepoint $d$ is of the form $f_w$, where $w\in d$, $w\neq 0$.
\item[(iii)] Let $v$ be a vector and $r$ a geodesic ray such that $f_v=f_r$. Then $f_{\lambda v}=f_r+ 2\sqrt{2}\ln
\lambda$.
\end{itemize}
\end{lemma}

\proof  (i) Since $L$ is non-degenerate, there exists $w\in Con_L$
such that $b_L(v,w)=1$. We proceed as in the proof of Lemma
\ref{geo} and consider a basis $\mathcal{B}$ with $v$ and $w$ the
first and respectively last vector, with respect to which $L$ can
be written as in (\ref{l0}). We consider the geodesic $\frak G$
joining $d$ and $\R w$, ${\frak G}(t) = (\Pi ^{-1})^T \mbox{diag
}(e^{\frac{t}{\sqrt{2}}},1,\dots ,1, e^{-\frac{t}{\sqrt{2}}})\Pi
^{-1}$, where $\Pi $ is the matrix having the vectors of
$\mathcal{B}$ as columns. The geodesic ray corresponding to $d$ is
$r(t)= {\frak G}(-t)\, ,\, t\geq 0$, and its horospherical
subgroup is $\Pi P_{e_1}^0 \Pi^{-1}$. For every $p\in \Pi
P_{e_1}^0 \Pi^{-1}$, $f_v ( Q[p] )= \sqrt{2} \ln Q(p(v))= \sqrt{2}
\ln Q(v)$, since $p(v)=v$. Also, $f_v({\frak G}(t))= t$. Lemma
\ref{Binv} allows to conclude.

(ii) is proved as in Lemma \ref{bvector} and (iii) follows
immediately from the formula of $f_v$.\endproof

\me

\Notat For every $v\in Con_L \setminus \{ 0\}$ we denote by
$H_v^a$, $Hb_v^a$ and $Hbo_v^a$ the horosphere defined by $f_v=a$,
and the horoball and open horoball defined by $f_v\leq a$ and
$f_v< a$, respectively. For $a=0$ we simply write $H_v$, $Hb_v$
and $Hbo_v$.

\me

\begin{lemma}
Let $v,w$ be two vectors in $Con_L$ such that $b_L(v,w)\neq 0$.
The oriented distance between the horospheres $Hb_v$ and $Hb_w$ is
$2\sqrt{2}\ln |b_L(v,w)| $.
\end{lemma}

\proof Let $w_1=\frac{1}{\chi } w$ with $\chi = b_L(v,w)$. The
Busemann function $f_{w}$ is equal to $f_{w_1} + 2\sqrt{2}\ln
|b_L(v,w)| $. Therefore it suffices to prove the statement of the
lemma when $b_L(v,w)=1$. By Remark \ref{d} we may suppose that
$L=L_0$. Moreover, by Witt Theorem (\cite{Be}, \cite[$\S
4.G$]{Wi}), we may suppose that $v=e_1$ and $w=e_s$. A geodesic
joining $\R e_1$ and $\R e_s$ is ${\frak G}(t)= \mbox{diag
}(e^{-\frac{t}{\sqrt{2}}}, 1,\dots ,1, e^{\frac{t}{\sqrt{2}}} )$.
We have that ${\frak G} \cap H_{e_1}={\frak G} \cap H_{e_s}=
\{Id_s \}$, which finishes the proof. \endproof

\begin{cor}\label{norm}
Let $v_0\in Con_L\setminus \{ 0\}$ be fixed and let
$\mathcal{H}_0$ be the hyperplane defined by $b_L(v_0, \cdot )=0$.
For every compact set $K$ in $\proj (Con_L \setminus
\mathcal{H}_0)$, we have that $|{\mathrm{odist}}\,
(Hb_{v_0},Hb_{w}) - 2\sqrt{2} \ln \| w \|\; |$ is bounded
uniformly in $w\in \R K$.
\end{cor}

\me

Corollary \ref{norm} and the counting result Corollary \ref{lSu}
give the following.

\begin{proposition}\label{counting1}
 Let $\Gamma$ be an irreducible lattice in
    $SO_{I}(L)$ and let $\bar{r}$ be a maximal singular cusp ray in
    $\mathcal{P}_s(L)/\Gamma $ such that if $r$ is a lift of
    it in $\mathcal{P}_s(L)$, then $r$ is of type $\wp$.
     Let $r(\infty )=d \in \proj Con_L$ and let $v$ be a non-zero vector on
    $d$.
    Let $\Omega$ be a relatively compact open subset of $\proj
    Con_L$ such that its closure does not intersect $\proj \ker b_L(v_0,\cdot
    )$ for some $v_0\in Con_L \setminus \{ 0\}$. Let $a>1$. For every open subset
$\oo$ of $\Omega$ we denote by $N(k\, ;\, \oo )$  the cardinal of
the set of vectors $$\left\{ v \gamma \; ;\; \gamma \in \Gamma \:
,\: \R v \gamma\in \oo \: ,\:\| v\gamma \| \in
[a^k,a^{k+1})\right\}\, .$$

For any $a\geq a_0(L,\Omega )$ we have that
$$
\mathbf{K}_1a^{(k+1)(s-2)}\nu (\oo )\leq N(k\, ;\, \oo )\leq
\mathbf{K}_2 a^{(k+1)(s-2)}\nu (\oo )\, ,\; \mbox{ for every }
k\geq k_0(\oo , \Omega , v)\, ,
$$ where $\nu $ is the canonical measure on
$\proj Con_L$ and $\mathbf{K}_i=\mathbf{K}_i(L,\Gamma , \Omega )$.
\end{proposition}

\me

\begin{lemma}\label{uang}
Let $d_0$ be a fixed line in $Con_L$ and let $\mathcal{H}_0$ be
the hyperplane defined by $b_L(d_0, \cdot )=0$. With the notation
following Remark \ref{ud}, on every compact subset $\mathcal{K}$
of $\proj (Con_L \setminus \mathcal{H}_0)$ the angle between two
lines $d_1$ and $d_2$ in $\mathcal{K}$ is bi-Lipschitz equivalent
to the distance between $\mathbf{u}_{d_1}$ and $\mathbf{u}_{d_2}$
in $U(d_0)$.
\end{lemma}

\proof Up to isometry, we can reduce to the case when $L=L_0$,
$d_0=\R e_1$ and $\mathcal{H}_0=\ker e_s^*$. In this case $U(d_0)
= U_{e_1}$. The Riemannian distance on $U_{e_1}$ coming from the
Lie group isomorphism
$$
\R^{s-2}  \to  U_{e_1}
$$
$$
\bar{b} \mapsto \mathbf{u}(\bar{b})
$$ is invariant. With the notation in the proof of Remark \ref{ud},
 the angle between two lines
 $d(\bar{b})$ and  $d(\bar{b}')$ is bi-Lipschitz equivalent to
  $\|\bar{b} -\bar{b}'\|_e$, if $\bar{b}$ and $\bar{b}'$ are in a
  compact set of $\R^{s-2}$.\endproof

\subsection{Traces of horoballs on unipotent orbits}\label{trc}

Throughout this section we fix $v_0$ a vector in $Con_L\setminus
\{ 0\}$ and a geodesic $\frak G$ in $\calp_s(L)$ with ${\frak
G}(-\infty )=d_0$, where $d_0=\R v_0$, and ${\frak G}(0)$ in
$H_v$. Let $\mathcal{H}_0$ be the hyperplane $\ker b_L(v_0,\cdot
)$, and let $P(d_0)$ be the parabolic group corresponding to
$d_0$. This group has a Langlands decomposition, $P(d_0)=MAU$ such
that $\frak G$ is an orbit of $A$. We denote by $r$ the geodesic
ray ${\frak G}|_{[0,+\infty)}$.

Let $w$ be an arbitrary vector in $ Con_L\setminus \mathcal{H}_0$
and let $D$ be the oriented distance between the horoballs
$Hb_{v_0}$ and $Hb_w$. We wish to study the trace of the horoball
$Hb_w$ on $U$ identified with the orbit of $r(t)$ under $U$, that is
the set
\begin{equation}\label{trt}
Tr_t(w)=\{ \mathbf{u}\in U \; ;\;  r(t)[\mathbf{u}]\in Hb_w\}\; .
\end{equation}

We put the condition $t=D+\tau \, ,\, \tau \geq 0$, otherwise
$Tr_t(w)$ is empty. We have the following lemma.

\begin{lemma}\label{bilevoltr}
 We consider the unipotent group $U$ endowed with an
invariant metric and with the Haar measure $\nu$. For any vector
$w$ in $Con_L\setminus \hh_0$ the following holds. Let
$D=\mathrm{odist}(Hb_{v_0},Hb_w)$.
\begin{itemize}
  \item[(1)] For any $\tau \geq 0$
  \begin{equation}\label{trbile}
Tr_{D+\tau } (w) \subset B\left(\mathbf{u}_w,\kappa_0 \,
e^{-\frac{D}{2\sqrt{2}}}\right)\; ,
\end{equation}
 where $\mathbf{u}_w\in U$ is such that ${\frak G}[\mathbf{u}_w](+\infty )= \R
 w$, and $\kappa_0$ is a constant depending on $v_0$ and on the metric chosen on $U$.
  \item[(2)] If $L$ has signature $(a,b)$ with $\min(a,b)\geq 2$
  then
\begin{displaymath}
C_1 e^{-\frac{D(s-2)+\tau }{2\sqrt{2}}}\leq \nu (Tr_{D+\tau }(w))
\leq C_2 e^{-\frac{D(s-2)+\tau }{2\sqrt{2}}}\, ,
\end{displaymath}
 where $C_1$ and $C_2$ are constants depending on $v_0$.
\end{itemize}
\end{lemma}

\proof \textit{Step} 1. First we consider the particular case when
$L=L_0$, $v_0=e_1$ and ${\frak G}_0(t)=Q_t=\mathrm{diag}\, \left(
e^{\frac{t}{\sqrt{2}}},1,\dots ,1,e^{-\frac{t}{\sqrt{2}}}\right)$.
 Via the isometry $\bar{b}\mapsto \mathbf{u}(\bar{b})$ we identify
 $U_{e_1}$ to $\R^{s-2}$, its Haar measure $\nu $ is the Lebesgue
 measure, and we choose as invariant metric the Euclidean metric.
 We take the vector $w$ to be $\lambda e_s$, with $\lambda =
e^{\frac{D}{2\sqrt{2}}}$, for an arbitrary $D\geq 0$. Let $r(t)
=Q_t$ for $t\geq 0$.

(1) We have
$$
Tr_t (\lambda e_s)= \left\{ \mathbf{u}\in U_{e_1}\; ;\;  Q_t
[\mathbf{u}] \in Hb_{\lambda e_s} \right\}= \left\{ \mathbf{u}\in
U_{e_1}\; ;\; Q_t [\mathbf{u}](e_s) \leq \frac{1}{\lambda^2}
\right\}=
$$
$$
\left\{ \mathbf{u}(\bar{b})\in U_{e_1}\; ;\;
\frac{1}{4}e^{\frac{t}{\sqrt{2}}} [L_0'(\bar{b})]^2 + \|
\bar{b}\|_e^2 + e^{-\frac{t}{\sqrt{2}}}\leq
e^{-\frac{D}{\sqrt{2}}} \right\}=
$$
$$
\left\{ \mathbf{u}(\bar{b})\in U_{e_1}\; ;\;
\frac{1}{4}e^{\frac{\tau }{\sqrt{2}}} \left[L_0'\left(
e^{\frac{D}{2\sqrt{2}}}\bar{b}\right) \right]^2 + \left\|
e^{\frac{D}{2\sqrt{2}}} \bar{b}\right\|_e^2 + e^{-\frac{\tau
}{\sqrt{2}}}\leq 1 \right\}\; .
$$

We denote by $\mathcal{O}_\alpha $ the homothety in $\R^{s-2}$ of
center the origin and of factor $\alpha $. We identify
$Tr_t(\lambda e_s)$ with $\mathcal{O}_\eta (\mathcal{M}_\tau )$,
where $\eta = e^{-\frac{D}{2\sqrt{2}}}$ and
$$
\mathcal{M}_\tau = \left\{ \bar{b}\in \R^{s-2}\; ;\;
\frac{1}{4}e^{\frac{\tau }{\sqrt{2}}} [L_0'(\bar{b})]^2 + \|
\bar{b}\|_e^2 + e^{-\frac{\tau }{\sqrt{2}}}\leq 1 \right\}\; .
$$

The set $\mathcal{M}_\tau $ is quasi-conformal to the set
$$
\mathcal{M}_\tau' = \left\{ \bar{b}\in \R^{s-2}\; ;\;  \|
\bar{b}\|_e^2 \leq 1-e^{-\frac{\tau }{\sqrt{2}}}\, ,\,
\frac{1}{4}e^{\frac{\tau }{\sqrt{2}}} [L_0'(\bar{b})]^2 \leq
1-e^{-\frac{\tau }{\sqrt{2}}} \right\}\; ,
$$ in the sense that
 $\mathcal{O}_{1/\sqrt{2}}(\mathcal{M}_\tau ') \subset \mathcal{M}_\tau  \subset \mathcal{M}_\tau
 '$. We have that $\mathcal{M}_\tau '\subset B(0,1)$. We conclude that
\begin{displaymath}
Tr_t(\lambda e_s) \subset B\left(
0,e^{-\frac{D}{2\sqrt{2}}}\right)\, .
\end{displaymath}

(2) We can rewrite $\mathcal{M}_\tau' $ as $\mathcal{O}_\chi\left(
\mathcal{M}_\tau''\right)$, where $\chi = \sqrt{2} \left(
1-e^{-\frac{\tau }{\sqrt{2}}}
\right)^{\frac{1}{4}}e^{-\frac{\tau}{4\sqrt{2}}}$ and
\begin{displaymath}
\mathcal{M}_\tau'' = \left\{ \bar{b}\in \R^{s-2}\; ;\;
|L_0'(\bar{b})| \leq 1\, ,\, \| \bar{b} \|_e\leq
\frac{1}{\sqrt{2}} \left( 1-e^{-\frac{\tau }{\sqrt{2}}}
\right)^{\frac{1}{4}} e^{\frac{\tau}{4\sqrt{2}}} \right\}\; .
\end{displaymath} Lemma 3.8 from \cite{EMM} implies that
$$
\nu \left( \mathcal{M}_\tau'' \right) \sim C' e^{\frac{\tau
(s-4)}{4\sqrt{2}}} \mbox{ as } \tau \to \infty \; ,
$$ where $C'$ is an absolute constant. Hence
$$
\nu \left( \mathcal{M}_\tau'\right) \sim C e^{-\frac{\tau
}{2\sqrt{2}}} \mbox{ as } \tau \to \infty \; ,
$$ and
$$
C_1 e^{-\frac{\tau }{2\sqrt{2}}}\leq \nu \left(\mathcal{M}_\tau
\right)\leq C_2 e^{-\frac{\tau }{2\sqrt{2}}}\, ,
$$ where $C_1$ and $C_2$ are universal constants. This and the fact that
 $Tr_{D+\tau }(\lambda e_s)$ is isometric to
 $\mathcal{O}_{e^{-D/2\sqrt{2}}} \mathcal{M}_\tau$
 yields the conclusion.

\me

\noindent \textit{Step} 2. We place ourselves in the general case.
There exists $B\in GL(n,\R )$ such that $\Phi (B)(L_0)=L$. Remark
\ref{d} implies that $\Phi (B)$ is an isometry between
$\calp_s(L_0)$ and $\calp_s(L)$. The fact that $SO_{I}(L)$ acts
transitively on geodesics with both points at infinity lines in
$Con_L$ implies that we may suppose that $\Phi (B) \left( {\frak
G}_0 \right) ={\frak G}$. Let $\mathbf{u}_w\in U$ be such that
${\frak G}[\mathbf{u}_w](+\infty )= \R w$. Since $\mathbf{u}_w$
acts by isometry on $U$ to the right, it suffices to prove the
result in the particular case when $\mathbf{u}_w=id$.

We then have that $B\iv U_{e_1}B=U$ and that $B\iv Tr_t(\lambda
e_s )B=Tr_t(w)$. The conjugation by $B$ transforms the Haar
measure on $U_{e_1}$ into the Haar measure on $U$ and the
Euclidean metric on $U_{e_1}$ into an invariant metric on $U$,
bi-Lipschitz equivalent to the one that was chosen.\endproof

\subsection{Quotient spaces, equidistribution of rational vectors}\label{bparab}

Let $\calp_s$ be the ambient symmetric space defined in Section
\ref{str}. Let $\Gamma$ be the lattice $\slsz$. Its $\Q$-rank
$\mathbf{r}$ is equal to the $\R$-rank of $\slsr$ and with $s-1$.

\me

\noindent \textit{Notation}: We denote by $\mathcal{T}_s$ the
quotient space ${\mathcal{P}}_{s}/ \Gamma$. In accordance with the
notation introduced in Section \ref{lssp}, we denote by $\pr $
   the projection of ${\mathcal{P}}_{s}$ onto $\mathcal{T}_s$.

\me

The projection $\overline{W}_0=\pr (W_0)$ is isometric to $W_0$.
Moreover, $\mathcal{T}_s$ is at finite Hausdorff distance of
$\overline{W}_0$. We denote by $\bar{r}_i$ the projection of the
ray $r_i$ defined in (\ref{rays}). According to Remark
\ref{horob}, (1), and to Lemma \ref{bvector}, for $a<0$ with $|a|$
large enough, the projection of $Hb_{e_s}^a$ is $Hb_a (\bar{r}_1)$
and its pre-image is $\bigcup_{v\in \pri^s} Hb_v^a$. Likewise
$Hb_{e_1^*}^a$ projects onto $Hb_a (\bar{r}_{s-1})$ and its
pre-image is $\bigcup_{v\in \pri^s} Hb_{v^*}^a$.

\me

Let $L$ be a non-degenerate rational quadratic form on $\R^s$. The
group $SO_\Z(L)$ is a lattice, which we denote by $\Gamma_L$. It
has $\Q$-rank $\mathbf{r}$ equal to the dimension of the maximal
rational linear subspace totally isotropic with respect to $L$
(that is, contained in $Con_L$).

\me

\noindent \textit{Notation}: We denote by $\calv_L$ the quotient
space ${\mathcal{P}}_{s}( L)/ \Gamma_L$. We denote by $\pr_L$
   the projection of ${\mathcal{P}}_{s}( L)$ onto $\calv_L$.

\me

The manifold $\calv_L$ is a locally symmetric space of finite
volume, at finite Hausdorff distance of a finite union of
Euclidean sectors of dimension $\mathbf{r}$, which are projections
of $\Q$-Weyl chambers in ${\mathcal{P}}_{s}( L)$ (\cite{BoS},
\cite{Le}). Let $\bar{r}_1,\, \bar{r}_2, \dots ,\, \bar{r}_k$ be
all the maximal singular cusp rays in $\calv_L$ such that their
lifts $r_1,\, r_2,\, \dots ,\, r_k$ in ${\mathcal{P}}_{s}( L )$
are of type $\wp \,$.
  %For $\tau_i $ sufficiently large,
  %$\mathcal{V}_L$ can be written as $ K \cup \bigsqcup_{i=1}^k \mathcal{V}_{L ,i}^{\tau_i}$, where $K$ is a compact subset and the subset
%$\mathcal{V}_{L,i}^{\tau_i} $ is the projection of $Hb_{-\tau_i}(\rho_i)$ in $\calv_L$ (see Section \ref{ssg}).
The set $\{ \bar{r}_1 (\infty ),\, \bar{r}_2 (\infty ),\, \dots
,\, \bar{r}_k (\infty )\}$ can be identified with the quotient under
the action of $\Gamma_L$ of the set of all rational lines in
$Con_L$. The latter set can also be seen as the set $\left(
\pri^s\cap Con_L\right) /\pm 1 = \pri_+^s\cap Con_L$.
 Let $r_i(\infty )=v_i\in \pri_+^s\cap Con_L$. By the
previous considerations, $\{ v_1,v_2,\dots , v_k \}$ can also be
identified with $ \left( \pri_+^s\cap Con_L \right) / \Gamma_L$.
 By Lemma \ref{bup}, $f_{r_i}=f_{\lambda_iv_i}$, where $\lambda_i\in
(0,\infty )$.

According to Remark \ref{horob}, (1), if $a<0$ with $|a|$ large
enough then for any $i\in \{ 1,2,\dots ,k \}$,
$\pr_L\left(Hb^{a}_{v_i}\right)$ coincides with
$Hb_{a_i}(\bar{r}_i)$, for some $a_i<0$, and the pre-image of it
is $Hb^{a}_{v_i}\Gamma_L =\bigcup_{v\in v_i\Gamma_L}Hb^{a}_{v}$.
Therefore the projection of $\bigcup_{i=1}^k Hb^{a}_{v_i}$ is
$\bigcup_{i=1}^k Hb_{a_i}(\bar{r}_i)$ and it has the pre-image
$\bigcup_{v\in \pri_+^s\cap Con_L }Hb^{a}_{v}$.

%The pre-image of the set $ \pi^L_X\left(Hb^{\tau}_{\rho_i}\right)$ is
%$$
%Hb^{\tau}_{\rho_i} \Gamma_L = \{ Q\in {\mathcal{P}}_{s}( L) \; ;\;
%f_{v_i} (Q) \leq - \tau-2\sqrt{2}\ln a_i \} \Gamma_L=\bigcup_{v\in
%v_i\Gamma_L}\{ Q\in {\mathcal{P}}_{s}( L) \; ;\; f_{v} (Q) \leq -
%\lambda \}\, ,
%$$ where we
%denote $\tau+2\sqrt{2}\ln a_i$ by $\lambda$.

%This implies that the pre-image of the set $ \pi^L_X \left(
%\bigcup_{i=1}^kHb^{\tau_i}_{\rho_i}\right)$ with
%$\tau_i=\lambda-2\sqrt{2}\ln a_i$ is $$\bigcup_{v\in \pri_+^s\cap
%Con_L} \{ Q\in {\mathcal{P}}_{s}( L) \; ;\;  f_{v} (Q) \leq -
%\lambda \} = \bigcup_{v\in \pri_+^s\cap Con_L} Hb_v^{-\lambda }\,
%.$$

The application of the Proposition \ref{counting1} to each of the
rays $\bar{r}_i$ gives the following.

\begin{proposition}\label{counting2}
 Let $\Omega$ be a relatively compact open subset of $\proj
    Con_L$ such that its closure does not intersect $\proj \ker b_L(v_0,\cdot
    )$ for some $v_0\in Con_L \setminus \{ 0\}$. Let $a>1$. For every open subset
$\oo$ of $\Omega$ we denote by $N(k\, ;\, \oo )$  the cardinal of
the set of vectors $$\left\{ v \in \pri^s\cap Con_L \; ;\; \R v\in
\oo \: ,\:\| v \| \in [a^k,a^{k+1})\right\}\, .$$

For any $a\geq a_0(L, \Omega )$ we have that
$$
\mathbf{K}_1 a^{(k+1)(s-2)}\nu (\oo )\leq N(k\, ;\, \oo )\leq
\mathbf{K}_2 a^{(k+1)(s-2)}\nu (\oo )\, ,\; \mbox{ for every
}k\geq k_0(\oo , \Omega )\, ,
$$ where $\nu $ is the canonical
measure on $\proj Con_L$ and $\mathbf{K}_i=\mathbf{K}_i(L,\Omega
)$.
\end{proposition}

A consequence of this proposition is Corollary \ref{counting3} in
the introduction.

\section{Diophantine approximation on a rational
quadric}\label{elip}

\subsection{Some preliminary considerations}

Let $\q :\R^n \to \R$ be a non-degenerate quadratic form with
rational coefficients and let ${\mathfrak{Q}}_\q$ be the quadric
defined by $\q=1$. Before beginning the proof of Theorem \ref{T2},
we wish to point out that an argument
 with a projection on a rational
hyperplane does not work. This can be illustrated on the example
of ${\mathfrak{Q}}_\q =\sph^n(0\, ,\, 1)\subset \R^{n+1}$. For
simplicity we replace ${\mathfrak{Q}}_\q$ by $\sph^n(e_{n+1}\, ,\,
1)$, which we denote in what follows by $T\sph^n$. We recall that
the stereographic projection with respect to $2e_{n+1}$ is
$$
  \begin{array}{cccc}
    pr\; : & T\sph^n & \rightarrow & \R^n \\
      & \bar{x} & \mapsto &
      \frac{2}{2-x_{n+1}}\bar{x}-\frac{2x_{n+1}}{2-x_{n+1}}e_{n+1}\; .
  \end{array}
$$

Its inverse is
$$
  \begin{array}{cccc}
    inv : & \R^n & \rightarrow & T\sph^n \\
      & \bar{y} & \mapsto &
      \frac{4}{4+\| \bar{y} \|_e^2}\bar{y}+\frac{2\| \bar{y} \|_e^2}{4+\| \bar{y} \|_e^2 }e_{n+1}\; .
  \end{array}
$$

We have that

\me

$\bullet$  $pr \left( \mathcal{S}_\alpha (T\sph^n ) \right)\subset
\mathcal{S}_{\alpha-\epsilon }( \R^n )$, for any $\alpha $ and
$\epsilon $ ;

\me

$\bullet$  $inv \left( \mathcal{S}_\alpha (\R^n ) \right)\subset
\mathcal{S}_{\frac{\alpha -1}{2}-\epsilon }(T\sph^n )$, for any
$\alpha $ and $\epsilon $.

It follows that
\begin{equation}\label{projst}
inv \left( \mathcal{S}_{1+2\alpha +2\epsilon }(\R^n )
\right)\subset \mathcal{S}_{\alpha }(T\sph^n )\subset inv \left(
\mathcal{S}_{\alpha -\epsilon } (\R^n ) \right)\; .
\end{equation}

On the other hand, by Jarn\'{\i}k Theorem, $\dim_H
\mathcal{S}_{\alpha } (\R^n ) =\frac{n+1}{\alpha +1}$ for all
$\alpha \geq \frac{1}{n}$. This and relation (\ref{projst}) imply
that
\begin{equation}\label{incl}
  \frac{n+1}{2(\alpha +1)}\leq \dim_H \mathcal{S}_{\alpha }(T\sph^n
  )\leq \frac{n+1}{\alpha +1}\; .
\end{equation}

For $n=1$ we obtain $\dim_H \mathcal{S}_{\alpha }(T\sph^1
  )\geq \frac{1}{1+\alpha }$, which together with the inequality
  of Melnichuk, $\dim_H \mathcal{S}_{\alpha }(T\sph^1
  )\leq \frac{1}{1+\alpha }$ \cite{DD}, imply the result of H.
  Dickinson and M.M. Dodson. For $n>1$ both bounds given by (\ref{incl}) are not sharp.

\me

The first step in the proof of Theorem \ref{T2} is the following
Lemma.

\begin{lemma}\label{aprox}
Let $\psi$ be an approximating function such that $\lim_{x\to
\infty } x\psi (x) =0$. Let $\bar{x}\in {\mathfrak{Q}}_\q$ and let
$\frac{1}{q}\bar{p}\in \Q^n$ be such that
\begin{equation}\label{(rel)}
\left\| \bar{x}- \frac{1}{q}\bar{p} \right\| \leq \frac{\psi
(q)}{q}\; .
\end{equation}

If $q$ is large enough then $\frac{1}{q}\bar{p}\in
{\mathfrak{Q}}_\q$.
\end{lemma}

\proof We have $\q (\bar{x})=\sum_{1\leq i\leq j \leq n}
a_{ij}x_ix_j =1 $, where $\bar{x}=(x_1\, ,\, \dots ,\, x_n )$ and
$a_{ij}\in \Q , \forall i,j$. Relation (\ref{(rel)}) implies that
$qx_i = p_i + \varepsilon_i$, where $\bar{p}=(p_1\, ,\, p_2\,
,\dots ,\, p_n )$ and $\varepsilon_i=O\left( \psi (q) \right)$. It
follows that $q^2=\sum_{1\leq i\leq j \leq n} a_{ij}qx_i\,
qx_j=\sum_{1\leq i\leq j \leq n}
a_{ij}(p_i+\varepsilon_i)(p_j+\varepsilon_j)=\q
(\bar{p})+\sum_{1\leq i\leq j \leq n} a_{ij}(p_i \varepsilon_j+p_j
\varepsilon_i) +\sum_{1\leq i\leq j \leq n} a_{ij}\varepsilon_i
\varepsilon_j$. We have that $\sum_{1\leq i\leq j \leq n}
a_{ij}(p_i \varepsilon_j+p_j \varepsilon_i)=O\left( q\psi (q)
\right)$ and $\sum_{1\leq i\leq j \leq n} a_{ij}\varepsilon_i
\varepsilon_j =O\left( \psi (q)^2 \right)$. Both sums tend to $0$
when $q\to \infty$. Since $\q $ has rational coefficients, $\q
(\Z^n )\subset \frac{1}{N} \Z$ for some $N\in \N $. It follows
that $q^2 -\q (\bar{p})\in \frac{1}{N} \Z $. On the other hand,
for $q$ large enough $| q^2 -\q (\bar{p}) | < \frac{1}{N}$. We
conclude that $\q \left( \frac{1}{q}\bar{p} \right) =1$ for $q$
large enough.
\endproof

\me

The previous lemma implies in particular that if
${\mathfrak{Q}}_\q \cap \Q^n = \emptyset $ then $\mathcal{S}_\psi
( {\mathfrak{Q}}_\q)= \emptyset $ for $\psi $ such that
$\lim_{x\to \infty }x\psi (x)=0$. In what follows we work under
the hypothesis that ${\mathfrak{Q}}_\q \cap \Q^n \neq \emptyset $.

\me

%Another key lemma is the following.

%\begin{lemma}\label{c} Let $M$ be a submanifold in $\R^n$ and let
%$$\mathcal{S}^\kappa_\psi (M) =\left\{ \bar{x}\in M \; ;\; \left\|\bar{x}-\frac{1}{q}\bar{p} \right\| \leq \kappa \psi (q) \mbox{has infinitely many solutions } \frac{1}{q}\bar{p}\in \Q^n\right\}\; ,$$ where $\kappa >0$.
% We denote $\mathcal{S}^1_\psi(M)$ by $\mathcal{S}_\psi (M)$ as previously. If there exists
%$\kappa_0$ and a continuous function $f:]1,\infty [ \to [0,\infty[ $ such that
%$\dim_H \mathcal{S}^{\kappa_0}_\psi (M)=f(\alpha),\; \forall \alpha>1$, then for any
%$\kappa >0$, $\dim_H\mathcal{S}^{\kappa}_\alpha (M)=f(\alpha ),\; \forall \alpha>1$.
%\end{lemma}

%\proof It is an immediate consequence of the fact that for any $0<
%\kappa \leq \kappa'$ and for any $\varepsilon >0$,
%$$\mathcal{S}^{\kappa'}_{\alpha +\varepsilon }(M) \subset\mathcal{S}^\kappa_\alpha (M)\subset \mathcal{S}^{\kappa '}_\alpha(M)\; .$$

%\endproof

\subsection{Generalized notion of Hausdorff measure, Hausdorff dimension}\label{ghm}

\begin{definitions}\label{dimf}
\begin{itemize}
    \item[(1)] A \textit{dimension function} is a function $\varphi : \R_+ \to
\R_+$ which is increasing and continuous and such that $\lim_{x\to
0}\varphi (x) =0$.
    \item[(2)] We say that  a dimension function $\varphi$\textit{ dominates the function }$x^\delta $
     for some $\delta >0$, if
\begin{itemize}
    \item $x\mapsto x^{-\delta } \varphi (x)$ is a decreasing
    function;
    \item $\lim_{x\to 0} x^{-\delta } \varphi (x) = \infty$.
\end{itemize}
\end{itemize}
\end{definitions}

\me

Let $(\mm , d)$ be a metric space and $F$ a non-empty subset of
$\mm$. For $\epsilon
>0$ we call $\epsilon$-\textit{cover of }$F$ a countable collection of
balls $\{ B_j\}_{j\in J}$ of radii $r_j$ at most $\epsilon $ for
every $j\in J$, such that $F\subset \bigcup_{j\in J} B_j$. Define
$$
\hh^\varphi_\epsilon (F) = \inf \left\{ \sum_{j\in J}\varphi
(r_j)\: :\: \{ B_j\}_{j\in J}\: \mbox{ an } \epsilon-\mbox{cover
of }F\right\}\, .
$$

The \textit{Hausdorff measure of $F$ with respect to the dimension
function }$\varphi$ is defined by
$$
\hh^\varphi (F) = \lim_{\epsilon \to 0 } \hh^\varphi_\epsilon
(F)=\sup_{\epsilon > 0 } \hh^\varphi_\epsilon (F)\, .
$$

\begin{remark}
%If $\psi : \R_+ \to \R_+$ is such that $\lim_{x\to 0}\frac{\psi (x)}{\varphi (x)}=0$ then $\hh^\varphi (F) < \infty$ implies that $\hh^\psi (F) =0$.
When $\varphi (x)=x^s$, $\hh^\varphi$ becomes the usual Hausdorff
measure $\hh^s$. Recall that one defines the Hausdorff dimension
$\dim_H F$ by
    $$
\dim_H F = \inf \{ s\geq 0 \; ;\;  \hh^s (F) =0 \}=\sup \{ s\geq 0
\; ;\;  \hh^s (F) =\infty \}\, .
    $$
\end{remark}

\subsection{Ubiquitous systems}\label{us}

We need the notion of ubiquitous system as introduced in
\cite{BDV2}, and one of its main properties, which we recall in
the sequel. We shall not need the notion in its most general form,
as the resonant sets we work with are points. See \cite{BDV2} for
a more general and detailed presentation, as well as for proofs.
We note that in the particular case when the compact metric space
considered below $(\mm , \dist)$ is a bounded subset of
 $\R^n$ with the Euclidean metric, and $m$ is the Lebesgue
 measure, the notion of ubiquitous system and a weaker version of Theorem
 \ref{minf} were formulated in \cite{DRV2} (see also \cite{BD}). In the
 case when the resonant sets are points this notion
 coincides with the notion of regular system from \cite{BSch}. Also, a variant of the notion has been defined and used in
 \cite{Bu}.

Let $(\mathcal{M} , \dist , m )$ be a compact metric space with a
probability measure. Assume that the measure $m$ satisfies the
following condition.
\begin{center}
$(M)$ There exists $\delta >0$ and $R_0>0$ such that
    for any $x\in \mm$ and $R\leq R_0$,
    $$
a R^\delta \leq m(B(x,R))\leq bR^\delta \, .
    $$

    The constants $a$ and $b$ are independent of the ball and can
    be assumed to satisfy $$0<a<1<b \, .$$
\end{center}

\begin{remark}
The condition $(M)$ implies that the Hausdorff dimension $\dim_H
\mm =\delta $.
\end{remark}

Let $I$ be an infinite countable family of indices and let $\varpi
:I \to \R_+$ be a weight function on it. Assume that for every
$M>0$, the set $\{ i\in I\; ;\;  \varpi (i)\leq M \}$ is finite.
Let $\Re =\{ p_i \; ;\;  i\in I \}$ be a collection of points in
$\mm$, called \textit{resonant points}.

Let $\rho :\R_+ \to \R_+$ be a function such that $\lim_{x\to
\infty }\rho (x)=0$. It will be called \textit{the ubiquitous
function}. Let $u=(u_n)_{n\in \N}$ be an increasing sequence of
positive real numbers such that $\lim_{n\to \infty} u_n =\infty \,
.$ We assume that the function $\rho $ is $u$-\textit{regular},
that is there
 exists a constant $0< \lambda < 1$ such that for $n\in \N$
 sufficiently large,
 $$
\rho (u_{n+1})\leq \lambda \rho (u_{n})\, .
 $$

The pair $(\Re , \varpi )$ is said to be \textit{a local
$m$-ubiquitous system relative to $(\rho , u ) $} if the following
condition is satisfied. There exists $R_1 >0$ such that an
arbitrary ball $B$ in $\mm$ of radius $R\leq R_1$ satisfies
$$
m\left( B\cap  \bigcup_{\varpi (i)\leq u_n} B\left( p_i\, ,\, \rho
\left( u_n\right)\right) \right)\geq \kappa\,  m(B)\, ,
$$ for every $n\geq n_0(B)$, where $\kappa >0$ is an absolute
constant.

We consider the lim-sup set
$$
\Lambda (\Re \, ,\, \psi) =\left\{ x\in \mm \; ;\; \dist (x,p_i) <
\psi (\varpi (i)) \mbox{ for infinitely many } i\in I \right\}\; .
$$

\begin{theorem}[\cite{BDV2}, Theorem 3]\label{minf}
Let $(\mm , \dist ,m)$ be a compact metric space equipped with a
measure satisfying the property $(M)$. Let
$\rho $ and $u$ be a function and respectively a sequence as
above. Let $(\Re , \varpi )$ be a local $m$-ubiquitous system
relative to $(\rho , u ) $, let $\psi$ be an approximating
function and let $\varphi$ be a dimension function dominating
$x^\delta$.

If $\sum_{n=1}^\infty \varphi (\psi (u_n))\, \rho (u_n)^{-\delta
}=\infty$ then $$\hh^{\varphi }
    \left( \Lambda (\Re \, ,\, \psi) \right)= \infty \, .$$
\end{theorem}

\begin{cor}[\cite{BDV2}, $\S 1.4.2$]\label{cbdv}
Let $(\mm , \dist ,m)$, $\rho $, $u$, $(\Re , \varpi )$ be as
above. Let $s\in [0, \delta )$.
\begin{itemize}
\item[(1)] If $\sum_{n=1}^\infty \psi (u_n)^s\, \rho (u_n)^{-\delta
}=\infty$ then $\hh^s
    \left( \Lambda (\Re \, ,\, \psi) \right)= \infty \, .$

\item[(2)] If $\lim_{n\to \infty } \frac{\psi (u_n)}{\rho
    (u_n)}=0$ then
    $$
\dim_H \Lambda (\Re \, ,\, \psi) \geq \sigma \delta \, ,
    $$ where $\sigma = \limsup_{n\to \infty } \frac{\ln \rho (u_n)}{\ln \psi
    (u_n)}$.

    Moreover, if $\liminf_{n\to \infty } \frac{\rho (u_n)}{\psi
    (u_n)^\sigma } < \infty $ then $\hh^{\sigma\delta } (\Lambda (\Re \, ,\, \psi) )=\infty
    $.
\end{itemize}
\end{cor}

\subsection{A geometric definition of $\mathcal{S}_\psi ({\mathfrak{Q}}_\q )$
}\label{trei}

We proceed with the proof of Theorem \ref{T2}. With a countable
covering argument, we reduce the problem to the study of
$\mathcal{S}_{\psi } (\Omega )$, where $\Omega$ is a relatively
compact open subset of ${\mathfrak{Q}}_{\q}$ such that its closure
does not intersect $T_{\bar{x}_0}{\mathfrak{Q}}_{\q}$ for some
$\bar{x}_0\in {\mathfrak{Q}}_{\q}$. Let $L_\q :\R^{n+1} \to \R$ be
defined as in the Introduction.

\medskip

\textit{Convention}: In what follows we shall drop the index of
the form $L_\q$ and we shall adopt for it all the notation
introduced in the Sections \ref{parab} to \ref{bparab}.

\medskip

We note that $\bar{y}\in {\mathfrak{Q}}_{\q}$ if and only if
$(\bar{y},1)\in Con_L$. Thus we may identify $\Omega$ to an open
subset of $\proj Con_L$, and consider $\R \Omega$. We denote by
$v_0$ the vector $(\bar{x}_0,1)\in Con_L $ and by
${\mathcal{H}}_0$ the hyperplane $\ker b_L(v_0, \cdot )$. The
condition that the closure of $\Omega $ does not intersect
$T_{\bar{x}_0}{\mathfrak{Q}}_{\q}$ is equivalent to the condition
that the closure of $\R \Omega $ does not intersect
${\mathcal{H}}_0$. Let $d_0=\R v_0$. Let $ P(d_0)$ be the
parabolic subgroup corresponding to $d_0$ in $SO_{I}(L)$ and let
$M_0A_0U_0$ be a Langlands decomposition of $P(d_0)$. In this case
the group $U_0$ is Abelian and $A_0$ is a maximal singular torus
$(\mathbf{a}_t)_{t\in \R }$. We take its parametrization such that
$U_0=U_+(A_0^+)$, where $A_0^+=(\mathbf{a}_t)_{t\geq 0 }$. Let
$\frak G$ be the maximal singular geodesic in $\calp_{n+1}(L)$
which is an orbit of $A_0$, such that ${\frak G}(0)\in H_{v_0}$
and ${\frak G} (-\infty )=d_0$. The geodesic ray $r_0={\frak
G}|_{[0,+\infty )}$ has $U_+(r_0)=U_0$.

 According to
Remark \ref{ud}, to every $d\in \proj \left( Con_L \setminus
{\mathcal{H}}_0 \right)$ corresponds a unique $\mathbf{u}_d \in
U_0$ such that $d=r_0(\infty ) \mathbf{u}_d$. By means of this, we
identify $\Omega$ to a relatively compact open subset of $U_0$.

\medskip

\noindent {\it Convention}: We also denote by $\Omega $ the open
subset of $\proj Con_L$ to which $\Omega \subset
{\mathfrak{Q}}_{\q} $ is identified via the map $\bar{y}\in
{\mathfrak{Q}}_{\q} \mapsto (\bar{y},1)\in Con_L$.
 We likewise denote by $\Omega $ the subset in $U_0$
  to which the previous subset in $\proj Con_L$ is identified via the map
   defined above.

\medskip

We want to study the set of vectors $\bar{y}\in \Omega$ such that
$\| \bar{y} - \frac{1}{q}\bar{p}\| < \frac{\psi (q)}{q}$ for
infinitely many $\frac{1}{q}\bar{p} \in {\mathfrak{Q}}_{\q}\cap
\Q^n$. Without loss of generality we may suppose moreover that
$\frac{1}{q}\bar{p} \in \Omega\cap \Q^n$.

Let $\bar{y}\, ,\, \bar{y}'$ be points in $\Omega $ and let $d\,
,\, d'$ be the lines in $Con_L$ containing $(\bar{y}, 1)$ and
$(\bar{y}', 1)$, respectively. Lemma \ref{uang} implies that the
distance in $U_0$ between $\mathbf{u}_d$ and $\mathbf{u}_{d'}$ is
bi-Lipschitz equivalent to the angle between $d$ and $ d'$, which
in its turn is bi-Lipschitz equivalent to $\|\bar{y}-\bar{y}' \|$
for $\bar{y}\, ,\, \bar{y}' \in \Omega$. Corollary \ref{norm}
implies that $\| (\bar{p},q)\|_e$  is bi-Lipschitz equivalent to
$e^{{\mathrm{odist}}\, (Hb_{v_0},Hb_{(\bar{p},q)})/2\sqrt{2}}$.
Also, since $\frac{1}{q}\bar{p}\in \Omega$, and $\Omega$ is
relatively compact, $\| (\bar{p},q)\|_e$  is bi-Lipschitz
equivalent to $q$.

Thus we have
\begin{displaymath}
\frac{1}{L} \, e^{{\mathrm{odist}}\,
(Hb_{v_0},Hb_{(\bar{p},q)})/2\sqrt{2}}\leq q\leq L\,
e^{{\mathrm{odist}}\, (Hb_{v_0},Hb_{(\bar{p},q)})/2\sqrt{2}}\, ,
\end{displaymath} for some $L>1$ depending on $\Omega$.

Let us denote by $\widetilde{\mathcal{S}}^0_\Psi (\Omega ) $ the
following subset of $U_0$.
$$
\widetilde{\mathcal{S}}^0_\Psi (\Omega ) = \{ \mathbf{u} \in
\Omega \; ;\;  \dist_{U_0}(\mathbf{u},\mathbf{u}_{(\bar{p},q)})
\leq \Psi \left( {\mathrm{odist}}\,
(Hb_{v_0},Hb_{(\bar{p},q)})\right)\mbox{ for infinitely many }
$$
\begin{equation}\label{defintr}
(\bar{p}, q)\in \pri_+^{n+1} \cap Con_L \}\, .
\end{equation}

The previous considerations imply that for some $L_1>1$ we have
the inclusions
\begin{equation}\label{dincls}
 \widetilde{\mathcal{S}}^0_{\Psi_1} (\Omega )\subset
\mathcal{S}_\psi(\Omega ) \subset
\widetilde{\mathcal{S}}^0_{\Psi_2} (\Omega )\, ,\mbox{ where
}\Psi_1(x)=\frac{\psi \left( L e^{x/2\sqrt{2}} \right)}{L_1
e^{x/2\sqrt{2}}}\mbox{ and }\Psi_2(x)=L_1\frac{\psi \left(
\frac{e^{x/2\sqrt{2}}}{L} \right)}{e^{x/2\sqrt{2}}}\, .
\end{equation}

Therefore we may replace in our study the set
$\mathcal{S}_\psi(\Omega )$ by the set
$\widetilde{\mathcal{S}}^0_\Psi (\Omega ) $ for some approximating
function $\Psi$.

\me

Let $\Gamma_L$, $\mathcal{V}_L$ and $\{ v_1\, ,\, v_2\, ,\, \dots
\, ,\, v_k \} \subset \pri_+^{n+1} \cap Con_L$ be defined as in
Section \ref{bparab}. Every $(\bar{p}, q)\in \pri_+^{n+1} \cap
Con_L$ is in the $\Gamma_L$-orbit
   of one of the vectors $\{ v_1\, ,\, v_2\, ,\, \dots \, ,\, v_k
  \}$. We fix a vector $v\in \{ v_1\, ,\, v_2\, ,\, \dots \,
,\, v_k \}$ and we consider the set of vectors $Csp=v \Gamma_L$.
%According to Lemma \ref{c}, we may replace in the
%definition of $\widetilde{\mathcal{S}}^0_\alpha (\Omega )$
  %the horoballs $Hb_{(\bar{p},q)}$ with
  %$Hb_{(\bar{p},q)}^{-\tau }$, with $\tau $ chosen as in Section
  %\ref{bparab}. We do this replacement, in order to have pairwise disjoint horoballs.

  %Let $w_i= e^{\tau /2\sqrt{2}} v_i$, for $i\in \{ 1,2, \dots k\}$. We introduce the new vectors $w_i$ for simplicity, as $Hb^\tau_{v_i}=Hb_{w_i}$.

 We restrict our attention to a subset of
  $\widetilde{\mathcal{S}}^0_\Psi (\Omega )$ defined as follows.
\begin{equation}\label{i}
\widetilde{\mathcal{S}}_\Psi (\Omega ) = \left\{ \mathbf{u} \in
\Omega \; ;\;  \dist_{U_0}(\mathbf{u},\mathbf{u}_{w}) \leq
\Psi\left( {\mathrm{odist}}\, (Hb_{v_0},Hb_{w})\right)\mbox{ for
infinitely many }w\in Csp \right\} \, .
\end{equation}

Since the set $\widetilde{\mathcal{S}}^0_\Psi (\Omega )$ is the
union of the sets $\widetilde{\mathcal{S}}_\Psi (\Omega )$ defined
as in (\ref{i}) for $i\in \{ 1,2, \dots , k\}$, without loss of
generality we may replace in our argument
$\widetilde{\mathcal{S}}^0_\Psi (\Omega )$ by
$\widetilde{\mathcal{S}}_\Psi (\Omega )$.

\me

\noindent \textit{Notation}: In the particular case when $\Psi
(x)=e^{-\frac{(1+\alpha )x}{2\sqrt{2}}}$, with $\alpha
>0$, in all the previous notation the index $\Psi $ is replaced
by the index $\alpha$.

\me

If the quadratic form $\q$ is positive definite or of signature
$(1,n-1)$, $\mathcal{P}_s(L)$ is a model of the $n$-dimensional
hyperbolic space. Via the map $\bar{x}\mapsto (\bar{x},1)$,
${\mathfrak{Q}}_{\q}$ can be identified with $\proj \left( Con_L
\setminus \ker e_{n+1}^* \right)$. If $\q$ is positive definite
then $\proj \left( Con_L \setminus \ker e_{n+1}^* \right) =\proj
\left( Con_L \right)$, which is the whole boundary at infinity of
the hyperbolic space $\mathcal{P}_s(L)$. If $\q$ has signature
$(1,n-1)$ then $\proj \left( Con_L \setminus \ker e_{n+1}^*
\right)$ is an open Zariski dense subset of the boundary at
infinity of the hyperbolic space. We change the model
$\mathcal{P}_s(L)$ with the half-space model $\hip^n$ of the
hyperbolic space, and we suppose that $\infty =v_0$.

The set $ \widetilde{\mathcal{S}}_\Psi (\Omega )$ can be
identified with the set of points $\vartheta \in \Omega \subset
\R^{n-1}\subset
\partial_\infty \hip^n $ such that there
are infinitely many $\xi \in Csp$ satisfying the inequality
\begin{displaymath}
\dist_{U_0}(\mathbf{u}_\vartheta\, ,\, \mathbf{u}_\xi ) \leq
\Psi\left( {\mathrm{odist}}\, (Hb_{\infty }\, ,\, Hb_\xi)\right)\;
.
\end{displaymath}

The term on the left is $\|\vartheta-\xi \|_e$ and the term on the
right is, up to some insignificant changes of the function $\Psi$,
$\Psi (-\ln h_\xi)$, where $h_\xi $ is the Euclidean height of the
horoball $Hb_\xi$. The set $\widetilde{\mathcal{S}}_\Psi (\Omega
)$ becomes the set of points $\vartheta \in \Omega$ such that the
inequality
\begin{displaymath}
\|\vartheta-\xi \|_e \leq \Psi (-\ln h_\xi)
\end{displaymath}
has an infinity of solutions $\xi \in Csp$. The equality $\dim_H
\widetilde{\mathcal{S}}_\alpha (\Omega )=\frac{n-1}{1+\alpha }$ is
in this case a consequence of a more general result of R. Hill and
S. L. Velani \cite{HV}. Moreover, all the results stated in this
paper are proved in this particular case in \cite{BDV2}.

\me

From now on we may therefore suppose that the form $L =L_\q$ is of
signature $(a,b)$, with $\min (a,b)\geq 2$.

\me

\noindent \textit{Notation}: We denote $\dim U_0$ by $\Delta $.

\me

Since $s=a+b\geq 4$, we have $\Delta = s-2\geq 2 $.

\subsection{Sets of points on a quadric very well approximable by cusp points}\label{nonhip}

Let $L$ be a non-degenerate quadratic form of signature $(a,b)$,
with $\min (a,b) \geq 2 ,\, a+b=s$, and let $\Gamma$ be an
irreducible lattice in $SO_{I}(L)$. Let $\mathcal{V}$ be the
quotient $\mathcal{P}_{s}(L)/\Gamma$ and let $\pr$ denote the
projection of $\mathcal{P}_{s}(L)$ onto $\calv$. As in Section
\ref{bparab}, we consider the set $\{ \bar{r}_1\, ,\, \bar{r}_2\,
,\, \dots \, ,\, \bar{r}_k \}$ of all the maximal singular cusp
rays in $\calv $ such that their lifts $\{ r_1\, ,\, r_2\, ,\,
\dots \, ,\, r_k \}$ in $\calp_s(L)$ are of type $\wp $. Let $v_i$
be a vector on the line $r_i(\infty )$, $i=1,2,\dots ,k$. The
choice of $v_i$ shall be made more precise later. We fix a vector
$v\in \{ v_1\, ,\, v_2\, ,\, \dots \, ,\, v_k \}$ and we consider
the set of vectors $Csp=v \Gamma$.

Consider the set $\Omega$, the vector $v_0\in Con_L \setminus \{
0\} $, the geodesic $\frak G$, the ray $r_0$ and all the notation
and properties as in Section \ref{trei}. We define the projection
$\pi_0 : SO_I(L) \to \calp_s(L),\, \pi_0 (g)={\frak G}(0) g$ and
the projection $\pi : SO_I(L)/\Gamma \to \calv$ induced by
$\pi_0$.

Identified with a subset of $U_0$, $\Omega$ can be equipped with
the restrictions of a left-invariant metric $\dist $ and of the
Haar measure $\nu$. The space $(\Omega ,\dist , \nu )$ satisfies
the condition $(M)$ with $\delta = \Delta $,
provided that $\Omega $ is a ball in $U_0$. Without loss of
generality we assume that $\Omega$ is indeed a ball.

We consider $Csp\cap \R \Omega $ as a countable family of indices
and we define $$\varpi : Csp\cap \R \Omega \to \R_+\, ,\; \;
\varpi (w)=\mathrm{odist}\, (Hb_{v_0}\, ,\, Hb_{w})\, .$$

\me

\noindent  \textit{Notation}:  We denote the oriented distance
$\mathrm{odist}\, (Hb_{v_0}\, ,\, Hb_{w})$ by $d_w$.

\me

For each $w\in Csp$ let $\mathbf{u}_w$ be the unipotent in $U_0$
such that $r_0[\mathbf{u}_w](\infty )=\R v$. We consider the
collection of resonant points $\Re = \{ \mathbf{u}_w \: ;\: w\in
Csp\cap \R \Omega \}$. Finally we consider the sequence
$u=(u_n)_{n\in \N }\, $, $u_n=nT$, where $T>0$ is large enough,
and the ubiquitous function $\rho : \R_+ \to \R_+$, $\rho (x)=
\varkappa e^{-\frac{x}{2\sqrt{2}}}$, where $\varkappa $ is a
constant to be chosen later. The function $\rho $ is $u$-regular.

\begin{proposition}\label{hn}
The pair $(\Re , \varpi )$ is a local $\nu $-ubiquitous system
relative to $(\rho , u)$.
\end{proposition}

\proof We need to prove that for any ball $B$ in $\Omega $ of
radius at most $R_1$, for some $R_1>0$, we have
\begin{displaymath}
\nu \left( B\cap \bigcup_{w\in Csp\cap \R \Omega ,\, d_w\leq nT}
B\left( \mathbf{u}_w\, ,\, \varkappa e^{-\frac{nT}{2\sqrt{2}}}
\right) \right) \geq \kappa \, \nu (B)\, ,
\end{displaymath} for $n\geq n_0(B)$ and $\kappa$ an absolute
constant. Without loss of generality we may suppose that $B$ is an
open ball in $U_0$ which is entirely contained in $\Omega$. We
denote by $2B$ the ball with same center as $B$ and with double
radius.

According to Lemma \ref{bup}, (iii), multiplying the vector $v$
with a large positive constant $\eta $ means adding to $f_v$ a
large positive constant $2\sqrt{2}\ln \eta $. Thus $Hb_{\eta v}
=Hb_v^{-2\sqrt{2}\ln \eta}$. Suppose that $v$ has been re-scaled
so that Remark \ref{horob}, (1), applies to $Hb_{v}$. Let $f:
\calv \to \R $ be a $C^\infty$-function taking values in $[0,1]$,
such that $f=1$ on $\pr (Hb_{ \lambda v})$ for some $\lambda >1$
and close to $1$, and $f=0$ outside $\pr (Hb_v)$. According to
Proposition \ref{equidunip}, for any $\bar{g}_0 \in SO_{I}(L) /
\Gamma $,
\begin{equation}\label{e1apl}
  \oint_{B} f\circ \pi (\mathbf{a}_t\mathbf{u}\bar{g}_0 )\, d\nu (\mathbf{u})
  \rightarrow \oint_{SO_{I}(L) / \Gamma  } f\circ \pi \, d \mu \mbox{ as } t \to +\infty \; ,
\end{equation} where $\mu $ is the measure on $SO_{I}(L) / \Gamma $ induced
by the Haar measure on $SO_{I}(L)$. We take $g_0=id$. The
convergence in (\ref{e1apl}) implies that for $t\geq t_0(v,\lambda
)$,
$$
\frac{1}{\nu (B)}\int_B \unu_{\pr (Hb_v )}\circ \pi
(\mathbf{a}_t\mathbf{u}\bar{g}_0 )\, d\nu (\mathbf{u}) \geq
(1-\epsilon )c(L)Vol_v\, ,
$$ where $c(L)$ is the total measure of every maximal compact subgroup of $SO_I(L)$, $Vol_v$ is
the volume of $\pr (Hb_v)$ and
 $\epsilon $ is a small positive constant such that $\lim_{\lambda \to 1}\epsilon
 =0$. We choose and fix $\lambda $ such that $\epsilon
 =\frac{1}{4}$. The inequality above is equivalent to
\begin{equation}\label{lowm}
\nu \left( \left\{ \mathbf{u}\in B \; ;\;
\pi_0(\mathbf{a}_t\mathbf{u}) \in \bigcup_{w\in Csp}Hb_w \right\}
\right)\geq \frac{3}{4} \nu (B) c(L)Vol_v \, .
\end{equation}

In (\ref{lowm}) it is enough to take the subset of vectors $w$
from $Csp\cap \R \Omega$ with $d_w\leq t$. We note that
$\pi_0(\mathbf{a}_t\mathbf{u})=r_0(t)[\mathbf{u}]$ for $t\geq 0$.
In accordance with the notation in (\ref{trt}), let us denote
$$
Tr_t(w)=\{ \mathbf{u}\in U_0 \; ;\; r_0(t)[\mathbf{u}]\in Hb_w\}\,
.
$$

We may then write
$$
\left\{ \mathbf{u}\in B \: ;\: r_0(t)[\mathbf{u}] \in
\bigcup_{w\in Csp\cap \R \Omega ,\, d_w\leq t} Hb_w \right\} =
B\cap \bigcup_{w\in Csp\cap \R \Omega ,\, d_w\leq t} Tr_t(w)\, .
$$

Consider now the subset corresponding to horoballs at distance
less than $t -\tau $ for some $\tau >0$ to be chosen later, that
is
\begin{equation}\label{less}
B\cap \bigcup_{w\in Csp\cap \R \Omega ,\, d_w< t-\tau } Tr_t(w)\,
.
\end{equation}

Suppose that $t-\tau  = k T$ for some $k \in \N$. Then we can
write
$$
\bigcup_{w\in Csp\cap \R \Omega ,\, d_w< t-\tau } Tr_t(w) =
\bigsqcup_{j=-k_0}^{k} \bigcup_{w\in Csp\cap \R \Omega ,\, d_w\in
[(j-1)T, jT) } Tr_t(w)\, .
$$

We have that
$$
\nu \left( B\cap \bigsqcup_{j=-k_0}^{k} \bigcup_{w\in Csp\cap \R
\Omega ,\, d_w\in [(j-1)T, jT) } Tr_t(w) \right) \leq
\sum_{j=-k_0}^{k}\; \sum_{w\in Csp\cap \R \Omega ,\, d_w\in
[(j-1)T, jT)} \nu \left( B\cap Tr_t(w) \right)\, .
$$

 According to Lemma \ref{bilevoltr}, (1), $Tr_t(w)\subset B \left( \mathbf{u}_w, \kappa_0 e^{-\frac{d_w}{2\sqrt{2}}}
 \right)$, where $\kappa_0$ depends on $v_0$ and on the metric
 $\dist$. Hence there exists $J_0=J_0(B, \kappa_0)$ such that for $j\geq J_0$
  the intersection $B\cap Tr_t(w)$ is
non-empty only for $w\in Csp$ with $\mathbf{u}_w\in 2B$. To
designate the property that $w\in Csp$ with $\mathbf{u}_w\in 2B$
we use the notation $w\in Csp \cap \R (2B)$. Thus for $j\geq J_0$
the sum is taken over the $w\in Csp \cap \R (2B),\, d_w\in
[(j-1)T, jT)$. Lemma \ref{bilevoltr}, (2), implies that
\begin{equation}\label{bound}
\sum_{w\in Csp \cap \R (2B) ,\, d_w\in [(j-1)T, jT)} \nu \left(
B\cap Tr_t(w) \right)\leq C \sum_{w\in Csp \cap \R (2B),\, d_w\in
[(j-1)T, jT)} e^{-\frac{(j-1)T \Delta + t - jT}{2\sqrt{2}}}\, ,
\end{equation}
 where $C$ is a constant depending on $v_0$. Corollary
 \ref{lSu} now gives that for any $T\geq T_0(L,v_0, \Omega )$,
  the term in (\ref{bound}) is smaller
 than
 $$
\mathbf{K} e^{-\frac{(j-1)T \Delta + t - jT}{2\sqrt{2}}}
e^{\frac{j T\Delta}{2\sqrt{2}}}\nu (2 B) = \mathbf{K}'
e^{\frac{T\Delta +jT -t}{2\sqrt{2}}} \nu (B)\, ,
 $$ for every $j\geq J_1$, where $J_1=J_1(B,v,v_0)$ and $\mathbf{K}'=\mathbf{K}'(L,\Gamma , v_0,\Omega)$.

Let $J_2 = \max (J_0 , J_1 )$. The considerations above and Lemma
\ref{bilevoltr}, (2), imply that the measure of the set in
(\ref{less}) is at most
\begin{equation}\label{bound2}
    \mathbf{K}' e^{\frac{T\Delta -t}{2\sqrt{2}}} \nu (B)
\sum_{j=J_2}^{k}e^{\frac{jT}{2\sqrt{2}}} +C \sum_{j=-k_0}^{J_2}
\sum_{w\in Csp\cap \R \Omega , d_w\in
[(j-1)T,jT)}e^{-\frac{t+(\Delta -1)(j-1)T}{2\sqrt{2}}}\, .
\end{equation}

The set of $w\in Csp\cap \R \Omega$ with $d_w<J_2 T$ is finite, of
cardinal $N$. Hence the second term in the sum above is less than
$$
C N \sum_{j=-k_0}^{J_2} e^{-\frac{t+(\Delta
-1)(j-1)T}{2\sqrt{2}}}\leq C'e^{-\frac{t}{2\sqrt{2}}}\, ,
$$ where $C'=C'(v_0,v,\dist , T, B)$.

On the whole, for $T \geq 2 \sqrt{2}\ln 2$, the sum in
(\ref{bound2}) is at most
$$
2\mathbf{K}' e^{\frac{T\Delta -t}{2\sqrt{2}}} \nu (B)
e^{\frac{kT}{2\sqrt{2}}} + C'e^{-\frac{t}{2\sqrt{2}}}\; \leq \;
2\mathbf{K}' e^{\frac{T\Delta -\tau }{2\sqrt{2}}} \nu (B)
 + C'e^{-\frac{t}{2\sqrt{2}}}\, ,
$$

We choose $\tau$ such that $2\mathbf{K}' e^{\frac{T\Delta -\tau
}{2\sqrt{2}}}=\frac{1}{4} Vol_v c(L)$. Note that it depends on
$L,\Gamma, \Omega ,v,v_0$ and $T$. Then for $t\geq t_1$ the
measure of the set in (\ref{less}) is smaller than $\frac{1}{2}
Vol_v c(L)\nu (B)$, where $t_1$ depends on $L,v_0,v,\dist , T,B$.

Let $t_2=\max (t_0,t_1)$. It follows that for $t\geq t_2$ the set
\begin{displaymath}
B\cap \bigcup_{w\in Csp\cap \R \Omega ,\, t-\tau \leq d_w\leq t}
Tr_t(w)
\end{displaymath} has measure at least $\frac{1}{4} Vol_v c(L)\nu
  (B)$. By Lemma \ref{bilevoltr}, (1), this set is included in
\begin{displaymath}
B\cap \bigcup_{w\in Csp\cap \R \Omega ,\, t-\tau \leq d_w\leq t}
B\left( \mathbf{u}_w \, ,\, \kappa_0 e^{-\frac{d_w}{2\sqrt{2}}}
\right) \; \subset \; B\cap \bigcup_{w\in Csp\cap \R \Omega ,\,
t-\tau \leq d_w\leq t} B\left( \mathbf{u}_w \, ,\, \kappa_0
e^{\frac{\tau }{2\sqrt{2}}} e^{-\frac{t}{2\sqrt{2}}} \right)\, .
\end{displaymath}

We choose the constant $\varkappa = \kappa_0 e^{\frac{\tau
}{2\sqrt{2}}}$ and $t = nT$. We have obtained that
$$
\nu \left( B\cap \bigcup_{w\in Csp\cap \R \Omega ,\, d_w\leq nT}
B\left( \mathbf{u}_w , \varkappa e^{-\frac{nT}{2\sqrt{2}}} \right)
\right)\geq \frac{1}{4} Vol_v c(L)\nu
  (B)\, ,
$$ for $n\geq n_0$, with $n_0$ depending on the data of the ubiquitous system and on $B$.
 This finishes the proof.\endproof

\me

For an approximating function $\psi$ we can define the set
\begin{equation}\label{igen}
\widetilde{\mathcal{S}}_\psi (\Omega ) = \left\{ \mathbf{u} \in
\Omega \; ;\;  \dist_{U_0}(\mathbf{u},\mathbf{u}_{w}) \leq \psi
(d_w)\mbox{ for infinitely many }w\in Csp \right\} \, .
\end{equation}

When $\psi (x)=e^{-\frac{(1+\alpha )x }{2\sqrt{2}}}$, with $\alpha
>0$, we replace the index $\psi$ by the index $\alpha $.

Theorem \ref{minf} implies the following.

\me

\begin{theorem}\label{minfapl}
Let $\varphi$ be a dimension function dominating $x^\Delta$. The
measure $\hh^\varphi\left(\widetilde{\mathcal{S}}_\psi
\left(\Omega
    \right)\right)$ is $+\infty$ if and only if for some/for every $T>0$ large enough
    $\sum_{n=1}^\infty \varphi (\psi (nT))e^{\frac{nT\Delta}{2\sqrt{2}}} =\infty
    $.
\end{theorem}

\begin{remark}
The significance of the alternative use of the conditions ``for
some/for every'' is the following: the ``if'' part holds under the
weaker condition that the sum is $+\infty$ for some $T>0$ large
enough, while in the ``only if'' part we obtain that the sum is
$+\infty$ for every $T>0$ large enough.
\end{remark}

\proof  The ``if'' part follows from Theorem \ref{minf}. We show
that if for some $T>0$ large enough
    $\sum_{n=1}^\infty \varphi (\psi (nT))e^{\frac{nT\Delta}{2\sqrt{2}}} <\infty
    $ then $\hh^\varphi\left(\widetilde{\mathcal{S}}_\psi
\left(\Omega
    \right)\right)<\infty$.

The set $\widetilde{\mathcal{S}}_\psi (\Omega )$ is covered by
$$
\bigcup_{w\in Csp \cap \R \nn_\varepsilon (\Omega )}
B\left(\mathbf{u}_w\, ,\, \psi (d_w) \right)\, .
$$

We have $\sum_{w\in Csp \cap \R \nn_\varepsilon (\Omega )} \varphi
\left( \psi (d_w) \right) =\sum_{n=1}^{\infty }\sum_{w\in Csp \cap
\R \nn_\varepsilon (\Omega )\, ,\, d_w\in [nT,(n+1)T)} \varphi
\left( \psi (d_w) \right) \ll \\
\sum_{n=1}^{\infty }\varphi \left( \psi (nT ) \right)
e^{\frac{\Delta nT}{2\sqrt{2}}}$. The last inequality follows from
Corollary \ref{lSu}.\endproof

\me

\begin{cor}\label{cbdvapl}
\begin{itemize}
    \item[(1)] Let $s\in [0,\Delta)$. We have $\hh^s
    \left( \widetilde{\mathcal{S}}_\psi (\Omega
    ) \right)= \infty $ if and only if
    for some/for every $T>0$ large enough
    $$\sum_{n=1}^\infty \psi (nT)^s\, e^{\frac{nT\Delta}{2\sqrt{2}}}=\infty\, .
    $$
    \item[(2)] If for some $T>0$ large enough $\lim_{n\to \infty } \psi (nT)e^{\frac{nT}{2\sqrt{2}}}=0$ then
    $$
\dim_H \widetilde{\mathcal{S}}_\psi (\Omega
    ) = \sigma \Delta \, ,
    $$ where $\sigma = \limsup_{n\to \infty } \frac{-nT}{2\sqrt{2}\ln \psi
    (nT)}$.

    Moreover, if $\limsup_{n\to \infty } e^{\frac{nT}{2\sqrt{2}}}\psi
    (nT)^\sigma  >0 $ then $\hh^{\sigma\Delta } \left(\widetilde{\mathcal{S}}_\psi (\Omega
    )\right)=\infty
    $.
\end{itemize}
\end{cor}

\proof Statement (1) follows from Theorem \ref{minfapl} applied to
$\varphi (x)=x^s$.

\me

\noindent (2) The definition of $\sigma$ implies that for any
$\epsilon
>0$ the following holds.
\begin{itemize}
    \item[(a)] For $n$ large enough,
    $\psi (nT) \leq e^{-\frac{nT}{2\sqrt{2}(\sigma + \epsilon
    )}}$;
    \item[(b)] For infinitely many $n$,
    $\psi (nT) \geq e^{-\frac{nT}{2\sqrt{2}(\sigma - \epsilon
    )}}$.
\end{itemize}

According to (a), for every $s > \Delta (\sigma +\epsilon )$, $
\sum_{n=1}^\infty \psi (nT)^s e^{\frac{nT\Delta}{2\sqrt{2}}}\ll
\sum_{n=1}^\infty e^{\frac{nT}{2\sqrt{2}}\left(\Delta -
\frac{s}{\sigma + \epsilon }\right)} < + \infty \, .$ Statement
(1) implies that $\hh^s \left(\widetilde{\mathcal{S}}_\psi (\Omega
    )\right)< \infty$.

Property (b) implies that $\psi (nT) e^{\frac{nT}{2\sqrt{2}(\sigma
- \epsilon
    )}}\geq 1$ for infinitely many $n$.
Statement (1) implies that $\hh^{\Delta (\sigma -\epsilon )}
\left(\widetilde{\mathcal{S}}_\psi (\Omega
    )\right)= \infty$.

We thus obtain that $\Delta (\sigma -\epsilon ) \leq \dim_H
\widetilde{\mathcal{S}}_\psi (\Omega )\leq \Delta (\sigma
+\epsilon )$, for any $\epsilon >0$, which implies that $\dim_H
\widetilde{\mathcal{S}}_\psi (\Omega )= \Delta \sigma $.

The last statement in (2) follows from (1) applied to $s= \sigma
\Delta$.\endproof

\me

\begin{cor}\label{stildea}
\begin{itemize}
    \item[(i)] The set $\widetilde{\mathcal{S}}_\alpha (\Omega )$ has Hausdorff
dimension $d=\frac{\dim U_0}{1+\alpha}$ for every $\alpha >0$, and
$\hh^d (\widetilde{\mathcal{S}}_\alpha (\Omega ))=\infty$.
    \item[(ii)] Both statements hold also for
$\mathcal{E}\widetilde{\mathcal{S}}_\alpha (\Omega ) =
\widetilde{\mathcal{S}}_\alpha (\Omega )\setminus \bigcup_{\beta
> \alpha }\widetilde{\mathcal{S}}_\beta (\Omega )$.
\end{itemize}
\end{cor}

\proof (i) is a consequence of Corollary \ref{cbdvapl}, (2).

\me

(ii) We can write $\mathcal{E}\widetilde{\mathcal{S}}_\alpha
(\Omega ) = \widetilde{\mathcal{S}}_\alpha (\Omega )\setminus
\bigcup_{n\in \N }\widetilde{\mathcal{S}}_{\alpha +\frac{1}{n}}
(\Omega )$. According to (i) we have, for $d=\frac{\dim
U_0}{1+\alpha}$, that $\hh^d (\widetilde{\mathcal{S}}_\alpha
(\Omega ))=\infty$ and $\hh^d \left(
\widetilde{\mathcal{S}}_{\alpha+\frac{1}{n}} (\Omega )\right)=0$.
Hence $\hh^d (\mathcal{E}\widetilde{\mathcal{S}}_\alpha (\Omega
))=\infty$.\endproof

\me

\begin{remark}[possible generalizations]\label{gen}
One might work in the general setting, that is when instead of
$\mathcal{P}_s(L)$ and $SO_{I}(L)$ there is a general symmetric
space of non-compact type $X$ and the connected component of the
identity of its group of isometries $G$, and when $\Gamma $ is an
irreducible lattice in $G$. One can consider a maximal singular
cusp ray $\bar{r}$ in $X/\Gamma$, a lift $r$ in $X$, $Csp =
r(\infty )\Gamma $ and for every $\xi \in Csp $ the horoball
$Hb_\xi$ of basepoint $\xi$ projecting onto $Hb(\bar{r})$.

For an arbitrary ray $\rho $ in the orbit $rG$ and $U_+=U_+(\rho
)$, the set $\rho (\infty )U_+$ is open and Zariski dense in
$r(\infty )G$ and contains infinitely many $\xi \in Csp$ to which
therefore one can associate unipotents $\mathbf{u}_\xi \in U_+$. A
set $\widetilde{\mathcal{S}}_\psi $ can be defined as before, that
is as the set
$$
\widetilde{\mathcal{S}}_\psi = \left\{ \mathbf{u} \in U_+ \; ;\;
\dist_{U_+}(\mathbf{u},\mathbf{u}_{\xi}) \leq \psi \left(
{\mathrm{odist}}\, (Hb_0,Hb_{\xi}) \right)\mbox{ for infinitely
many }\xi \in Csp \right\} \, ,
$$ where $Hb_0$ is the horoball determined by the ray opposite to $\rho$.

Let $(\mathbf{a}_t)_{t\in \R}$ be the maximal singular torus such
that $\rho (t)=\rho(0)\mathbf{a}_t$ and suppose that
$\mathbf{u}\mapsto
    \mathbf{a}_{t}\mathbf{u}\mathbf{a}_{-t}$ is a homothety of $U_+$ of factor $e^{\lambda
    t}$, $\lambda >0$. In order to prove that
    $$
    d=\dim_H \widetilde{\mathcal{S}}_\psi=\mathrm{dim}\, U_+\, \sigma
    \; \mbox{ and }\; \hh^d(\widetilde{\mathcal{S}}_\psi )=\infty\, ,\; \mbox{ for
    }\; \sigma=\limsup_{n\to \infty }\frac{-\lambda nT}{\ln \psi
    (nT)}\, ,
    $$ the following conditions are sufficient:
\begin{itemize}
    \item the equidistribution results given in Propositions
    \ref{equidunip} and \ref{equiduc} and the counting result Corollary \ref{lSu}; these hold in general;
    \item for $Tr_{D+\tau }$ as defined in (\ref{trt}), a formula
    of the measure of the form $$\nu (Tr_{D+\tau })=e^{-D\lambda \dim U_+ - f(\tau
    )}\mbox{ with }\lim_{\tau \to \infty } f(\tau )=\infty \, ;
    $$
    \item the inclusion $Tr_{D+\tau } \subset B(\mathbf{u}_\xi , C e^{-\lambda
    D})$, with $C$ an absolute constant.
\end{itemize}
\end{remark}

A consequence of Theorem \ref{minfapl} is the following.

\begin{theorem}\label{resda}
Let $\q :\R^m \to \R$ be a non-degenerate quadratic form with
rational coefficients, let ${\mathfrak{Q}}_\q$ be the quadric
defined by $\q=1$ and let $\psi $ be an approximating function
such that $\lim_{x\to \infty } x\psi (x)=0$. Let $\varphi$ be a
dimension function
    dominating $x^{m-1}$. Then $\hh^\varphi(\mathcal{S}_\psi (\qq_\q
    ))=\infty$ if and only if for some/for every $T>0$ large
    enough
$\sum_{n=1}^\infty \varphi \left(\frac{\psi
(e^{nT})}{e^{nT}}\right)e^{nT(m-1)} =\infty $.
\end{theorem}

Corollaries \ref{cbdvapl} and \ref{stildea} applied in this
setting yield Theorem \ref{T2}.

\begin{remark}
The results in Corollary \ref{stildea} and in Theorem \ref{T2}
concerning the sets $\mathcal{E}\mathcal{S}_\alpha $ of vectors
for which the exact order of approximation is known, can be
formulated in the more general context of approximating functions.
See for example \cite[$\S 8.8$]{BDV2} and \cite{BDV3} for such
results. For the sake of simplicity we have not done it here.
\end{remark}

\section{Sets of geodesic rays moving away into the
cusp}\label{sectrays}

\subsection{The case of $\mathcal{P}_s(L)$ and of the geodesic rays of type $\wp$}\label{sectrays1}

We consider $L$, $\Gamma$, $\calv$, $\pr$, $\bar{r}_i$, $r_i$ and
$v_i$, $i\in \{1,2,\dots ,k \}$, with the same significance as in
the beginning of Section \ref{nonhip}, with the only difference
that the condition $\min (a,b)\geq 2$ is replaced by $\min
(a,b)\geq 1$. Without loss of generality we may suppose that the
vector $v_i$ is such that $f_{r_i}=f_{v_i}$. Let $\varrho $ be an
arbitrary geodesic ray of type $\wp$ in $\mathcal{P}_s(L)$ and let
$U_0=U_+(\varrho )$. We denote $\dim U_0$ by $\Delta$.

Consider a function $\phi : [a,+\infty ) \to [b,+\infty )\, ,\,
a,b\in \R $, and define the set of unipotents
\begin{equation}\label{rfi}
    \mathcal{R}_{\phi }^i = \left\{ \mathbf{u}\in U_0 \; ;\;
-f_{\bar{r}_i}\left(\pr\left(\varrho (t)\mathbf{u}\right)\right)
\geq t-\phi (t)\mbox{ infinitely many times as }t\to \infty
\right\}\, .
\end{equation}

\begin{remark}\label{mdepth}
The maximal possible depth for
$\pr\left(\varrho(t)\mathbf{u}\right)$, measured with respect to
the ray $\bar{r}_i$, is $t+c$, where $c$ is a constant. Such a
depth can occur infinitely many times if and only if the ray
$\pr\left(\varrho (t)\mathbf{u}\right)$ is asymptotic to
$\bar{r}_i$ (see Corollary \ref{corrbeta}, (c)). Therefore, it
makes sense to put $t-\phi(t)$ with $\phi$ a function bounded
below near $+\infty $, in the definition of $\mathcal{R}_{\phi
}^i$.
\end{remark}

In the particular case when $t-\phi (t)=\beta t$, with $\beta \in
[0,1],$ we  replace in our notation in (\ref{rfi}) the index
$\phi$ by the index $\beta$.

\begin{theorem}\label{diva}
Suppose that $\phi$ and $id-\phi $ are increasing functions, and
that $\phi$ is a bijection.
\begin{itemize}
    \item[(1)] Let $s\in [0,\Delta )$. $\hh^s (\mathcal{R}^i_{\phi } )=\infty$ if and only if for some/for every
    $T>0$ large enough,
    $$\sum_{n=1}^\infty e^{\frac{\Delta nT -s\phi\iv (nT)}{2\sqrt{2}}}=\infty \,
    .$$
    \item[(2)] If for some $T>0$ large enough, $\lim_{n\to \infty } [nT -\phi\iv (nT)]=-\infty$ then
    $$d=\dim_H \mathcal{R}^i_{\phi } = \sigma \cdot \Delta \mbox{, where }
\sigma = \limsup_{n\to \infty }\frac{nT}{\phi\iv (nT)}\, .
    $$

    If moreover $\limsup_{n\to \infty} [nT - \sigma \phi\iv (nT)] >
    -\infty$ then $\hh^d (\mathcal{R}^i_{\phi } )=\infty$.
\end{itemize}
\end{theorem}

\proof  We denote by $Csp$ the set $ v_i\Gamma$.
%Let ${A}=(\mathbf{a}_t)_{t\in \R }$ be a maximal singular
%torus in $SO_{I}(L)$ and $\mathcal{A}_+=(\mathbf{a}_t)_{t\geq 0 }$ a
%sub-semigroup such that $r(t)= r(0)\mathbf{a}_t\, ,\, t\geq 0$.
We may restrict our study to a set $\mathcal{R}^i_{\phi } (\Omega
)=\mathcal{R}^i_{\phi }\cap \Omega $, where $\Omega $ is an open
relatively compact subset of $U_0$. Let $d_0$ be the line in
$Con_L$ which appears as point at infinity of the geodesic ray
$\varrho^{op}$ opposite to $\varrho$. We choose the vector $v_0$
on $d_0$ so that $f_{\varrho^{op}}=f_{v_0}$. For every vector
$v\in Con_L$ such that $b_L(v,v_0)\neq 0$, we denote by
$\mathbf{u}_v$ the element in $U_0$ such that the geodesic ray
$\varrho \mathbf{u}_v$ has as point at infinity the line $\R v$.
%According to Section \ref{lssp}, we denote by $\pr_{r(0)}$ the
%projection $\pr_{r(0)}: SO_{I}(L) \to \mathcal{P}_s(L)\,
%,\,\pr_{r(0)}(g)=r(0)g$.

\me

\noindent \textbf{I.} We show that $\mathcal{R}^i_{\phi}(\Omega)
\subset \widetilde{\mathcal{S}}_{ \Phi_1 }(\Omega )$, where
$\widetilde{\mathcal{S}}_{ \Phi_1 }(\Omega )$ is defined as in
(\ref{i}) for the approximating function $\Phi_1 (x) = \kappa_0
e^{-\frac{\phi \iv (x)}{2\sqrt{2}}}$. Here $\kappa_0$ is the
constant appearing in the inclusion (\ref{trbile}).

Let $\mathbf{u}\in \mathcal{R}^i_{\phi}(\Omega)$. It follows that
the inequality
\begin{equation}\label{first}
f_w\left( \varrho(t)\mathbf{u}\right) \leq -(t-\phi(t))\; \;
\end{equation} has infinitely many solutions $w\in Csp$ and $t\in (0,+\infty )$.
Let $w$ and $t$ be two such solutions. Then, with the notation
$d_w= {\mathrm{odist}}\, (Hb_{v_0}\, ,\, Hb_w)$, we have
\begin{equation}\label{dv1}
d_w \leq t - |f_w(\varrho(t)\mathbf{u})| \leq \phi(t)\;
\Rightarrow \; \phi\iv (d_w)\leq t \; .
\end{equation}

The fact that $\varrho(t)\mathbf{u}\in Hb_w^{-(t -\phi (t))}$
implies by Lemma \ref{bilevoltr} that $$\dist_{U_0}(\mathbf{u}\,
,\, \mathbf{u}_w)\leq \kappa_0
e^{-\frac{1}{2\sqrt{2}}\left(d_w+t-\phi(t)\right)}\, .$$ The last term
of the inequality is by (\ref{dv1}) smaller than $\kappa_0
e^{-\frac{\phi\iv (d_w)}{2\sqrt{2}}}$. We conclude that
$\mathcal{R}^i_{\phi}(\Omega) \subset \widetilde{\mathcal{S}}_{
\Phi_1 }(\Omega )$.

\me

\noindent \textbf{II.} We show that $\mathcal{R}^i_{\phi}(\Omega)$
contains $\widetilde{\mathcal{S}}_{\Phi_2}(\Omega) $, where
$\Phi_2(x)=ce^{-\frac{\phi\iv (x+1)}{2\sqrt{2}}}$ for an
appropriate constant $c$. Let $\mathbf{u} \in
\widetilde{\mathcal{S}}_{\Phi_2 }(\Omega )$. Let $w\in Csp$ be
such that $\dist_{U_0}(\mathbf{u}\, ,\, \mathbf{u}_w)\leq c
e^{-\frac{\phi\iv (d_w+1)}{2\sqrt{2}}}$. We consider $t= \phi\iv
(d_w+1)$. We have that
\begin{displaymath}
\left| f_w\left( \varrho(t)\mathbf{u}\right) -f_w\left( \varrho(t)
\mathbf{u}_w \right)\right| \leq \dist \left(
\varrho(t)\mathbf{u}\, ,\, \varrho(t) \mathbf{u}_w \right)
$$
$$\leq e^{\frac{t}{2\sqrt{2}}} \dist \left( \varrho(0)\mathbf{u}\, ,\,
\varrho(0)\mathbf{u}_w \right)\leq 1 \, ,
\end{displaymath} if we choose the constant $c$ properly, depending on
$\varrho(0)$ and on the metric chosen on $U_0$. Since $f_w\left(
\varrho(t) \mathbf{u}_w\right) = d_w-t =\phi(t)-1-t$, this implies
that $f_w\left( \varrho(t) \mathbf{u}\right)\leq \phi(t)-t$. We
conclude that
$$
-f_{\bar{r}_i}(\pr (\varrho(t) \mathbf{u} )) \geq t-\phi (t)\mbox{
infinitely many times as }t\to \infty \, .
$$

We have obtained that
\begin{displaymath}
\widetilde{\mathcal{S}}_{\Phi_2}(\Omega) \subset
\mathcal{R}^i_{\phi}(\Omega) \subset \widetilde{\mathcal{S}}_{
\Phi_1 }(\Omega )\, ,
\end{displaymath} where $\Phi_1 (x) = \kappa_0 e^{-\frac{\phi \iv
(x)}{2\sqrt{2}}}$ and $\Phi_2(x)=ce^{-\frac{\phi\iv
(x+1)}{2\sqrt{2}}}$. We apply Corollary \ref{cbdvapl} and we
obtain the desired conclusion.\endproof

\me

\begin{remark}\label{ficonstq}
When defining the set $\calr^i_\phi$, one can replace $\phi$ by
$\phi_c=\phi -c$, where $c $ is a constant, and Theorem \ref{diva}
still holds. In order to see this it suffices to show, using the
monotonicity of $\phi$, that all the conditions in Theorem
\ref{diva} are satisfied by $\phi$ if and only if they are
satisfied by $\phi_c$. We leave this as an exercise to the reader.
See also Remark \ref{ficonst} where a similar statement is proved
in full detail.
\end{remark}

\me

\begin{remark}\label{last}
The set $ P(\varrho) U_0$ is open Zariski dense in $G$, hence the
projection of $U_0$ is open Zariski dense in $P(\varrho)\backslash
G$. We note that $P(\varrho)\backslash G $ is the stratum $\wp$,
in the terminology of the Introduction. We also note that if a
geodesic ray has a projection on $\mathcal{V}$ moving away into
the cusp infinitely many times with depth measured by the function
$id -\phi $ with respect to the ray $\bar{r}_i$, then any geodesic
ray asymptotic to it has the same property, up to a bounded
perturbation of the depth. The previous Theorem therefore gives
the formula of the Hausdorff dimension of the set of points of
type $\wp$ in the boundary at infinity corresponding to rays
moving away into the cusp at depth at least $id -\phi $ with
respect to $\bar{r}_i$, infinitely many times.
\end{remark}

\me

\begin{cor}\label{corrbeta}
Consider the set
\begin{equation}\label{rbeta}
    \mathcal{R}^i_{\beta } = \left\{ \mathbf{u}\in
U_0 \; ;\; -f_{\bar{r}_i}\left(\pr
\left(\varrho(t)\mathbf{u}\right)\right)\geq \beta t \mbox{
infinitely many times as }t\to \infty \right\}\, .
\end{equation}
\begin{itemize}
        \item[(a)] For any $\beta \in
(0,1)$ the set $ \mathcal{R}^i_{\beta }$ has Hausdorff dimension
$d=\Delta (1-\beta )$ and $\hh^d \left( \mathcal{R}^i_{\beta }
\right)=\infty$.
    \item[(b)] Both statements in (a) also hold for the set
    $$
\mathcal{E}\mathcal{R}^i_{\beta } = \mathcal{R}^i_{\beta
}\setminus \bigcup_{\beta' > \beta } \mathcal{R}^i_{\beta' } =
\left\{ \mathbf{u}\in \mathcal{R}^i_\beta \; ;\; \limsup_{t\to
+\infty } \frac{-f_{\bar{r}_i}\left(\pr
\left(\varrho(t)\mathbf{u}\right)\right)}{t} =\beta \right\}\, .
$$
\item[(c)] The set $\mathcal{R}^i_{0 }$ coincides with $U_0$ and the set
$\mathcal{R}^i_{1}$ is contained in $\left\{ \mathbf{u}_w\; ;\;
w\in Csp \right\}$.
\end{itemize}
\end{cor}

\proof (a) follows from Theorem \ref{diva}, (2), and (b)
immediately follows from (a).

\me

(c) For $\beta =0$ it is a consequence of the logarithm law
\cite{KM}. Suppose that $\beta =1$. Let $\mathbf{u}\in
\mathcal{R}^i_{1}$. Then for infinitely many $w\in Csp$ and $t\in
(0,+\infty )$ the following inequality holds:
\begin{displaymath}
f_w\left( \varrho(t)\mathbf{u}\right) \leq -t\, .
\end{displaymath}

As in (\ref{dv1}) we obtain that for every such $w$ we have
$d_w\leq 0$. The inclusion $\varrho(t)\mathbf{u}\in Hb_w^{-t}$
implies by Lemma \ref{bilevoltr} that
\begin{equation}\label{dist1}
    \dist_{U_0}(\mathbf{u}\, ,\, \mathbf{u}_w)\leq \kappa_0
e^{-\frac{1}{2\sqrt{2}}\left(d_w+t\right)}\, .
\end{equation}
Thus for $t$ large enough we may suppose that the corresponding
$w\in Csp$ satisfies $\mathbf{u}_w \in B(\mathbf{u},1)$. On the
other hand, the number of $w\in Csp$ with $\mathbf{u}_w \in
B(\mathbf{u},1)$ and $d_w\leq 0$ is finite. Hence, by eventually
taking a subsequence we may suppose that $w$ is fixed. By letting
$t\to \infty $ in (\ref{dist1}) we obtain that
$\mathbf{u}=\mathbf{u}_w$. Thus $\mathcal{R}^i_{1}\subset \left\{
\mathbf{u}_w\; ;\; w\in Csp \right\}$.
\endproof

\begin{remark}
Both Theorem \ref{diva} and Corollary \ref{corrbeta} follow from
Corollary \ref{cbdvapl} and inclusion (\ref{trbile}). Consequently
the conditions in Remark \ref{gen} are sufficient also for the
generalization of these two results.
\end{remark}

\subsection{The symmetric space $\mathcal{P}_{n+1}$ and the rays of slope $r_i(\infty )$, $i=1,n$}

Consider the symmetric space $\mathcal{P}_{n+1}\, $, with group of
isometries $\slnpr$, and the lattice $\slnpz$. Let
$\mathcal{T}_{n+1} = \calp_{n+1}/\slnpz$ and let $\pr$ be the
projection of $\calp_{n+1}$ onto $\mathcal{T}_{n+1}\, $.

Let $r_1$ and $r_n$ be the geodesic rays in $\calp_{n+1}$ defined
as in (\ref{rays}). The ray $r_i,\, i=1,n,$ projects onto a
maximal singular cusp ray in $\mathcal{T}_{n+1}$, which we denote
by $\bar{r}_i$. The point at infinity $r_1(\infty )$ is the
projective point $\langle e_{n+1}\rangle $. The point at infinity
$r_n (\infty )$ is the hyperplane in $\mathbb{P}^n\R$ defined by
$x_1=0$. We denote it by $e^*_1$.

The set $\mathcal{S}_\psi (\R^n )$ can be related to sets of
geodesic rays similar to $\mathcal{R}^i_\phi $ from (\ref{rfi}).
The formula (\ref{jar}) will then imply a result similar to
Theorem \ref{diva}. More precisely, define
\begin{equation}\label{rr1}
\mathcal{R}_{\phi }^1 = \left\{ \mathbf{u}\in U_+(r_1) \; ;\;
-f_{\bar{r}_1}\left(\pr \left(r_1(t)\mathbf{u}\right)\right)\geq t
- \phi (t) \mbox{ infinitely many times as }t\to \infty \right\}\,
.
\end{equation}

One can define a similar set for the ray $r_n$, that is
\begin{equation}\label{rrn}
\mathcal{R}_{\phi }^n = \left\{ \mathbf{u}\in U_+(r_n) \; ;\;
-f_{\bar{r}_n}\left(\pr \left(r_n(t)\mathbf{u}\right)\right)\geq t
- \phi (t) \mbox{ infinitely many times as }t\to \infty \right\}\,
.
\end{equation}

A remark similar to Remark \ref{mdepth} justifies the way
$\calr_\phi^i ,\, i=1,n,$ are defined. In the case when $t-\phi
(t)=\beta t$ we replace in (\ref{rr1}) and in (\ref{rrn}) the
index $\phi$ by $\beta$.

To simplify the formulas, we use the notation $\eta_n$ for the
constant $\sqrt{\frac{n+1}{n}}$.

\begin{theorem}\label{thri}
Let $\phi :[a,+\infty ) \to [b,+\infty )$ be a bijective function
such that $\phi$ and $\eta_n^2 \, id -\phi$ are increasing
functions. Then for $i=1,n$,
\begin{equation}\label{rfirn}
    \hh^s \left( \calr_\phi^i \right) =\left\{%
\begin{array}{ll}
    0\, , & \mbox{ if } \sum_{k=1}^\infty k^{n} e^{-\frac{s\eta_n }{2}\phi \iv \left( 2\eta_n \ln k \right)} < \infty \, , \\
    \infty \, , & \hbox{ if } \sum_{k=1}^\infty k^{n}e^{-\frac{s\eta_n }{2}\phi \iv \left( 2\eta_n \ln k \right)} = \infty \, . \\
\end{array}%
\right.
\end{equation}
\end{theorem}

\proof Consider $\mathcal{S}^e_\psi (\R^n )$ the set defined as
$\mathcal{S}_\psi (\R^n )$, but with the max-norm replaced by the
Euclidean norm. We have that $\mathcal{S}_{\frac{1}{\sqrt{n}}\psi}
(\R^n ) \subset \mathcal{S}^e_\psi (\R^n )\subset \mathcal{S}_\psi
(\R^n )$. This easily implies that (\ref{jar}) holds with
$\mathcal{S}_\psi (\R^n )$ replaced by $\mathcal{S}^e_\psi (\R^n
)$.

\me

\textit{Step }1. We first consider the set $\calr_\phi^1$. We
recall that according to (\ref{uri}),
\begin{displaymath}
U_+(r_1)=\left\{ \mathbf{u}_{\bar{x}}=\left(
\begin{array}{cc}
Id_{n} & \bar{x} \\
0 & 1
\end{array}
\right) \; ;\; \bar{x}\in \R^n \right\}\; .
\end{displaymath}

We may therefore identify $\R^n$ with $U_+(r_1)$ and thus identify
$\mathcal{S}^e_\psi (\R^n )$ to a subset of $U_+(r_1)$, which we
denote by $\ss_\psi$.

\me

\noindent \textbf{I.} We prove that
\begin{equation}\label{fi1}
 \calr_\phi^1 \subset
\ss_{\psi_1},\mbox{ where }  \psi_1 (x)=x
e^{-\frac{\eta_n}{2}\phi\iv (2\eta_n \ln x)}\, .
\end{equation}

Let $\ux \in \calr_\phi^1$. Infinitely many times as $t\to \infty$
we have that
\begin{equation}\label{if1}
 -f_{\bar{r}_1}\left(\pr \left(r_1(t)\ux \right)\right)\geq t -
\phi (t)\, .
\end{equation}

By the discussion in the beginning of Section \ref{bparab}, this
is equivalent to the statement that for infinitely many $t$ as
$t\to \infty$ and infinitely many $(\bar{p},q)\in \pri^{n+1}$,
$$
f_{(\bar{p},q)}(r_1(t)\ux )\leq \phi (t) -t\, .
$$

The last inequality is equivalent to $r_1(t) \circ \ux
(\bar{p},q)\leq e^{\frac{\phi(t)-t}{\eta_n}}$, which in its turn
writes as
$$
e^{\frac{t}{\sqrt{n(n+1)}}}\| \bar{p} + q \bar{x}\|_e^2
+e^{-\frac{nt}{\sqrt{n(n+1)}}} q^2\leq
e^{\frac{\phi(t)-t}{\eta_n}}\, .
$$

The inequality $e^{-\frac{nt}{\sqrt{n(n+1)}}} q^2\leq
e^{\frac{\phi(t)-t}{\eta_n}}$ is equivalent by monotonicity of
$\phi$ with
\begin{equation}\label{t1}
t\geq \phi\iv (2\eta_n \ln q )\, .
\end{equation}

The inequality $e^{\frac{t}{\sqrt{n(n+1)}}}\| \bar{p} + q
\bar{x}\|_e^2 \leq e^{\frac{\phi(t)-t}{\eta_n}} $ then implies
that $\| \bar{p} + q \bar{x}\|_e^2 \leq e^{\frac{\phi(t)
}{\eta_n}-\eta_n t}$. This, inequality (\ref{t1}) and the fact
that the function $\eta_n^2 id -\phi $ is increasing, imply that
$$
\| \bar{p} + q \bar{x}\|_e \leq \psi_1 (q)\, .
$$

\me

\noindent \textbf{II.} We prove that
\begin{equation}\label{fi2}
\ss_{\psi_2}\subset \calr^1_\phi,\mbox{ where }\psi_2(x)= x
e^{-\frac{\eta_n}{2}\phi\iv (2\eta_n \ln (\sqrt{2}x))}\, .
\end{equation}

Let $\ux \in \ss_{\psi_2}$. Then for infinitely many
$(\bar{p},q)\in \pri^{n+1}$,
$$
\| \bar{p} + q \bar{x}\|_e \leq \psi_2 (q)\, .
$$

We take $t=\phi\iv (2\eta_n \ln (\sqrt{2}q))$ and consider
$$
e^{\frac{t}{\sqrt{n(n+1)}}}\| \bar{p} + q \bar{x}\|_e^2
+e^{-\frac{nt}{\sqrt{n(n+1)}}} q^2\, .
$$

By the choice of $t$, the second term of the sum is equal to
$\frac{1}{2}e^{\frac{\phi(t)-t}{\eta_n}}$. The first term is at
most
$$
e^{\frac{t}{\sqrt{n(n+1)}}} \psi_2
(q)^2=e^{\frac{t}{\sqrt{n(n+1)}}} q^2 e^{-\eta_n
t}=\frac{1}{2}e^{\frac{\phi(t)-t}{\eta_n}}\, .
$$

We conclude that
$$
e^{\frac{t}{\sqrt{n(n+1)}}}\| \bar{p} + q \bar{x}\|_e^2
+e^{-\frac{nt}{\sqrt{n(n+1)}}} q^2\leq
e^{\frac{\phi(t)-t}{\eta_n}}
$$

This implies inequality (\ref{if1}) for the chosen $t=t(q)$. Hence
inequality (\ref{if1}) is satisfied infinitely many times as $t\to
\infty $, consequently $\ux \in \calr_\phi$.

The double inclusion $\ss_{\psi_2}\subset \calr^1_\phi \subset
\ss_{\psi_1}$ and (\ref{jar}) with $\mathcal{S}_\psi (\R^n )$
replaced by $\mathcal{S}^e_\psi (\R^n )$ imply the conclusion.
Note that in the divergence part, what appears is the sum in
(\ref{jar}) for the function $\psi_2$. The function $\psi_2$
differs from $\psi_1$ by a factor $\sqrt{2}$ in front of the
variable, in the argument of $\ln$. Nevertheless, it is easy to
eliminate this factor from the sum with an argument as in the
proof of Remark \ref{ficonst}.

\me

\noindent \textit{Step }2. We now consider the set $\calr_\phi^n$.
By (\ref{uri}),
\begin{displaymath}
U_+(r_n)=\left\{ \mathbf{u}_{\bar{x}}=\left(
\begin{array}{cc}
1 & \bar{x}^T \\
0 & Id_{n}
\end{array}
\right) \; ;\; \bar{x}\in \R^n \right\}\; .
\end{displaymath}

We identify $\R^n$ with $U_+(r_n)$ and thereby $\mathcal{S}^e_\psi
(\R^n )$ to a subset of $U_+(r_n)$, denoted by $\ls_\psi$.

The pre-image of $Hb_a(\bar{r}_n)$ is $\bigcup_{v\in
\pri^{n+1}}Hb_{v^*}^a$. Also, for every $v=(q,\bar{p})\in
\pri^{n+1}$,
$$
f_{v}^*(r_n(t)\ux )=\eta_n \ln  \left[ \left( r_n(t)\ux
\right)^*(q,\bar{p})\right] =\eta_n \ln \left(
e^{-\frac{nt}{\sqrt{n(n+1)}}}q^2 + e^{\frac{t}{\sqrt{n(n+1)}}} \|
\bar{p}-q\bar{x} \|_e^2\right)\, .
$$

An argument almost identical to the one of Step 1 gives that
$$
\ls_{\psi_2}\subset \calr^n_\phi \subset \ls_{\psi_1}\, .
$$
This together with (\ref{jar}) implies the conclusion.\endproof

\begin{remark}\label{ficonst}
In the definition of the set $\calr^i_\phi$, one can replace
$\phi$ by $\phi_c=\phi -c$, where $c $ is a constant, and the
conclusion of Theorem \ref{thri} still holds.
\end{remark}

\proof Without loss of generality we may assume that $c>0$ (the
case $c<0$ is obtained by intertwining $\phi$ with $\phi_c$).
Theorem \ref{thri} applied to the function $\phi_c$ gives
(\ref{rfirn}) with $\phi_c$ instead of $\phi$. The sum appearing
in (\ref{rfirn}) is $\Sigma_c=\sum_{k=1}^\infty k^{n}
e^{-\frac{s\eta_n \phi \iv \left( 2\eta_n \ln k +c
\right)}{2}}\leq \Sigma_0=\sum_{k=1}^\infty k^{n}
e^{-\frac{s\eta_n \phi \iv \left( 2\eta_n \ln k  \right)}{2}}$. If
$\Sigma_0<\infty$ then $\Sigma_c<\infty$.

Suppose that $\Sigma_0=\infty$. Let $p>0$ such that
$e^{\frac{c}{2\eta_n}} \leq 2^p$. There exists $r\in \{ 0,1,\dots
, 2^p-1\}$ such that $\Sigma_0(r)=\sum_{k\in 2^p\Z +r } k^{n}
e^{-\frac{s\eta_n \phi \iv \left( 2\eta_n \ln k
\right)}{2}}=\infty$. On the other hand $$\Sigma_c\geq
\sum_{k=1}^\infty k^{n} e^{-\frac{s\eta_n \phi \iv \left( 2\eta_n
\ln \left( 2^p k+r \right) \right)}{2}}\, ,$$ and the latter sum is
$+\infty$ because $\Sigma_0(r)=+\infty $.\endproof

\begin{cor}\label{corr11}
\begin{itemize}
    \item[(i)] For any $\beta \in (0,1)$, the set $\calr_\beta^i$, $i=1,n$, has Hausdorff dimension $d=n(1-\beta
    )$ and $\hh^d (\calr_\beta^i )= \infty$.
    \item[(ii)] Both statements remain true for the set
    $$
    \mathcal{E}\calr_\beta^i = \calr_\beta^i \setminus \bigcup_{\beta' > \beta
    }\calr_{\beta'}^i =\left\{ \mathbf{u}\in
\calr_\beta^i \; ;\; \limsup_{t\to \infty }
\frac{-f_{\bar{r}_i}\left(\pr
\left(r_i(t)\mathbf{u}\right)\right)}{t}=\beta \right\}\, .
    $$
\item[(iii)] The set $\calr_0^i$ coincides with $U_+(r_i)$. The
set $\calr_1^1$ is contained in $\{ \mathbf{u} \in U_+(r_1) \; ;\;
r_1 \mathbf{u}(\infty )\in \pri^{n+1} \}$ and the set $\calr_1^n$
is contained in $\{ \mathbf{u} \in U_+(r_n) \; ;\; r_n
\mathbf{u}(\infty )\in (\pri^{n+1} )^*\}$.
\end{itemize}
\end{cor}

\proof (i) follows from Theorem \ref{thri}, (ii) follows from (i),
(iii) is obtained as Corollary \ref{corrbeta}, (c). \endproof

\subsection{Case when the ray measuring the depth has a different
slope}\label{cds}

We now try to relate the set $\mathcal{L}_\psi (\R^n )$ to sets of
geodesic rays similar to $\mathcal{R}_\phi^i,\, i=1,n $, from
(\ref{rr1}) and (\ref{rrn}), and to reformulate (\ref{linrh}) in
terms of their Hausdorff measure. It turns out that the sets
$\mathcal{R}_\phi$ to be considered in this case are a bit
different from all the sets considered before. More precisely,
what has to be considered is either the set of rays of slope
$r_1(\infty )$ and their divergence measured with respect to
$\bar{r}_n$ or the set of rays of slope $r_n(\infty )$ and their
divergence measured with respect to $\bar{r}_1$. Before defining
them, we remark that the maximal possible depth of $\pr
\left(r_i(t)\mathbf{u}\right)$ measured with respect to
$\bar{r}_j$, when $\{i,j\}=\{ 1,n \}$, is $\frac{1}{n}t+c$, where
$c$ is a constant. Such a depth occurs infinitely many times if
and only if the ray $\pr \left(r_i\mathbf{u}\right)$ is contained
in the projection of a Weyl chamber with a 1-dimensional face
asymptotic to $\bar{r}_j$ (see Corollary \ref{corr1n}, (iii)).
Hence in this case in the definitions of the sets of rays moving
away in the cusp one must put $\frac{1}{n}t-\phi(t)$ with
$\phi:[a,+\infty )\to [b,+\infty )$.

Thus we define
\begin{displaymath}
\mathcal{R}_{\phi }^{1n} = \left\{ \mathbf{u}\in U_+(r_1) \; ;\;
-f_{\bar{r}_n}\left(\pr \left(r_1(t)\mathbf{u}\right)\right)\geq
\frac{1}{n} t - \phi (t) \mbox{ infinitely many times as }t\to
\infty \right\}\, .
\end{displaymath}

Similarly we define
\begin{displaymath}
\mathcal{R}_{\phi }^{n1} = \left\{ \mathbf{u}\in U_+(r_n) \; ;\;
-f_{\bar{r}_1}\left(\pr \left(r_n(t)\mathbf{u}\right)\right)\geq
\frac{1}{n} t - \phi (t) \mbox{ infinitely many times as }t\to
\infty \right\}\, .
\end{displaymath}

Let $\cl^e_\psi (\R^n )$ be the set defined as $\cl_\psi (\R^n )$
but with the max-norm replaced by the Euclidean norm. We need to
replace in (\ref{linrh}) $\cl_\psi (\R^n )$ by $\cl_\psi^e (\R^n
)$.

\begin{lemma}\label{em}
Let $\psi$ be an approximating function such that
$\psi_1(x)=\frac{\psi (x)}{x}$ is a decreasing function. Then
(\ref{linrh}) holds with $\cl_\psi (\R^n )$ replaced by
$\cl_\psi^e (\R^n )$.
\end{lemma}

\proof The hypothesis on $\psi$ implies that
$$
\cl_{\frac{1}{\sqrt{n}}\psi (\sqrt{n}\: \cdot\: ) }(\R^n ) \subset
\cl_\psi^e (\R^n ) \subset \cl_{\sqrt{n} \psi } (\R^n )\, .
$$ This and (\ref{linrh}) yield
\begin{equation}\label{line}
    \hh^s (\mathcal{L}^e_\psi(\R^n ))=\left\{%
\begin{array}{ll}
    0\, , & \mbox{ if } \sum_{k=1}^\infty k^{n}\psi_1 (k)^{s-(n-1)} < \infty \, ,\\
    \infty \, , & \hbox{ if } \sum_{k=1}^\infty k^{n}\psi_1 (\sqrt{n} k)^{s-(n-1)} = \infty \, . \\
\end{array}%
\right.
\end{equation}

Let $p\in \N$ be such that $2^p\geq \sqrt{n}$. The hypothesis
$\sum_{k=1}^\infty k^{n}\psi_1 (k)^{s-(n-1)} = \infty$ implies
that there exists $r\in \{ 0,1,\dots ,2^p-1 \}$ such that
$\sum_{k\in 2^p\Z +r} k^{n}\psi_1 (k)^{s-(n-1)} = \infty$. We have
that $\sum_{k=1}^\infty k^{n}\psi_1 (\sqrt{n} k)^{s-(n-1)}\geq
\sum_{k=1}^\infty k^{n}\psi_1 (2^p k+r )^{s-(n-1)}$, which implies
that $$\sum_{k=1}^\infty k^{n}\psi_1 (\sqrt{n} k)^{s-(n-1)}=\infty
\, .$$

Hence (\ref{line}) can be re-written such that the sum
$\sum_{k=1}^\infty k^{n}\psi_1 (k)^{s-(n-1)}$ also appears on the
second line.\endproof

\begin{theorem}\label{thrij}
Let $\phi :[a,+\infty ) \to [b,+\infty )$ be a bijective function
such that $\phi$ and $\eta_n^2 \, id -\phi$ are increasing
functions. Then for $\{ i,j\} = \{ 1,n \}$,
\begin{displaymath}
    \hh^s \left( \calr_\phi^{ij} \right) =\left\{%
\begin{array}{ll}
    0\, , & \mbox{ if } \sum_{k=1}^\infty k^{n} e^{-(s-n+1)\frac{\eta_n }{2}\phi \iv \left( 2\eta_n \ln k \right)} < \infty \, , \\
    \infty \, , & \hbox{ if } \sum_{k=1}^\infty k^{n}e^{-(s-n+1)\frac{\eta_n }{2}\phi \iv \left( 2\eta_n \ln k \right)} = \infty \, . \\
\end{array}%
\right.
\end{displaymath}
\end{theorem}

\proof We give a proof only for $i=1,j=n$, the argument in the
other case is similar. As in Step 1 of the proof of Theorem
\ref{thri}, we can identify $\R^n$ to $U_+(r_1)$ and thus identify
$\cl_\psi^e (\R^n )$ to a subset of $U_+(r_1)$, which we denote by
$\tl_\psi$. We prove that
\begin{equation}\label{dincl}
\tl_{\psi_2}\subset \calr_\phi^{1n} \subset \tl_{\psi_1}\, ,
\end{equation} where $\psi_1$ and $\psi_2$ are the functions
defined in (\ref{fi1}) and (\ref{fi2}), respectively.

\me

\noindent \textbf{I.} We prove the second inclusion in
(\ref{dincl}). Let $\ux \in \calr_\phi^{1n}$. Then for infinitely
many $t>0$ and $v=(\bar{p},q)\in \pri^{n+1}$ we have that
$$
f_v^*\left( r_1(t)\ux \right)\leq \phi (t)-\1n t\; \Leftrightarrow
\; \eta_n \ln \left[ \left( r_1(t)\ux \right)^*(\bar{p},q)\right]
\leq \phi (t)-\1n t \, .
$$

The inequality above is equivalent to
$$
e^{-\frac{t}{\sqrt{n(n+1)}}} \| \bar{p}\|^2_e +
e^{\frac{nt}{\sqrt{n(n+1)}}} (\bar{p}\cdot \bar{x} - q)^2 \leq
e^{\frac{\phi (t)}{\eta_n}-\frac{t}{\sqrt{n(n+1)}}}\, .
$$

The inequality $e^{-\frac{t}{\sqrt{n(n+1)}}} \| \bar{p}\|^2_e \leq
e^{\frac{\phi (t)}{\eta_n}-\frac{t}{\sqrt{n(n+1)}}}$ is equivalent
to
\begin{equation}\label{tlow}
    t\geq \phi\iv\left( 2\eta_n \ln \| \bar{p} \|_e \right)\, .
\end{equation}

The inequality $e^{\frac{nt}{\sqrt{n(n+1)}}} (\bar{p}\cdot \bar{x}
- q)^2 \leq e^{\frac{\phi (t)}{\eta_n}-\frac{t}{\sqrt{n(n+1)}}}$
implies that
\begin{equation}\label{aprlow}
   (\bar{p}\cdot \bar{x}
- q)^2\leq e^{\frac{\phi (t)}{\eta_n} - \eta_n t}\, .
\end{equation}

The hypothesis that $\eta_n^2 \, id -\phi$ is an increasing
function, together with (\ref{tlow}) and (\ref{aprlow}) imply that
$|\bar{p}\cdot \bar{x} - q|\leq \psi_1(\|\bar{p}\|_e )$.

\me

\noindent \textbf{II.} We prove the first inclusion in
(\ref{dincl}). Let $\ux \in \tl_{\psi_2}$. It follows that for
infinitely many $(\bar{p},q)\in \pri^{n+1}$, we have
$|\bar{p}\cdot \bar{x} - q|\leq \psi_2(\|\bar{p}\|_e )$. Consider
$t=\phi\iv\left( 2\eta_n \ln (\sqrt{2} \| \bar{p} \|_e) \right)$
and the expression
\begin{equation}\label{expr}
e^{-\frac{t}{\sqrt{n(n+1)}}} \| \bar{p}\|^2_e +
e^{\frac{nt}{\sqrt{n(n+1)}}} (\bar{p}\cdot \bar{x} - q)^2 \, .
\end{equation}

The choice of $t$ implies that $e^{-\frac{t}{\sqrt{n(n+1)}}} \|
\bar{p}\|^2_e = \frac{1}{2}e^{\frac{\phi
(t)}{\eta_n}-\frac{t}{\sqrt{n(n+1)}}}$. We also have
$$
e^{\frac{nt}{\sqrt{n(n+1)}}} (\bar{p}\cdot \bar{x} - q)^2\leq
e^{\frac{nt}{\sqrt{n(n+1)}}}\psi_2(\|\bar{p}\|_e
)^2=\frac{1}{2}e^{\frac{\phi
(t)}{\eta_n}-\frac{t}{\sqrt{n(n+1)}}}\, .
$$

Thus the expression in (\ref{expr}) is at most $e^{\frac{\phi
(t)}{\eta_n}-\frac{t}{\sqrt{n(n+1)}}}$, which implies that
$$
f_{(\bar{p},q)}^*\left( r_1(t)\ux \right)\leq \phi (t)-\1n t \, .
$$

Since this holds for infinitely many $(\bar{p},q)\in \pri^{n+1}$,
we obtain that $\ux \in \calr_\phi^{1n}$.

\me

The double inclusion (\ref{dincl}) and Lemma \ref{em} yield the
conclusion.\endproof

\begin{remark}\label{ficonstl}
In the definition of the sets $\calr^{ij}_\phi$, the function
$\phi$ can be replaced by $\phi_c=\phi -c$, where $c $ is a
constant, and Theorem \ref{thrij} still holds.
\end{remark}

\proof Similar to the one of Remark \ref{ficonst}.\endproof

\begin{cor}\label{corr1n}
\begin{itemize}
    \item[(i)] For any $\beta \in \left( 0,\frac{1}{n} \right)$, the set $\calr_\beta^{ij}$, $\{ i,j \} = \{ 1,n \}$, has Hausdorff dimension $d=n(1-\beta
    )$ and $\hh^d (\calr_\beta^{ij} )= \infty$.
    \item[(ii)] Both statements remain true for the set
    $$
    \mathcal{E}\calr_\beta^{ij} = \calr_\beta^{ij} \setminus \bigcup_{\beta' > \beta
    }\calr_{\beta'}^{ij} =\left\{ \mathbf{u}\in
\calr_\beta^{ij} \; ;\; \limsup_{t\to \infty }
\frac{-f_{\bar{r}_j}\left(\pr
\left(r_i(t)\mathbf{u}\right)\right)}{t}=\beta \right\}\, .
    $$
\item[(iii)] The set $\calr_0^{ij}$ coincides with $U_+(r_i)$.

The set $\calr^{1n}_{\frac{1}{n}}$ is a subset of the set of
$\mathbf{u} \in U_+(r_1)$ such that the projective point $r_1
\mathbf{u}(\infty )$ is contained in one of the hyperplanes of
equation $x_i=q$, where $q\in \Z$ and $i\in \{ 1,2,\dots ,n\}$.

 The
set $\calr^{n1}_{\frac{1}{n}}$ is a subset of the set of
$\mathbf{u} \in U_+(r_n)$ such that the hyperplane $r_n
\mathbf{u}(\infty )$ contains one of the vectors $e_i+qe_{n+1}$,
where $q\in \Z$ and $i\in \{ 1,2,\dots ,n\}$.
\end{itemize}
\end{cor}

\proof (i) follows from Theorem \ref{thrij} and (ii) follows from
(i).

The first statement in (iii) is again a consequence of the
logarithm law \cite{KM}.

We prove the second statement. The proof of the third is similar.

Let $\ux \in \calr^{1n}_{\frac{1}{n}}$. As in the proof of Theorem
\ref{thrij}, I, we deduce that there exist infinitely many
$(\bar{p},q)\in \pri^{n+1}$ and $t>0$ such that
$$
e^{-\frac{t}{\sqrt{n(n+1)}}} \| \bar{p}\|^2_e +
e^{\frac{nt}{\sqrt{n(n+1)}}} (\bar{p}\cdot \bar{x} - q)^2 \leq
e^{-\frac{t}{\sqrt{n(n+1)}}}\, .
$$

It follows that $\| \bar{p}\|^2_e \leq 1$ and that $(\bar{p}\cdot
\bar{x} - q)^2\leq e^{-\eta_n t}$. The first inequality implies
that $\bar{p}=e_i\in \R^n$ for some $i\in \{ 1,2,\dots ,n\}$. We
may suppose that for infinitely many $t>0$ it is the same $i$. The
second inequality gives $(x_i-q)^2\leq e^{-\eta_n t}$, for
infinitely many $t$, as $t\to \infty$. There are finitely many
possibilities for $q$, so again we may suppose that in the
previous sequence of inequalities $q$ is fixed. Then as $t\to
\infty$, this gives $x_i=q$.\endproof

\end{document}